\documentclass[12pt]{article}

\usepackage{vmargin,epsfig}
\usepackage[francais]{babel}
\usepackage[T1]{fontenc}
\usepackage{amsmath,amsthm,amssymb}

\newtheorem{deftn}{Définition}[subsection]

\newtheorem{theo}[deftn]{Théorème}
\newtheorem{prp}[deftn]{Propriété}
\newtheorem{prop}[deftn]{Proposition}

\newtheorem{lemme}[deftn]{Lemme}

\newtheorem{cor}[deftn]{Corollaire}

\newenvironment{preuve}{
  \noindent \textbf{Démonstration.}}{\hfill $\square$ }

\setmarginsrb{2.3cm}{2.5cm}{2.3cm}{2.5cm}{0cm}{0cm}{2cm}{1.5cm}

\def\leq{\leqslant}
\def\geq{\geqslant}

\def\N{\mathbb{N}}
\def\Z{\mathbb{Z}}
\def\Q{\mathbb{Q}}

\def\C{\mathbb{C}}
\def\F{\mathbb{F}}
\def\O{\mathcal{O}}
\def\A{\mathcal{A}}
\def\B{\mathcal{B}}
\def\m{\mathfrak{m}}

\def\pa#1{\left(#1\right)}
\def\cro#1{\left[ #1 \right]}
\def\ccro#1{\left[ \! \left[ #1 \right] \! \right]}
\def\acco#1{\left\{ #1 \right\}}
\def\brac#1{\left< #1 \right>}

\def\epsilon{\varepsilon}
\def\ncr#1#2{\text{C}_{#1}^{#2}}
\def\hom{{\text{Hom}}}
\def\ext{{\text{Ext}}}
\def\End{{\text{End}}}
\def\ext{{\text{Ext}}}

\def\ker{\text{ker}\,}

\def\im{\text{im}\,}

\def\Fil{\text{Fil}\,}
\def\gal{\text{Gal}}
\def\card{\text{Card}\,}
\def\rg{\text{rg}\,}
\def\tr{\text{Tr}}

\def\M{\mathcal{M}}
\def\Mr{\underline{\mathcal{M}}^r}
\def\Mrtilde{\underline{\widetilde {\mathcal{M}}}^r}
\def\MFrtilde{\underline{\widetilde {\text{MF}}}^r}

\def\MSK0{\underline{\mathcal{M}}^{S_{K_0}}}

\def\pd{\text{PD}}
\def\pr{\text{pr}}

\def\Tst{T_{\text{st}}}
\def\Acris{A_{\text{cris}}}

\def\Ast{\hat A_{\text{st}}}
\def\Bst{\hat B^+_{\text{st}}}
\def\Ass{{A_{\text{ss}}}}
\def\Astinf{\hat A_{\text{st}, \infty}}

\def\up{\underline p}
\def\upi{\underline \pi}
\def\ueps{\underline \epsilon}

\def\It{I_t}
\def\Is{I_s}

\def\rat{\mathcal{R}}
\def\nr{\text{nr}}
\def\mr{\text{mr}}
\def\et{\text{ét}}

\def\t#1{\vphantom{#1}^{\text t}\! {#1}}

\def\Ker{\ker}
\def\Im{\im}
\def\Phi{\phi}
\def\lg{\text{long}\,}

\def\xyhat#1{\smash{\hat{#1}}\vphantom{#1}}

\def\calM{\mathcal{M}}
\def\calN{\mathcal{N}}
\def\calX{\mathcal{X}}
\def\calY{\mathcal{Y}}
\def\calK{\mathcal{K}}
\def\calC{\mathcal{C}}

\include{xy}
\xyoption{all}

\title{Représentations semi-stables de torsion \\ dans le cas 
$er < p-1$}
\author{Xavier Caruso}
\date{Mai 2004}

\begin{document}

\maketitle

\renewcommand{\abstractname}{Abstract}
\begin{abstract}
Let $K$ be a local field of mixed characteristic not absolutely 
ramified. Fontaine-Laffaille theory (see \cite{fontaine-laffaille}) 
gives a description of the torsion crystalline $\Z_p$ representations of 
the absolute Galois group of $K$ ($p$ denotes the characteristic of the 
residual field). Improving the former works, Breuil introduced new 
modules 
and obtained an integer and torsion thoery for the semi-stable 
representations (see \cite{breuil-ens}).

In this paper, we follow Breuil's works and adapt them to the case where
the local field $K$ can be absolutely ramified. However, we would have a
limitation on the index of absolute ramification.
\end{abstract}

\medskip

\renewcommand{\abstractname}{Résumé}
\begin{abstract}
Soit $K$ un corps local de caractéristique mixte non absolument ramifié.
La théorie de Fontaine-Laffaille (voir \cite{fontaine-laffaille}) permet 
de décrire les $\Z_p$-représentations galoisiennes 
cristallines entières de torsion ($p$ désigne la caractéristique du 
corps résiduel).
Poursuivant les précédents travaux, Breuil a introduit de nouveaux 
modules et a
obtenu une théorie entière et de torsion pour les représentations
semi-stables (voir \cite{breuil-ens}).

Dans cet article, nous reprenons les travaux de Breuil et les adaptons
dans le cas où le corps local $K$ peut être absolument ramifié. Nous 
aurons toutefois une contrainte sur l'indice de ramification absolu.
\end{abstract}

\hfill

\tableofcontents

\section{Introduction}

Dans toute la suite de ce papier, $p$ désigne un nombre premier 
et $k$ un corps parfait de caractéristique $p$. On note $\bar k$ une 
clôture algébrique de $k$, $\F_p$ le sous-corps premier de $k$ et si $q 
= p^h$ est une puissance de $p$, $\F_q$ l'ensemble des solutions dans 
$\bar k$ de l'équation $x^q = x$.

\medskip

On désigne par 
$W$ l'anneau des vecteurs de Witt à coefficients dans $k$. On rappelle 
que comme $k$ est parfait, cet anneau est un anneau de valuation 
discrète complet de caractéristique nulle dont $p$ est une uniformisante 
et dont le corps résiduel s'identifie canoniquement à $k$. On dispose 
en outre d'une application $\sigma : W \to W$ appelée Frobenius qui 
induit par passage au quotient le Frobenius classique sur $k$, 
c'est-à-dire l'élévation à la puissance $p$.

\medskip

On appelle $K_0$ le corps des fractions de $W$, c'est un corps local de 
caractéristique mixte. On prend $K$ une extension finie totalement 
ramifiée de $K_0$. On note $e$ le degré de l'extension $K/K_0$, c'est 
l'indice de ramification absolue de $K$. On appelle $\O_K$ 
l'anneau des entiers de $K$ et on choisit $\pi$ une uniformisante de cet 
anneau. On fixe $\bar K$ une clôture algébrique de $K$, on 
note $\O_{\bar K}$ l'anneau des entiers de $\bar K$ et $G_K = 
\gal\pa{\bar K / K}$ le groupe de Galois absolu du corps $K$. On note 
$I$ le groupe d'inertie (c'est un sous-groupe de $G_K$), $\Is$ le groupe 
d'inertie sauvage et $\It = I/\Is$ le groupe d'inertie modérée. Enfin, 
on appelle $v$ la valuation sur $\bar K$ normalisée par 
$v\pa{\pi}=1$ (et donc $v\pa p = e$).

\bigskip

Une $\Z_p$-représentation (resp. $\F_p$-représentation, resp. 
$\F_q$-représentation, resp. $\Q_p$-repré\-sentation) de $G_K$ est une 
action linéaire et continue de $G_K$ sur un $\Z_p$-module (resp. un 
$\F_p$-espace vectoriel, resp. un $\F_q$-espace vectoriel, resp.
un $\Q_p$-espace vectoriel). Afin d'étudier ces représentations, 
diverses catégories ont été introduites. Nous allons nous préoccuper 
dans ce papier des catégories $\Mr$ introduites par Breuil dans 
\cite{breuil-invent}, et nous montrerons 
comment il résulte de notre étude le théorème \ref{theo:serre} 
ci-dessous. 

\medskip

Avant de l'énoncer, faisons quelques rappels (pour plus de précisions, 
voir le paragraphe 1 de \cite{serre}). Soient $h$ un entier et $q = 
p^h$. 
Notons $\hat V = \{x \in \O_{\bar K} \, / \, x^{p^h} = \pi x\}$ et $V 
\subset \O_{\bar K}/p$ la réduction modulo $p$ de $\hat V$. L'espace 
$V$ est une $\F_p$-représentation de $G_K$. De plus $V$ hérite 
naturellement d'une structure de $\F_q$-espace vectoriel de dimension 
$1$ et fournit un caractère $I \to \F_q^\star$ qui se factorise 
par $\theta : \It \to \F_q^\star$. On pose $\theta_i = \theta^{p^i}$,
ce sont les \emph{caractères fondamentaux de niveau $h$}. Toute 
$\F_p$-représentation irréductible de dimension $d$ du groupe d'inertie 
modérée s'écrit comme un produit de caractères fondamentaux de 
niveau $h$ (voir la proposition 5 du paragraphe 1 de \cite{serre}).

\begin{theo}
\label{theo:serre}
Soit $X$ un schéma propre et lisse sur $K$ à réduction semi-stable 
sur l'anneau des entiers $\O_K$. On fixe $r$ un entier. Les exposants 
qui décrivent l'action de l'inertie modérée sur la semi-simplifiée 
modulo $p$ de $H^r_\et\pa{X_{\bar K}, \Q_p}^\star$ (où $X_{\bar K}$ est 
l'extension des scalaires de $X$ à $\bar K$ et où \og $\star$ \fg\ 
signifie que l'on prend le dual) sont compris entre $0$ et $er$.
\end{theo}

Ce théorème est à rapprocher d'une conjecture formulée par Serre dans 
le paragraphe 1.13 de \cite{serre} qui prédit le même résultat pour
la représentation $H^r_\et\pa{X_{\bar K}, \Z/p\Z}^\star$. À l'heure 
actuelle, cette conjecture est connue dans le cas $r=1$ bonne réduction
(\cite{raynaud}), le cas non ramifié bonne 
réduction (\cite{fontaine-laffaille}, \cite{kato}), 
le cas non ramifié à réduction semi-stable (\cite{breuil-duke}) et 
le cas $r=1$ (\cite{breuil-annals}). Le résultat donné ici ne fait 
aucune hypothèse ni sur $e$, ni sur $r$. Remarquons toutefois qu'il est 
vide pour $er \geq p-1$.

\bigskip

Soit $r$ un entier vérifiant $er < p-1$. Nous présentons dans le 
chapitre 2, la catégorie $\Mr$ et
le foncteur $\Tst$ qui associe à tout objet de 
cette catégorie une $\Z_p$-représentation de torsion de $G_K$. Le 
chapitre 3 est consacré à l'étude de la catégorie $\Mr$. En particulier, 
toujours dans le cas $er < p-1$, on démontre qu'elle est abélienne et 
artinienne.

Nous donnons ensuite dans le chapitre 4 une description complète des 
objets simples de $\Mr$, lorsque le corps résiduel $k$ est supposé 
algébriquement clos. Plus précisément nous prouvons le théorème suivant :

\begin{theo}
\label{intro:simple}
Supposons $k$ algébriquement clos et $er < p-1$. Soit $\calM$ un objet 
simple de $\Mr$. Alors, il 
existe une base $\pa{e_1, \ldots, e_h}$ de $S$ et une suite d'entiers
$\pa{n_i}$ compris entre $0$ et $er$, périodique de période 
\emph{exactement} $h$, le tout tel que :
$$\Fil^r \calM = u^{n_1} e_1 + \ldots + u^{n_h} e_h + \Fil^p S \cdot 
\calM,$$
$\Phi_r \pa{u^{n_i} e_i} = e_{i+1}$ et $N\pa{e_i} = 0$ pour tout $i$
(considéré modulo $h$).

En outre, ces objets sont tous simples et deux à deux non isomorphes.
\end{theo}

\medskip

Par la suite, nous nous intéressons véritablement au foncteur $\Tst$.
On commence par déterminer son image sur les objets simples précédemment 
calculés. On obtient le théorème :

\begin{theo}
\label{intro:tstsimple}
Supposons $k$ algébriquement clos et $er < p-1$. Soit $\calM$ un objet 
simple de $\Mr$ comme dans le théorème \ref{intro:simple}. 
Alors la représentation galoisienne $\Tst\pa \calM$ est isomorphe à :
$$\theta_1^{m_1} \theta_2^{m_2} \ldots \theta_h^{m_h}$$
où $m_i$ est défini par $n_i + m_i = er$ et où les $\theta_i$ sont les 
caractères fondamentaux de niveau $h$.

En particulier, pour tout objet $\calM$ de $\Mr$ tué par $p$, les exposants 
qui décrivent l'action de l'inertie modérée sur la semi-simplifiée 
modulo $p$ de $\Tst\pa \calM$ sont compris entre $0$ et $er$.
\end{theo}

La conclusion des chapitres $5$ et $6$ est une réponse affirmative à
une conjecture formulée à la fin de \cite{breuil-invent}, énoncé que 
nous rappelons ici :

\begin{theo}
Supposons $er < p-1$, alors le foncteur $\Tst$ de la catégorie $\Mr$ 
dans la catégorie des représentations linéaires de $G_K$ est exact, 
pleinement fidèle, d'image essentielle stable par sous-objets et 
quotients et indépendante du choix de l'uniformisante $\pi$.
\end{theo}

\medskip

Le chapitre 7 étudie les conséquences de tout ce travail préliminaire. 
On commence par répondre à un cas particulier d'une conjecture formulée 
dans \cite{breuil-azumino} (conjecture 2.2.6) :

\begin{theo}
\label{theo:fortdivtst}
Supposons $er < p-1$. Alors le foncteur $\Tst$ réalise une 
anti-équivalence de catégories entre la catégorie des modules fortement
divisibles\footnote{Voir le paragraphe \ref{sec:fortdiv} pour une 
définition.} et la catégorie des réseaux stables par Galois dans les
$\Q_p$-représentations semi-stables de $G_K$ à poids de Hodge Tate 
compris entre $0$ et $r$.
\end{theo}

\noindent
On donne ensuite une preuve du théorème \ref{theo:serre}.

\bigskip

Ce travail a été accompli dans le cadre de ma thèse de doctorat 
en mathématique que je prépare sous la direction de Christophe 
Breuil. Je tiens à le remercier vivement ici pour les conseils, les 
explications et les réponses qu'il a toujours su me fournir, ainsi que
pour la relecture patiente des versions préliminaires de ce texte. Je
tiens à remercier également Florian Herzig pour avoir relu en profondeur 
cet article, et pour ses commentaires toujours très appropriés.

\section{Présentation des objets}

Les objets introduits dans cette partie ne sont pas nouveaux et décrits 
plus en détail dans les articles \cite{breuil-ens} et 
\cite{breuil-invent}. La première de ces références n'étudie que le cas 
$e=1$, et donc ne présente les objets que dans ce cas particulier.

\medskip

On fixe maintenant et jusqu'à la fin de cet article un entier 
$r$ positif ou nul vérifiant l'inégalité $er < p-1$. Les 
définitions que nous allons donner ont un sens pour tout entier $r < 
p-1$ mais certains théorèmes ne sont plus vérifiés lorsque $er \geq 
p-1$. 

\subsection[La catégorie $\mathcal M^r$ et ses variantes]{La catégorie 
$\Mr$ et ses variantes}
\subsubsection*{L'anneau $S$}
\label{sec:defS}

On commence par définir un anneau que l'on va munir de structures 
supplémentaires. Bien que ces structures dépendent
du corps $K$ et de l'uniformisante $\pi$ choisie, nous le notons 
simplement $S$ par la suite.

\medskip

Soit $W\cro u$ l'anneau des polynômes en une indéterminée $u$ à 
coefficients dans $W$. Soit $E\pa u$ le polynôme
minimal de l'élément $\pi$ sur $K_0$, c'est un polynôme d'Eisenstein. On 
considère l'enveloppe aux puissances divisées de $W\cro u$ 
par rapport à l'idéal principal engendré par $E\pa u$ compatibles aux 
puissances divisées canoniques sur $p W\cro u$. On rappelle que 
cela signifie que l'on ajoute formellement à l'anneau $W\cro u$ 
les éléments $\frac {\pa{E\pa u}^i}{i!}$. En tant qu'anneau, $S$ est
le complété $p$-adique de cette enveloppe aux puissances divisées. De 
façon plus terre à terre, $S$ est la sous-$W$-algèbre de $K_0\ccro u$ 
suivante :
$$S = \acco{\sum_{i=0}^{\infty} w_i \frac{\pa{E\pa u}^i}{i!}, \, w_i 
\in W\cro u, \, \lim_{i\to \infty} w_i = 0}$$
ou encore :
$$S = \acco{\sum_{i=0}^{\infty} w_i \frac{u^i}{q\pa i!}, \, w_i 
\in W, \, \lim_{i\to \infty} w_i = 0}$$
où $q\pa i$ désigne le reste de la division euclidienne de $i$ par $e$, 
$e$ étant l'indice de ramification absolue de corps $K$, également le 
degré du polynôme $E\pa u$.

\medskip

On prolonge le Frobenius à l'anneau $S$ en définissant 
l'application $\phi$ par :
$$\phi \pa{\sum_{i=0}^{\infty} w_i \frac{u^i}{q\pa i!}} = 
\sum_{i=0}^{\infty} \sigma \pa{w_i} \frac{u^{pi}}{q\pa i!}.$$
Il s'agit d'une application $\sigma$-semi-linéaire.

On munit en outre $S$ de l'application $W$-linéaire $N$ définie par :
$$N \pa{\sum_{i=0}^{\infty} w_i \frac{u^i}{q\pa i!}} = -
\sum_{i=1}^{\infty} i w_i \frac{u^i}{q\pa i!}.$$
Il s'agit d'une dérivation au sens classique
mais pas de la dérivation classique par rapport à $u$, le 
degré du polynôme n'étant pas abaissé.

On munit finalement $S$ d'une filtration : pour tout entier positif ou nul 
$n$, on définit $\Fil^n S$ comme le complété $p$-adique de l'idéal 
engendré par les éléments $\frac {\pa{E \pa u}^i}{i!}$ pour $i \geq n$.
On a donc :
$$\Fil^n S = \acco{\sum_{i=n}^{\infty} w_i \frac{\pa{E\pa u}^i}{i!}, 
\quad w_i \in W\cro u, \, \lim_{i\to \infty} w_i = 0}.$$
Il est évident que $\Fil^0 S = S$, que $\Fil^n S \subset \Fil^{n-1} S$ 
et que $\bigcap_{n \in \N} \Fil^n S = 0$. On vérifie de plus certaines 
compatibilités avec les opérateurs définis précédemment : $N\pa{\Fil^n 
S} \subset \Fil^{n-1} S$ et, pour $0 \leq n \leq p-1$, 
$\phi\pa{\Fil^n S} \subset p^n S$. Ainsi, si $0 \leq n \leq p-1$, 
on pose $\phi_n = \frac \phi {p^n} : \Fil^n S \to S$. L'élément  
$\phi_1\pa{E\pa u}$ est une unité de $S$, on le notera $c$ par la suite.

\subsubsection*{Définition des catégories}

On rappelle que $r$ est un entier fixé vérifiant $er < p-1$. 
Un objet de la catégorie $\Mr$ est la donnée :
\begin{enumerate}
\item d'un $S$-module $\calM$ isomorphe à une somme directe (finie) de 
$S/p^n S$ pour des entiers $n$ convenables ;
\item d'un sous-module $\Fil^r \calM$ de $\calM$ contenant $\Fil^r S \cdot \calM$ ;
\item d'une flèche $\phi$-semi-linéaire $\phi_r : \Fil^r \calM \to \calM$ 
vérifiant la condition :
$$\phi_r\pa {sx} = \frac 1 {c^r} \phi_r \pa{s} \phi_r\pa{\pa{E\pa u}^r 
x}$$
pour tout élément $s \in \Fil^r S$ et tout élément $x \in \calM$ et 
telle que $\im \phi_r$ engendre $\calM$ en tant que $S$-module ;
\item d'une application $W$-linéaire $N : \M \to \M$ telle que :
\begin{itemize}
\item pour tout $s \in S$ et tout $x \in \calM$, $N\pa{sx} = N\pa s x + s 
N\pa x$
\item $E\pa u N\pa{\Fil^r \calM} \subset \Fil^r \calM$
\item le diagramme suivant commute :
$$\xymatrix @C=50pt {
\Fil^r \calM \ar[r]^{\phi_r} \ar[d]_{E\pa u N} & \calM \ar[d]^{cN} \\
\Fil^r \calM \ar[r]^{\phi_r} & \calM}$$
\end{itemize}
\end{enumerate}

Une flèche entre deux objets $\calM$ et $\calM'$ de cette catégorie est 
un morphisme $S$-linéaire de $\calM$ dans $\calM'$ respectant la filtration et 
commutant aux applications $\phi_r$ et $N$.

\bigskip

On peut définir également la catégorie $\Mr_0$. Il s'agit de la 
même chose sauf que l'on ne fait pas cas de l'application $N$, les objets 
sont donc la donnée des trois premiers points exposés précédemment.

\subsection{Le foncteur vers les représentations galoisiennes}
\subsubsection*{L'anneau $\Acris$}

Soit $R$ l'anneau limite projective du diagramme :
$$\O_{\bar K} / p \O_{\bar K} \leftarrow
\O_{\bar K} / p \O_{\bar K} \leftarrow
\ldots \leftarrow
\O_{\bar K} / p \O_{\bar K} \leftarrow
\ldots \leftarrow$$
les applications de transition étant à chaque fois l'élévation à la 
puissance $p$. Un élément de $R$ est une suite $\pa{u^{\pa k}}_{k 
\geq 1}$ d'éléments de $\O_{\bar K} / p \O_{\bar 
K}$ telle que pour tout entier $k$, $\pa{u^{\pa{k+1}}}^p = u^{\pa k}$.

\medskip

On considère $W\pa R$ l'anneau des vecteurs de Witt à coefficients dans 
$R$ et l'application suivante :
$$\begin{array}{rcl}
\hat\theta : \quad W\pa R & \to & \O_{\C_p} \\
\pa{a_0, a_1, \ldots, a_n, \ldots} & \mapsto & \sum_{n \geq 0} p^n 
\hat x_n^{\pa n}
\end{array} $$
où $\C_p$ désigne le complété $p$-adique de $\bar K$ et où $\hat 
x_n^{\pa n}$ est la limite quand $m$ tend vers l'infini d'une suite 
$(\hat a_n^{\pa{n+m}})
^{p^m}$, $\hat a_i^{\pa j} \in \O_{\bar K}$ désignant un relevé 
quelconque de $a_i^{\pa j}$.

\medskip

On montre\footnote{Pour une preuve simple, voir le paragraphe II.2.2 de 
\cite{berger}} que le noyau de $\hat \theta$ est l'idéal principal de $W 
\pa R$ engendré par l'élément $\xi = \cro \up - 
p$, où $\cro \up$ est le représentant de Teichmüller de $\up \in R$ 
défini par $\up = \pa{p_1, \ldots, p_n, \ldots}$, les $p_n$ formant un 
système compatible de racines $p^n$-ièmes de $p$. L'anneau $\Acris$ 
s'obtient en introduisant des puissances divisées en $\xi$, et en 
complètant $p$-adiquement :
$$\Acris = \acco{\sum_{i \geq 0} a_i \frac {\xi^i}{i!} \, , \, a_i \in 
W\pa R, \, a_i \longrightarrow 0}.$$

L'anneau $\Acris$ hérite d'un Frobenius $\phi$ et d'une 
action du groupe de Galois $G_K$ définis \emph{via} leur action sur $W 
\pa R$. On munit également $\Acris$ d'une filtration décroissante 
définie de la façon suivante :
$$\Fil^n \Acris = \acco{\sum_{i \geq n} a_i \frac {\xi^i}{i!} \, , \, 
a_i \in W\pa R, \, a_i \longrightarrow 0} \subset \Acris.$$

\subsubsection*{L'anneau $\Ast$}

L'anneau $\Ast$ s'obtient en complétant $p$-adiquement la PD-algèbre 
polynomiale $\Acris\brac X$ :
$$\Ast = \acco{\sum_{i \geq 0} a_i \frac {X^i}{i!} \, , \, a_i \in 
\Acris, \, a_i \longrightarrow 0}.$$
On étend le Frobenius et l'action de Galois à $\Ast$ de la façon 
suivante. On pose $\phi\pa{X} = \pa{1+X}^p - 1$.
Soit $\pa{\pi_1, \pi_2, \ldots, \pi_n, \ldots}$ un système compatible de 
racines $p^n$-ièmes de l'uniformisante\footnote{Ainsi l'anneau $\Ast$ 
dépend \emph{a priori} du choix de ce système compatible de racines. 
Cependant, on prouve que ce n'est pas le cas.} $\pi$ et soit $g \in 
G_K$. Pour 
tout entier $n$, il existe $\epsilon_n \pa g$ une racine $p^n$-ième de 
l'unité telle que $g\pa{\pi_n} = \epsilon_n \pa g \pi_n$. La suite 
$\cro{\epsilon_n \pa g}$ définit un élément $\cro{\ueps \pa g} \in 
\Acris$. L'élément $g$ agit sur $X$ par $g\pa X = \cro{\ueps 
\pa g} X + \cro{\ueps \pa g} - 1$. La 
filtration sur $\Ast$ est obtenue en faisant le produit de convolution 
de la filtration de $\Acris$ par la filtration naturelle 
donnée par les puissances divisées en $X$ :
$$\Fil^n \Ast = \acco{\sum_{i \geq 0} a_i \frac {X^i}{i!} \, , \, 
a_i \in \Fil^{n-i} \Acris, \, a_i \longrightarrow 0} \subset \Ast$$
avec la convention $\Fil^k \Acris = \Acris$ si $k < 0$.
Pour $n \leq p-1$, on a $\phi(\Fil^n \Ast) \subset p^n \Ast$ et on pose
$\phi_n = \frac {\phi}{p^n}_{| \Fil^n \Ast}$.

\medskip

On munit finalement $\Ast$ d'un opérateur de monodromie $N$ défini comme 
l'unique dérivation continue $\Acris$-linéaire telle que $N\pa{X} = 
1+X$. 

\medskip

L'anneau $\Ast$ n'est pas sans lien avec $S$ : 
dans \cite{breuil-griffiths}, Breuil prouve que le morphisme de 
$W$-algèbres $S \to \Ast$, $u \mapsto \frac {\cro \upi} {1+X}$ ($\cro 
\upi$ désigne le représentation de Teichmüller de $\upi = \pa{\bar 
\pi_1, \ldots, \bar \pi_n, \ldots} \in R$, $\bar \pi_i$ étant la
réduction modulo $p$ de $\pi_i$) identifie $S$ avec l'ensemble
$\Ast^{G_K}$ des 
invariants de $\Ast$ sous l'action du groupe de Galois. En 
outre, ce morphisme fait de $\Ast$ un $S$-module. Toutefois, $\Ast$
ne vérifie pas les propriétés nécessaires pour être un objet de la
catégorie $\Mr$.

\subsubsection*{Le foncteur $\Tst$}

On pose $\Astinf = 
\Ast \otimes_W K_0/W$. L'action du groupe de Galois, le Frobenius, la 
filtration et la monodromie s'étendent à $\Astinf$ car $\Fil^n \Ast 
\cap p^r \Ast = p^r \Fil^n \Ast$. En outre, pour la même raison, si $r < 
p-1$, l'objet $\Astinf$ hérite de
$\phi_r$. Ce n'est toutefois pas un objet de la catégorie $\Mr$ : il 
n'est pas de longueur finie en tant que $S$-module, et 
l'image de $\phi_r$ n'engendre pas tout l'espace. Il est quand même 
légitime de 
considérer l'ensemble des morphismes d'un objet $\calM$ de $\Mr$ dans
$\Astinf$ et on définit :
$$\Tst \pa \calM = \hom (\calM, \Astinf)$$
le $\hom$ précédent signifiant que l'on prend les morphismes compatibles
au $\Fil^r$, au Frobenius et à l'opération de 
monodromie. Cet ensemble est naturellement une $\Z_p$-représentation 
galoisienne de torsion, tuée par la puissance de $p$ qui annule $\calM$.

\medskip

Notre but est 
principalement d'étudier le foncteur $\Tst$, et pour ce faire, nous 
allons
quasiment toujours procéder par dévissage en regardant dans un
premier temps les objets tués par $p$, que nous étudions dans le 
paragraphe suivant.

\subsection{Les objets tués par $p$}
\label{sec:deftilde}
\subsubsection*{Les catégories $\Mrtilde$}

L'anneau important ici est $k\cro u / u^{ep}$ qui est relié 
à $S/pS$ grâce à l'application de projection $\sigma : S/pS \to k \cro u / 
u^{ep}$ définie par $\sigma \pa u = u$ et $\sigma (\frac {u^{ei}}{i!}) = 
0$ pour $i \geq p$. Sur cet anneau, on définit une filtration par $\Fil^n 
k\cro u / u^{ep} = u^{en} k\cro u / u^{ep}$, un Frobenius $\phi$ par
$\phi\pa{\sum {w_i u^i}} = \sum w_i^p  u^{ip}$ (pour 
$w_i \in k$) et un opérateur de monodromie $N$ comme 
l'unique dérivation $k$-linéaire vérifiant $N\pa u = -u$

\medskip

On définit ensuite la catégorie $\Mrtilde$ : les objets sont les 
données des quatre points qui suivent :
\begin{enumerate}
\item un $k\cro u / u^{ep}$-module $\calM$ libre de rang fini ;
\item un sous-module $\Fil^r \calM$ de $\calM$ contenant $\Fil^r k\cro 
u / u^{ep} \cdot \calM = u^{er} \calM$ ;
\item une flèche $\phi$-semi-linéaire $\phi_r : \Fil^r \calM \to \calM$ 
telle que l'image de $\phi_r$ engendre $\calM$ en tant que $k 
\cro u / u^{ep}$-module ;
\item une application $k$-linéaire $N : \calM \to \calM$ telle que :
\begin{itemize}
\item pour tout $\lambda \in k\cro u / u^{ep}$ et tout $x \in \calM$, 
$N\pa{\lambda x} = N\pa \lambda x + \lambda N\pa x$
\item $u^e N\pa{\Fil^r \calM} \subset \Fil^r \calM$
\item le diagramme suivant commute :
$$\xymatrix @C=50pt {
\Fil^r \calM \ar[r]^{\phi_r} \ar[d]_{u^e N} & \calM \ar[d]^{c_\pi N} \\
\Fil^r \calM \ar[r]^{\phi_r} & \calM}$$
où $c_\pi$ est la réduction de $c$ dans $k\cro u / u^{ep}$.
\end{itemize}
\end{enumerate}

\medskip

On introduit également la catégorie $\Mrtilde_0$ définie comme 
$\Mrtilde$ 
sauf que l'on oublie la donnée de l'opérateur $N$.

\bigskip

On peut comparer les objets de $\Mr$ tués par $p$ et ceux de 
$\Mrtilde$. Si $\calM$ est un 
objet de $\Mr$ tué par $p$, c'est naturellement un $S/pS$-module (même 
libre de rang fini), et on peut donc considérer le produit tensoriel
$T \pa \calM = \calM \otimes_{\pa \sigma} k\cro u / u^{ep}$ qui hérite
d'une filtration, d'un Frobenius et d'une monodromie et 
dont on vérifie qu'il est dans $\Mrtilde$. Cette construction définit un 
foncteur $T$ allant de la sous-catégorie pleine de $\Mr$ formée des 
objets tués par $p$ dans la catégorie $\Mrtilde$.

\begin{prop}
Le foncteur $T$ défini précédemment est une équivalence de catégories.
\end{prop}

\begin{preuve}
Elle est en tout point similaire à celle donnée pour la 
proposition 2.2.2.1 de \cite{breuil-ens}.
\end{preuve}

\bigskip

On obtient ainsi une description plus simple des objets de $\Mr$ tués 
par $p$, les objets de $\Mrtilde$ pouvant être vus comme des 
$k$-espaces vectoriels de dimension finie.

\subsubsection*{Description du quotient $\Ast / p \Ast$}

Soit $\calM$ un objet de $\Mr$ tel que $p\calM = 0$. Alors $\Tst\pa 
\calM = \hom (\calM, \Ast/p)$. Nous allons dans un premier temps décrire 
explicitement le quotient $\Ast/p$.

\bigskip

On rappelle que l'on a défini deux éléments de $R$ qui sont $\up$ et 
$\upi$. On a le résultat suivant (voir paragraphe 3.7 de 
\cite{fontaine-lnim}) :

\begin{lemme}
Avec les notations précédentes, $\Acris / p \Acris$ s'identifie 
à l'enveloppe aux puissances divisées $R^\pd$ de $R$ par rapport à
l'idéal principal engendré par $\up$. En outre, on peut également 
identifier cet anneau à $R\cro{X_i} / \pa{\up ^p, X_i ^p} _{i \geq 1}$,
l'isomorphisme envoyant $X_i$ sur la $p^i$-ième puissance divisée 
$\frac{\cro \up^{p^i}}{(p^i)!} \in \Acris / p \Acris$.
\end{lemme}

La première projection $R \to \O_{\bar K} / p$ induit un 
isomorphisme $R/\up^p R \simeq \O_{\bar K} / p$. On déduit
du lemme précédent que $\Acris / p \Acris$ s'identifie 
canoniquement à $\O_{\bar K} / p \cro{X_i} / X_i^p$, $i$ 
décrivant l'ensemble des entiers strictement positifs. 
Finalement on voit que $\Ast / p$ s'identifie à l'anneau suivant : 
$$\pa{\O_{\bar K} \cro{X_i} \brac X} / \pa{p, X_i^p}_{i \geq 1}.$$

\medskip

On rappelle que $p_1$ est une racine $p$-ième de $p$. \emph{Via} les 
identifications précédentes, et pour $n < p$, $\Fil^n (\Ast/p)$ est le 
sous-$\O_{\bar K}/p$-module de $\Ast/p$ engendré par les $p_1^{n-i} 
\frac{X^i}{i!}$ (pour 
$i \leq n$), les $\frac {X^i}{i!}$ (pour $i > n$) et les $X_i$ (pour $i 
\geq 1$). On a $\phi_r \pa{X_i} = 0$ et $\phi_1 \pa X = \frac 
{\pa{1+X}^p - 1} p = Y$. La monodromie est l'unique dérivation 
$\pa{\Acris / p \Acris}$-linéaire et continue $N$ qui envoie $\frac 
{X^i} {i!}$ sur $\pa{1+X} \frac{X^{i-1}}{\pa{i-1}!}$.

\subsubsection*{Description du foncteur $\Tst$}

Nous cherchons à faire le transport 
\emph{via} le foncteur $T$ pour voir comment le foncteur $\Tst$ se 
réalise à travers la catégorie $\Mrtilde$. L'objet à calculer est le 
produit tensoriel $\Ast / p \otimes_{S/pS} k\cro u / u^{ep}$. Pour cela, 
on définit $\hat A = \pa{O_{\bar K}/p} \brac X$. On a un morphisme de 
$\O_{\bar K}/p$-modules :
$$\pr : \Ast/p\Ast \to \hat A$$
donné, \emph{via} la description précédente, par $\pr\pa X = X$ et $\pr 
\pa {X_i} = 0$ pour tout $i$. On vérifie que $\pr$ est 
$S/pS$-linéaire. On définit également $\Fil^r \hat A = \pr (\Fil^r 
\Ast/p\Ast)$ et on vérifie que l'on peut définir une unique application
$\phi_r : \Fil^r \hat A \to \hat A$ (resp. $N : \hat A \to \hat A$) 
vérifiant $\phi_r \circ \pr = \pr \circ \phi_r$ (resp. $N \circ \pr = 
\pr \circ N$). Notons qu'il faut faire attention lorsque l'on veut
effectuer des calculs avec $\phi_r$ : avant d'élever à la puissance
$p$, il faut toujours relever l'élément dans $\Fil^r \Ast$. Par exemple,
dans $\Acris / p \Acris$, on a $\phi_1 \pa{p_1} = X_1 - 1$ et donc si 
$x \in \O_{\bar K} / p$ est un multiple de $\pi_1^{er}$, on obtient
$\phi_r \pa x = \pa{-1}^r \frac {x^p}{p^r}$, avec un signe \emph{a 
priori} inattendu.

\medskip

De plus, on a une inclusion $S/pS$-linéaire :
$$i : k\cro u / u^{ep} \to \hat A$$
définie par $i\pa 1 = 1$. On peut former le produit :
$$\pr \cdot i : \Ast / p \Ast \otimes_{S/pS} k\cro u / u^{ep} \to \hat 
A.$$

\begin{lemme}
\label{lem:hata}
L'application précédente est un isomorphisme qui respecte les 
structures.
\end{lemme}

\begin{preuve}
La surjectivité et le respect des structures sont immédiats. Comme 
$\sigma : S/pS \to k\cro u / u^{ep}$ est surjectif, tout élément de 
$\Ast/p\Ast \otimes_{S/pS} k\cro 
u / u^{ep}$ s'écrit $x \otimes 1$ avec $x \in \Ast/p\Ast$. Pour vérifier 
l'injectivité, il suffit donc de voir que $\pa{\ker \pr} \otimes_{S/pS} 
k\cro u / u^{ep} = 0$ mais ceci résulte directement de :
$$X_i \otimes 1 = \frac{\pi_1^{ep^i}}{(p^i)!} \otimes 1 = 
\frac{u^{ep^i}}{(p^i)!} \otimes 1 = 1 \otimes 
\sigma\pa{\frac{u^{ep^i}}{(p^i)!}} = 0.$$
\end{preuve}

\bigskip

On construit une application :
$$\hom (\calM, \Astinf) \to \hom (T \pa \calM, \hat A)$$
déduite de la tensorisation par $k \cro u / u^{ep}$ au-dessus de $S/pS$
(et où tous les morphismes doivent commuter aux structures 
supplémentaires).

\begin{lemme}
\label{lem:tsttilde}
L'application précédente est un isomorphisme de $\Z_p$-modules 
galoisiens.
\end{lemme}

\begin{preuve}
Commençons par l'injectivité. Soit $\psi \in \hom (\calM, \Astinf)$ 
induisant par tensorisation l'application nulle $T\pa \calM \to \hat A$. Le 
diagramme suivant est commutatif :
$$\xymatrix @C=50pt {
\calM \ar[r]^-{\psi} \ar[d]_{x \mapsto 1 \otimes x} & \Astinf \ar[d]^\pr \\
T\pa \calM \ar[r]^-{0} & \hat A }$$
d'où $\im \psi \subset \ker \pr$. On vérifie facilement que $\phi_r \pa 
{\ker \pr} = 0$. Comme $\psi$ commute à $\phi_r$ et $\phi(\Fil^r \calM)$ 
engendre $\calM$, on en déduit $\psi = 0$. L'application $\hom (\calM, 
\Astinf) \to \hom (T \pa \calM, \hat A)$ est injective.

\medskip

Pour la surjectivité, notons $T\pa \calM _0$ l'image de $\phi_r$ sur $T\pa 
\calM$. La preuve de la proposition 2.2.2.1 de \cite{breuil-ens} fournit
l'isomorphisme :
$$\calM \simeq T\pa \calM _ 0 \otimes_{k\cro{u^p}/u^{ep}} S/pS.$$
Soit $\bar \psi : T\pa \calM \to \hat A$. D'après l'isomorphisme précédent, 
elle induit une application $S/pS$-linéaire $\calM \to \hat A  
\otimes_{k\cro{u^p}/u^{ep}} S/pS$, et ce dernier module s'envoie de 
façon naturelle dans $\Astinf$. On vérifie finalement que l'application 
composée commute à $\Fil^r$, $\phi_r$ et $N$ et relève $\bar \psi$.
\end{preuve}

\subsubsection*{Description de l'anneau $\hat A$}

\begin{lemme}
\label{lem:descy}
Soit $R$ un anneau dans lequel tous les entiers premiers à $p$ sont 
inversibles. Alors on a un isomorphisme :
$$\pa{R\cro{X'}\brac Y}/\pa{X'^p-1, p} \longrightarrow \pa{R \brac X}/p 
$$
envoyant $X'$ sur $X+1$ et $\frac{Y^i}{i!}$ sur 
$\frac 1 {i!} \pa{\frac{\pa{X+1}^p-1} p}^i$.
\end{lemme}

\begin{preuve}
D'abord, l'application précédente, disons $\psi$, est bien 
définie : on a $\pa{1+X}^p \equiv 1 + X^p \equiv 1 \pmod p$.

\medskip

Pour prouver que $\psi$ est un isomorphisme, on 
remarque que chacun des objets intervenant est un $R/p$-module libre et 
que $\psi$ est $R/p$-linéaire. Une base du module source est donnée 
par la famille $\pa{X'^i \cdot \frac{Y^j}{j!}}_{0 \leq i \leq p-1, \, j 
\geq 0}$. Le module but admet pour base la famille 
$\pa{\frac{X^n}{n!}}_{n \geq 0}$. L'image par 
$\psi$ de l'élément $X'^i \cdot \frac{Y^j}{j!}$ est :
$$\psi\pa{X'^i \cdot \frac{Y^j}{j!}} = \pa{1+X}^i \cdot 
\frac{\pa{\frac{\pa{1+X}^p-1} p}^j}{j!}.$$
Le terme dominant de cette 
dernière expression est $\frac{X^{pj+i}}{p^j j!}$ et si on note 
$v_p$ la valuation $p$-adique normalisée par $v_p\pa p = 1$, on a :
$$v_p \pa{\pa{pj+i}!} = j + v_p\pa{j!} = v_p\pa{p^j j!}$$
puisque $i < p$. Comme les entiers 
premiers à $p$ sont par hypothèse inversibles dans $R$, l'égalité 
précédente assure qu'il existe un élément inversible $\alpha \in R/p$ 
tel que $p^j j! = \alpha \pa{pj+i}!$. Ainsi la \og matrice \fg\ 
représentant l'application $\psi$ dans les bases données ci-dessus est 
triangulaire et les termes diagonaux sont tous inversibles. Cela prouve 
que $\psi$ est bijective.
\end{preuve}

\medskip

L'anneau $\O_{\bar K}$ vérifie les hypothèses du lemme que l'on vient de 
prouver ; on obtient donc le corollaire suivant qui donne une nouvelle
description relativement explicite de l'anneau $\hat A$ :

\begin{cor}
\label{cor:hata}
On a un isomorphisme :
$$\hat A \to \pa{\O _{\bar K} \cro{X'} \brac Y}/\pa{X'^p-1, p}$$
En outre l'opérateur de monodromie s'exprime simplement sur cette 
description : on a $N\pa{X'} = X'$ et $N\pa{\frac{Y^i}{i!}} = 
\frac{Y^{i-1}}{\pa{i-1}!}$.
\end{cor}

\subsubsection*{Action de Galois sur l'anneau $\hat A$.}

On va déterminer l'action de Galois 
sur les éléments $X'$ et $Y$. Pour $X'$ c'est facile puisque par 
définition on a $g\pa {X'} = \epsilon \pa g X'$ pour tout $g \in G_K$.

\medskip

Pour $Y$, on pourrait être tenté d'écrire :
$$g\pa Y = \frac {\epsilon\pa g ^p \pa{1+X}^p - 1} p = Y$$
mais on n'a pas le droit de faire ce calcul à cause de la division par 
$p$. Ce qu'il faut, c'est choisir un relevé de $Y$ dans $\Ast$, calculer 
l'action de Galois sur ce relevé et voir quel élément correspond dans 
$\hat A$.

Comme relevé, on pourrait choisir $\frac {\pa{1+X}^p - 1} p$ mais on 
choisit d'abord :
$$\log \pa{1+X} = X - \frac {X^2} 2 + \frac {X^3} 3 + \ldots + 
\pa{-1}^{i-1} \frac {X^i} i + \ldots \in \Ast$$
Soit $g \in G_K$. On a $g \log \pa{1+X} = \log \pa{g\pa{1+X}} = \log 
\pa{\cro{\ueps \pa g} \pa{1+X}} = g\pa Y = Y + \hat t\pa g$ où :
$$\hat t \pa g = \log \pa{\cro{\ueps \pa g}} =  \cro{\ueps \pa g} - 
\frac {\cro{\ueps \pa g}^2} 2 + \frac {\cro{\ueps \pa g}^3} 3 + \ldots + 
\pa{-1}^{i-1} \frac {\cro{\ueps \pa g}^i} i + \ldots \in \Acris$$

\medskip

Nous allons déterminer l'image $t \pa g$ de $\hat t\pa g$ dans $\hat A$.
Remarquons que comme $\hat t \pa g \in \Acris$, on a simplement 
$t \pa g \in \O_{\bar K} / p$. Nous allons prouver qu'il s'agit  
d'une racine $\pa{p-1}$-ième de $\pa{-p}$.

\begin{lemme}
\label{lem:tsigma0}
Avec les notations précédentes, $t\pa g$ est l'image dans $\O_{\bar 
K}/p$ de :
$$- \frac {\pa{\epsilon \pa g - 1}^p} p$$
où $\epsilon \pa g \in \O_{\bar K}$ est la racine $p$-ième de l'unité 
telle que $g\pa{\pi_1} = \epsilon \pa g \pi_1$.
\end{lemme}

\begin{preuve}
Il est plus pratique ici d'écrire les choses sous la forme suivante :
$$Y - \frac{X^p}p = \frac{X'^p - 1 - \pa{X'-1}^p} p$$
et de développer :
$$Y - \frac{X^p}p = a_1 X' + a_2 {X'}^2 + \ldots + a_{p-1} {X'}^{p-1}$$
avec $a_i = \frac {\pa{-1}^i \ncr p i} p$. En appliquant $g$, on 
obtient :
$$g Y - g\pa{\frac{X^p}p} = a_1 \cro{\epsilon \pa g} X' + a_2 
\cro{\epsilon \pa g}^2 {X'}^2 + \ldots + a_{p-1} \cro{\epsilon \pa g 
}^{p-1} {X'}^{p-1}$$
d'où dans $\hat A$ :
$$t\pa g \equiv g\pa{\frac{X^p}p} - \frac {\pa{\epsilon \pa g -1}^p} p
\pmod X$$
Comme on sait que $t\pa g \in \O_{\bar K}/p$, il suffit pour conclure de 
prouver que $g\pa{\frac{X^p}p}$ est nul modulo $X$. Mais dans $\Ast$, on 
a $g\pa{\frac{X^p}p} = \frac {\pa{\cro{\ueps\pa g} \pa{1+X} - 1}^p} p$ 
et donc modulo $X$, on obtient :
$$g\pa{\frac{X^p}p} \equiv \frac {\pa{\cro{\ueps\pa g} - 1}^p} p \pmod 
X$$
On conclut en remarquant que $\cro{\ueps\pa g} - 1 \in \ker \hat 
\theta$.
\end{preuve}

\begin{lemme}
\label{lem:tsigma}
L'élément $t \pa g$ est soit nul soit égal dans $\O_{\bar K}/p$ à 
la réduction modulo $p$ d'une racine $\pa{p-1}$-ième de $\pa{-p}$
dans $\O_{\bar K}$.
\end{lemme}

\begin{preuve}
Supposons $t\pa g \neq 0$.
Prouvons d'abord que $t\pa g^{p-1} \equiv -p \pmod {p^2}$. D'après le 
lemme \ref{lem:tsigma0}, cela revient à montrer que :
$$\pa{\epsilon \pa g - 1}^{p\pa{p-1}} \equiv -p^p \pmod {p^{p+1}}.$$
Modulo $1+X+ \ldots 
+ X^{p-1}$, le polynôme $\pa{X-1}^{p-1}$ s'écrit $a_0 + a_1 X + \ldots + 
a_{p-2} X^{p-2}$ avec $a_i = \pa{-1}^i \ncr {p-1} i - 1$. On vérifie que 
$a_i$ est un multiple de $p$ et on pose $b_i = \frac {a_i} p$. En 
élevant à la puissance $p$, on obtient :
$$\pa{X-1}^{p\pa{p-1}} \equiv p^p \pa{b_0 + b_1 X + \ldots + b_{p-2} 
X^{p-2}}^p \pmod {1+X+\ldots+X^{p-1}}$$
d'où
$$\pa{X-1}^{p\pa{p-1}} \equiv p^p \pa{b_0 + b_1 + \ldots + b_{p-2}}
\pmod {1+X+\ldots+X^{p-1}, p^{p+1}}.$$
Il ne reste plus qu'à vérifier que $b_0 + b_1 + \ldots + b_{p-2} = -1$ 
pour conclure.

\medskip

Notons $\eta_1, \ldots, \eta_{p-1} \in \O_{\bar K}$ 
les racines $\pa{p-1}$-ièmes de $\pa{-p}$. On a :
$$\pa{t\pa g-\eta_1} \ldots \pa{t\pa g-\eta_{p-1}} = 0 \pmod {p^2}$$
ou encore $v\pa{t\pa g-\eta_1} + \ldots + v \pa{t\pa g-\eta_{p-1}} 
\geq 2e$. Il existe donc $i$ tel que $v\pa {t\pa g-\eta_i} 
\geq \frac {2e} {p-1}$. De plus pour tout $i$, $v\pa{\eta_i} = \frac e 
{p-1}$ et pour $i \neq j$, $v\pa{\eta_i - \eta_j} = \frac e {p-1}$ car 
deux racines $\pa{p-1}$-ièmes de l'unité sont encore distinctes dans le 
corps résiduel. Il vient, si $j \neq i$, $v\pa{t\pa g - \eta_j} = 
v\pa{\pa{t\pa g-\eta_i} + \pa{\eta_i-\eta_j}} = \frac e {p-1}$, puis
$v\pa{t\pa g-\eta_i} \geq \pa{2 - \frac{p-2}{p-1}} e \geq e$. Cela 
conclut.
\end{preuve}

\section[Généralités sur les catégories $\mathcal M^r$ et 
$\widetilde {\mathcal M}^r$]{Généralités sur les catégories $\Mr$ et 
$\Mrtilde$}

Outre de nombreuses explicitations, cette partie a pour but de démontrer 
les deux résultats suivants. D'une part les catégories $\Mr$ définies 
précédemment ne dépendent pas du choix d'une uniformisante $\pi$. 
D'autre part, ces catégories sont abéliennes et même artiniennes.

\subsection{Indépendance du choix de l'uniformisante}
\label{sec:uniformi}
Considérons $\pi$ et $\pi'$ deux uniformisantes de $K$. Notons 
respectivement $E\pa u$ et $E'\pa u$ les polynômes minimaux de 
$\pi$ et $\pi'$.

Soit $P\pa u$ un polynôme à coefficients dans $W$ tel que $P\pa {\pi} 
= \pi'$ et $P \pa 0 = 0$. On définit 
une application $\nu : S \to S$ en posant $\nu\pa s = s \circ P$. C'est 
un morphisme d'anneaux, bijectif. Il n'est par contre compatible ni au 
Frobenius, ni à l'opérateur de monodromie, et nous allons dans un 
premier temps voir comment $\nu$ se comporte vis-à-vis de ces 
opérateurs.

\medskip

Plongeons $S$ dans $T = K_0 \ccro u$ et prolongeons les opétateurs
$\phi$ et $N$ à $T$. Ils vérifient la relation $N \phi = p 
\phi N$. De même la bijection $\nu$ s'étend en une bijection de 
$T$. Notons finalement $\m$ l'idéal maximal de $S$, c'est l'idéal 
engendré par $p$, $u$ et $\frac {u^{ei}}{i!}$ pour $i \geq 1$.

\begin{lemme}
\label{lem:serie}
Soit $t \in \m$. L'application de $T$ dans $T$ définie par :
$$x \mapsto \sum_{i=0}^\infty \frac{t^i} {i!} N^i\pa x$$
est l'unique morphisme d'anneaux qui envoie $u$ sur $u \exp\pa{-t}$.
\end{lemme}

\begin{preuve}
Puisque $t \in \m$, on n'a aucun souci de convergence 
dans $T$. En outre, comme $N\pa u = -u$, il vient $N^i\pa u = \pa{-1}^i 
u$ et donc $u$ est bien envoyé sur $u \exp\pa{-t}$.

Il reste à vérifier que l'on a bien affaire à un morphisme d'anneaux. La 
stabilité par addition est immédiate. Soient $x$ et $y$ dans $T$, 
calculons :
\begin{eqnarray*}
\sum_{i=0}^\infty \frac{t^i} {i!} N^i\pa {xy}
& = & \sum_{i=0}^\infty \frac{t^i} {i!} \sum_{k+l=i} \ncr i k N^k\pa x 
N^l\pa y \\
& = & \sum_{k,l \geq 0} \frac{t^i} {k! l!} N^k\pa x N^l\pa y \\
& = & \pa{\sum_{k=0}^\infty \frac{t^k} {k!} N^k\pa {x}} \cdot 
\pa{\sum_{l=0}^\infty \frac{t^l} {l!} N^l\pa {y}}
\end{eqnarray*}
ce qui conclut la preuve.
\end{preuve}

\begin{lemme}
\label{lem:nuphir}
Il existe un (unique) élément $t \in \m$ tel que l'application $\nu^{-1} 
\circ \phi \circ \nu : S \to S$ soit donnée par la formule :
$$x \mapsto\sum_{i=0}^\infty \frac {t^i} {i!} N^i \circ \phi \pa x.$$
\end{lemme}

\begin{preuve}
Faisons les calculs dans $T$ après avoir vérifié que si une suite 
d'éléments de $S$ admet une limite dans $S$, alors elle converge aussi
dans $T$, et vers la même limite.

Regardons d'abord le cas où $P$ s'écrit $u H \pa u$ avec $H \in 1 + \m$.
Dans ces conditions on est capable de définir $\log H \in T$. D'autre 
part, notons $u S\pa u$ l'image réciproque de $u$ par $\nu$. Notons 
$H^{\pa \phi}$ le polynôme déduit de $H$ en appliquant $\phi$ à chacun 
de ses coefficients: on a $\phi\pa H \pa u = H^{\pa \phi}\pa{u^p}$.

Dans l'anneau $T$, on a alors les égalités suivantes :
$$\nu\pa x = \sum_{i=0}^\infty \frac {\pa{-\log H\pa u}^i} {i!} N^i 
\pa x \quad \text{et} \quad \nu^{-1}\pa x = \sum_{i=0}^\infty \frac 
{\pa{-\log S\pa u}^i} {i!} N^i \pa x.$$
Un calcul donne :
$$\nu^{-1} \circ \phi \circ \nu \pa x 
= \sum_{i=0}^\infty \frac{\cro{-\log \pa{H^{\pa \phi} \pa{u^p S\pa 
u^p}}}^i} {p^i i!} \nu^{-1} \circ N^i\pa{\phi \pa x}.$$
On a d'autre part :
$$\nu^{-1} \circ N^i\pa{\phi \pa x} = \sum_{j=0}^\infty \frac {\pa{-\log 
S\pa u}^j}{j!} N^{i+j}\pa {\phi \pa x}$$
et donc en regroupant :
$$\nu^{-1} \circ \phi \circ \nu \pa x = \sum_{i,j \geq 0} 
\frac{\cro{-\log \pa{H^{\pa \phi} \pa{u^p S\pa u^p}}}} {p^i i!} \frac 
{\pa{-\log S\pa u}^j}{j!} N^{i+j}\pa {\phi \pa x}$$
ce que l'on réduit, grâce à la formule du binôme, en :
$$\nu^{-1} \circ \phi \circ \nu \pa x = \sum_{i=0}^\infty \frac 1 {i!} 
\cro {\frac{-\log H^{\pa \phi}\pa{u^p S\pa u ^p}} p - \log S\pa u} ^ i 
N^i \pa {\phi \pa x}.$$
On voit sur cette dernière écriture que l'on a trouvé un candidat 
pour $t$. Il se réécrit sous la forme plus sympatique suivante :
$$t = -\frac 1 p \log \cro{S\pa u^p H^{\pa \phi}\pa{u^p S\pa u ^p}}$$

Mais par définition de $S$ et de $H$, on a $S\pa u H\pa{u S\pa u} = 1$ 
et donc en appliquant $\phi$ et en regardant modulo $p$, on trouve $S\pa 
u^p H^{\pa \phi}\pa{u^p S\pa u ^p} \equiv 1 \pmod p$. On en déduit que 
$t \in S$ et vérifie les conditions du lemme.

\medskip

Si $P$ n'est pas de la forme précédente, on peut toujours décomposer 
$\nu : S \to S \to S$ où la première flèche $\nu_0$ est de la forme 
précédente et la seconde un morphisme d'anneaux envoyant $u$ sur 
$\cro\lambda u$, où $\cro \lambda$ est le représentant de Teichmüller 
d'un $\lambda \in k$. On vérifie que l'on a le diagramme commutatif 
suivant :
$$\xymatrix @C=50pt {
S \ar[r]_-{\nu_0} \ar@/^4mm/[rr]^-{\nu} \ar[d]_-{\nu^{-1}\phi\nu \, = \,
\nu_0^{-1}\phi\nu_0 }  & S \ar[d]^\phi \ar[r] & S \ar[d]^\phi \\
S \ar[r]_-{\nu_0} & S \ar[r] & S } $$
On est donc ramené au même problème avec $\nu_0$, déjà traité.
\end{preuve}

\begin{lemme}
\label{lem:nun}
Il existe un (unique) élément $n \in S$ tel que l'application $\nu^{-1} 
\circ N \circ \nu : S \to S$ soit donnée par la formule :
$$x \mapsto n N\pa x$$
\end{lemme}

\begin{preuve}
Rappelons que l'application $\nu$ était donnée par $x \mapsto x \circ 
P$, et que l'on peut décrire $N$ \emph{via} la formule $N\pa 
s = -us'$ où $s'$ désigne la dérivée usuelle de $s$ (par rapport à $u$).

On peut alors calculer :
$$\nu^{-1} \circ N \circ \nu \pa x = \nu^{-1}\cro{N\pa{x \circ P}} = 
\nu^{-1}\cro{-u P' \cdot \pa{x' \circ P}}$$
D'autre part, on a : 
$$\nu \pa{N\pa x} = \nu\pa {-u x'} = -P \cdot \pa{x' \circ P}$$
d'où :
$$\nu^{-1} \circ N \circ \nu \pa x = \nu^{-1}\pa{\frac{-uP'\pa u}{P\pa 
u}} N\pa x$$
et on a ainsi un candidat pour $n$. Or $\nu^{-1}\pa{P\pa u} = u$ 
par définition et par $\nu^{-1}$, $u$ s'envoie sur un multiple de $u$ : 
$n = \nu^{-1}\pa{\frac{-uP'\pa u}{P\pa u}} \in S$ et convient.
\end{preuve}

\subsubsection*{Construction du foncteur}

Notons $\Mr_\pi$ (resp. $\Mr_{\pi'}$) la catégorie obtenue en 
choisissant $\pi$ (resp. $\pi'$) comme uniformisante de $K$. On souhaite 
construire un foncteur (qui va s'avérer être une équivalence de 
catégories) entre les catégories $\Mr_\pi$ et $\Mr_{\pi'}$. Notons que 
si $r=0$, les catégories $\underline{\mathcal{M}}^0_\pi$ et 
$\underline{\mathcal{M}}^0_{\pi'}$ sont identiques. On peut supposer 
$r>0$ et donc $p>2$ (puisque $er < p-1$).

\medskip

Soit $\calM$ un objet de $\Mr$. L'application $\nu$ définie précédemment 
fait de $S$ une $S$-algèbre et on remarque que si l'on munit les anneaux 
des filtrations correspondant respectivement au choix des uniformisantes 
$\pi$ et $\pi'$, l'application $\nu$ est compatible aux filtrations.

Considérons les constantes $t$ et $n$ fournies par les lemmes 
\ref{lem:nuphir} et \ref{lem:nun} et définissons :
\begin{eqnarray*}
M' & = & S \,{}_{\pa \nu}\!\otimes M \\
\Fil^r M' & = & S \,{}_{\pa \nu}\!\otimes \Fil^r M \\
\phi'_r \pa{s \otimes x} & = & \phi\pa s \otimes \pa{\sum_{i=0}^\infty 
\frac{t^i}{i!} N^i \circ \phi_r \pa x} \\
N' \pa {s \otimes x} & = & N\pa s \otimes x + s \otimes n N\pa x
\end{eqnarray*}
les deux dernières égalités étant définies pour tout $s \in S$ et 
respectivement tout $x \in \Fil^r \calM$ et tout $x \in \calM$.

\begin{lemme}
\label{lem:phirni}
Pour tout entier $i \geq 1$, le diagramme suivant est commutatif :
$$\xymatrix @C=50pt {
\Fil^r \calM \ar[r]^-{\phi_r} \ar[d]_{E\pa u^i N^i} & \calM \ar[d]^-{c^i N^i} \\
\Fil^r \calM \ar[r]^-{\phi_r} & \calM }$$
\end{lemme}

\begin{preuve}
On prouve la propriété par récurrence. Pour $i=1$, elle est vraie par 
hypothèse. Pour l'hérédité, juxtaposons les deux diagrammes :
$$\xymatrix @C=50pt {
\Fil^r \calM \ar[r]^-{\phi_r} \ar[d]_{E\pa u^i N^i} & \calM \ar[d]^-{c^i N^i} \\
\Fil^r \calM \ar[r]^-{\phi_r} \ar[d]_{E\pa u N} & \calM \ar[d]^-{cN} \\
\Fil^r \calM \ar[r]^-{\phi_r} & \calM }$$
Le grand rectangle est commutatif puisque les deux carrés le sont. 
Soient $x \in \calM$ et $s \in S$. On a :
\begin{eqnarray*}
\pa{sN} \circ \pa{s^i N^i} \pa{\phi_r \pa x} & = & s \cro{N\pa{s^i} N^i 
\pa{\phi_r \pa x} + s^i N^{i+1} \pa{\phi_r\pa x}} \\
& = & s^i N\pa s N^i \pa x + s^{i+1} N^{i+1} \pa x.
\end{eqnarray*}
En appliquant le calcul précédent deux fois et en utilisant la 
commutativité du diagramme, on obtient, pour tout $x \in \Fil^r \calM$ :
$$\begin{array}{l}
c^i N\pa c N^i \pa {\phi_r \pa x} + c^{i+1} N^{i+1} \pa {\phi_r \pa x} \\
\hspace{1cm} = \phi_r \cro {E\pa u^i N \pa{E\pa u} N^i \pa x + E\pa 
u^{i+1} N^{i+1} \pa x} \\
\hspace{1cm} = \phi \pa{N \pa{E\pa u}} \phi_r \pa {E\pa u^i N^i \pa x} 
+ \phi_r \pa{ E\pa u^{i+1} N^{i+1} \pa x}.
\end{array}$$
On sait que $\phi \pa{N \pa{E\pa u}} = N\pa c$, ce qui permet de 
conclure en utilisant une dernière fois l'hypothèse de récurrence.
\end{preuve}

\begin{lemme}
L'application $\phi'_r$ est bien définie et est $\phi$-semi 
linéaire.
\end{lemme}

\begin{preuve}
Dans un premier temps, si $x \in \Fil^r \calM$, d'après le lemme 
\ref{lem:phirni} l'élément $\frac 1{i!} N^i \circ \phi_r \pa x$ est 
bien défini puisqu'égal à $\phi_r\pa{\frac{E\pa u^i}{i!} N^i \pa x}$. 
Remarquons que $\frac{E\pa u^i}{i!} N^i \pa x$ est toujours 
élément de $\Fil^r \calM$ : si $i < r < p$, c'est vrai car $i!$ est 
inversible et si $i \geq r$, c'est vrai par hypothèse.

D'autre part, pour $i \gg 0$, on a :
$$\phi_r \pa{\frac {E\pa u^i}{i!} N^i \pa x} = \frac 1{c^r} \phi_r \pa{ 
\frac {E\pa u^i}{i!} } \phi_r \pa{E\pa u^r N^i \pa x}$$
et le facteur $\phi_r \pa{ \frac {E\pa u^i}{i!} }$ est multiple de 
$\frac {p^{i-r}} {i!}$. Comme on a supposé $p > 2$, la valuation 
$p$-adique de ce 
dernier tend vers l'infini. Cela prouve que la suite des $\frac 
1 {i!} N^i \circ \phi_r \pa x$ converge vers $0$ et donc que la somme de 
la série est bien définie.

\medskip

Reste à voir que si $s \in S$ et $x \in \Fil^r \calM$, on a $\phi'_r \pa{1 
\otimes sx} = \phi'_r \pa{\nu\pa s \otimes x}$. Comme dans le lemme 
\ref{lem:serie}, on prouve :
$$\sum_{i=0}^\infty \frac{t^i}{i!} N^i \pa {\phi_r \pa {sx}} = \pa{ 
\sum_{i=0}^\infty \frac{t^i}{i!} N^i \pa{\phi\pa s} } \cdot \pa{
\sum_{i=0}^\infty \frac{t^i}{i!} N^i \pa {\phi_r \pa x} } .$$
Le premier facteur vaut $\nu^{-1} \circ \phi \circ \nu \pa s$ d'après 
le lemme \ref{lem:nuphir}. Cela conclut, le fait que $\phi'_r$ est 
$\phi$-semi-linéaire étant évident.
\end{preuve}

\begin{lemme}
L'application $N'$ est bien définie et vérifie la condition de Leibniz.
\end{lemme}

\begin{preuve}
Comme précédemment, il s'agit de vérifier que pour $s \in S$ et $x 
\in \calM$, on a $N' \pa{1\otimes sx} = N' \pa{\nu\pa s \otimes x}$. 
Calculons :
$$N' \pa{1\otimes sx} = 1 \otimes n N \pa {sx}
= 1 \otimes n N\pa s x + 1 \otimes n s N \pa x.$$
Or d'après le lemme \ref{lem:nun}, on a $n N\pa s = \nu^{-1} \circ N 
\circ \nu \pa s$ et donc $1 \otimes n N\pa s x = N \circ \nu \pa s 
\otimes x$. D'autre part, on a $1 \otimes n s N \pa x = \nu \pa s 
\otimes n N\pa x$. On en déduit que :
$$N' \pa{1\otimes sx} = N \circ \nu \pa s \otimes x + \nu \pa s \otimes 
n N \pa x$$
comme on voulait.
\end{preuve}

\begin{prop}
L'objet $\calM'$ muni de $\Fil^r \calM'$, de $\phi'_r$ et de $N'$ est un 
objet de la catégorie $\Mr_{\pi'}$.
\end{prop}

\begin{preuve}
La seule vérification qui pose problème est la commutativité du 
diagramme reliant $\phi'_r$ à $N'$. Par un simple calcul, on prouve dans 
un premier temps qu'il existe une constante $c'$ faisant commuter le
diagramme suivant :
$$\xymatrix @C=50pt {
\Fil^r \calM' \ar[r]^-{\phi'_r} \ar[d]_-{\nu\pa{E\pa u} N'} & \calM' \ar[d]^{c' 
N'} \\
\Fil^r \calM' \ar[r]^-{\phi'_r} & \calM' \\
}$$
Comme $\nu\pa{E\pa u}$ s'obtient à partir de $E' \pa u$ simplement par 
la multiplication par une unité de $S$, un diagramme équivalent, dans 
lequel on a remplacé $\nu\pa{E\pa u}$ par $E'\pa u$ et dans lequel la
constante $c'$ a été modifié, commute. D'autre part le calcul prouve 
que la constante $c'$ obtenue ne dépend pas de $\calM$.

\medskip

Soient $n$ un entier et $\calM = S/p^n S \cdot e_1$ muni de $\Fil^r \calM = 
\calM$, $\Phi_r \pa{e_1} = e_1$ et $N\pa{e_1} = 0$. On a $\pa{c'N'} \circ 
\phi'_r\pa{u \otimes e_1} = \phi'_r \circ \pa{E'\pa u N} \pa{u \otimes 
e_1}$, ce qui donne après calcul :
$$\cro{c' p u^p - u^p \phi\pa{E'\pa u}} \otimes e_1 = 0.$$
Ainsi $p^n$ divise $c' p u^p - u^p \phi\pa{E'\pa u}$ pour tout $n$ et 
finalement $c' = \phi_1\pa{E' \pa u}$.
\end{preuve}

\bigskip

On a ainsi défini un foncteur (la définition sur les flèches est 
évidente) $\Mr_\pi \to \Mr_{\pi'}$.

\subsubsection*{Canonicité et compatibilité}

\begin{prop}
Le foncteur défini précédemment ne dépend pas du choix de l'élément $P 
\in S$.
\end{prop}

\begin{preuve}
Avec les notations précédentes, il suffit de prouver que si $P = uH$ est 
tel que $P\pa \pi = \pi$, alors $\calM$ et $\calM'$ sont canoniquement 
isomorphes. Notons $\nu : S \to S$ le morphismes d'anneau tel que 
$\nu\pa 
u = P\pa u$. La condition implique $H \pa u - 1 \in \Fil^1 S$ et donc
l'élément $\log\pa{H\pa u}$ est bien défini dans $S$.

Si $\calM$ est un objet de $\Mr_\pi$, on peut définir l'application :
$$\begin{array}{ccl}
S {}_{\pa \nu}\!\otimes \calM & \to & S \\
s \otimes x & \mapsto & s \displaystyle \sum_{i=0}^\infty 
\frac{\pa{-\log H \pa u}^i}{i!} N^i \pa x
\end{array}$$

\noindent
Comme $\log\pa{H\pa u} \in \Fil^1 S$, l'élément $\frac{\pa{-\log
H \pa u}^i}{i!}$ est bien défini. En outre le fait que dans $T$, 
$\exp\pa{\log H \pa u} = H\pa u \in S$ prouve que la suite 
$\frac{\pa{\log H \pa u}^i}{i!}$ converge vers $0$ et finalement que 
l'application est bien définie.

Il ne reste plus qu'à voir que c'est un isomorphisme $S$-linéaire et 
compatible à toutes les structures ; c'est donc une flèche dans 
$\Mr_\pi$.
\end{preuve}

\begin{cor}
\label{cor:uniformi}
Le foncteur défini précédemment est une équivalence de catégorie.
\end{cor}

\medskip

Si, comme précédemment, $\pi$ et $\pi'$ sont deux uniformisantes de $K$, 
on peut définir $\Ast{}_\pi$ et $\Ast{}_{\pi'}$. Pour cela, rappelons 
que l'on avait besoin de choisir $\pa{\pi_1, \ldots, \pi_n, \ldots}$ 
(resp. $\pa{\pi'_1, \ldots, \pi'_n, \ldots}$) un système compatible de 
racines $p^n$-ièmes de $\pi$ (resp. de $\pi'$). On définit $\omega_n$ 
en imposant $\pi_n = \omega_n \pi'_n$, obtenant ainsi $\pa{\bar 
\omega_1, \ldots, \bar \omega_n, \ldots} \in R$ puis $\cro {\underline 
\omega} \in \Acris$ (notez que $\Acris$ ne dépend pas du choix d'une 
uniformisante).

L'unique morphisme de $\Acris$-algèbre $\Ast{}_\pi \to \Ast{}_{\pi'}$
envoyant $\pa{1+X}$ sur $\cro {\underline \omega} \pa{1+X}$ est un 
isomorphisme compatible à $\phi_r$, à $N$ et à l'action du groupe de 
Galois $G_K$.

\begin{prop}
\label{prop:tstunif}
Le diagramme suivant est commutatif :
$$\xymatrix @C=10pt {
\Mr_{\pi} \ar[dr]_{\Tst{}_{\pi}} \ar[rr] & & \Mr_{\pi'} 
\ar[ld]^{\Tst{}_{\pi'}} \\
& \text{Rep}_{\Z_p}\pa{G_K} }$$
où la flèche horizontale est le foncteur défini précédemment.
\end{prop}

\begin{preuve}
L'anneau $S$ s'identifie à la fois aux points fixes sous l'action de 
Galois de $\Ast{}_\pi$ et de $\Ast{}_{\pi'}$. Notons $\rho : S \to 
\Ast{}_\pi$ et $\rho' : S \to \Ast{}_{\pi'}$ les inclusions 
correspondantes. Il existe un unique morphisme de $\Acris$-algèbre, $\nu 
: \Ast{}_\pi \to \Ast{}_{\pi'}$ faisant commuter le diagramme suivant :
$$\xymatrix @C=50pt {
S \ar[r]^\nu \ar[d]_-{\rho} & S \ar[d]^-{\rho'} \\
\Ast{}_\pi \ar@{.>}[r]^\nu & \Ast{}_{\pi'} }$$
En effet, le diagramme impose la valeur de $\nu\pa{\frac{X^i}{i!}}$ et 
on vérifie que l'application ainsi définie convient. En outre, elle est 
$G_K$-équivariante et induit une flèche $\nu : \Astinf{}_\pi \to 
\Astinf{}_{\pi'}$ encore $G_K$-équivariante.

Soit $\calM$ un objet de $\Mr_\pi$ et $\calM'$ l'objet de $\Mr_{\pi'}$ qui lui 
est associé par le foncteur précédent. On rappelle qu'en tant que 
module, on a $\calM' = S \,{}_{\pa \nu}\!\otimes \calM$. Soit $f \in 
\Tst\pa \calM$. On 
lui associe l'application suivante :
$$\begin{array}{rcl}
\calM' & \to & \Astinf{}_{\pi'} \\
s \otimes x & \mapsto & \rho'\pa s \cdot \nu \circ f \pa x 
\end{array}$$
On vérifie qu'elle est $S$-linéaire et compatible aux structures 
définissant ainsi un élément de $\Tst\pa{\calM'}$.

On définit ainsi une application $\Tst\pa \calM \to \Tst\pa{\calM'}$. Elle est 
$\Z_p$-linéaire et bijective puisque l'on peut construire l'application 
réciproque de façon analogue. On vérifie qu'elle est compatible à 
l'action de Galois et donc qu'il s'agit d'un isomorphisme dans la
catégorie des $\Z_p$-représentations galoisiennes.
\end{preuve}

\subsection[Description des objets de $\widetilde {\mathcal 
M}^r$]{Description des objets de $\Mrtilde$}
\label{sec:mrtilde}
On considère dans ce paragraphe un objet $\calM$ de $\Mrtilde$. Il s'agit 
d'un $k \cro u / u^{ep}$-module libre de rang fini $d$ muni d'une 
filtration, d'une application $\Phi_r$ et d'une application $N$, le tout
vérifiant les propriétés données précédemment.

\subsubsection*{Bases adaptées}

On a dans un premier temps un résultat bien utile (et classique) qui 
est le suivant :

\begin{prop}
\label{prop:baseadapt}
Il existe une base $\pa{e_1, \ldots, e_d}$ de $\calM$ et des entiers 
$n_1, \ldots, n_d$ tels que :
$$\Fil^r \calM = \bigoplus_{i=1}^d u^{n_i} k\cro u/u^{ep} \cdot e_i.$$
Une telle base est par définition une \emph{base adaptée} de $\calM$.
\end{prop}

\begin{preuve}
Comme $\calM$ est supposé libre, il existe un $k\cro u$-module 
libre $\calM'$ tel que $\calM'/u^{ep}\calM' = \calM$. Autrement dit, il 
existe une flèche $f : \calM' \to \calM$ dont le noyau est exactement 
$u^{ep} \calM'$. Définissons $\Fil^r \calM' = f^{-1} \pa{\Fil^r \calM}$. 
C'est un sous-$k\cro u$-module de $\calM'$.

Puisque $k\cro u$ est un anneau principal, d'après le théorème de 
structure, il existe une suite de polynômes $P_1, \ldots, P_d$ tels que 
$P_i$ divise $P_{i+1}$ pour tout $i$ et une base $\pa{\hat e_1, \ldots, 
\hat e_d}$ de $\calM'$ telle que $\pa{P_1 \hat e_1, \ldots, P_d \hat e_d}$ 
soit une base de $\Fil^r \calM'$. D'autre part, $u^{er} \calM' \subset \Fil^r 
\calM'$ et donc tous les polynômes $P_i$ sont des diviseurs de $u^{er}$ ; 
ils sont donc de la forme $P_i = u^{n_i}$ pour certains entiers $n_i$.

En posant $e_i = f(\hat e_i)$, on a bien le résultat annoncé.
\end{preuve}

\medskip

\noindent {\it Remarque}. Les entiers $n_i$ ne dépendent pas à 
permutation près de la base considérée. En effet, la dimension en
tant que $k$-espace vectoriel du quotient $\Fil^r \calM / \pa{u^k \calM \cap 
\Fil^r \calM}$ est donnée par la somme des $k-n_i$, somme étendue à tous les 
$i$ pour lesquels $n_i \leq k$. On voit facilement que la 
connaissance de toutes ces sommes permet de déterminer les $n_i$ à 
permutation près.

\medskip

Fixons à présent $\pa{e_1, \ldots, e_d}$ une base adaptée de 
$\calM$. Nous allons essayer de décrire un peu mieux la fonction $\Phi_r$ et 
pour cela nous introduisons la définition suivante.

\medskip

\begin{deftn}
\label{def:psi}
Soit $x \in \calM \backslash u \calM$, et soit $n$ le plus 
petit entier tel que $u^n x \in \Fil^r \calM$. On pose $\varphi_r 
\pa x = \Phi_r \pa{u^n x}$.
\end{deftn}

Soit $x_i = \varphi_r \pa {e_i}$ pour $1 \leq i \leq d$. On rappelle que 
la famille des $e_i$ est la base adaptée que l'on s'est fixée 
précédemment.

\begin{prop}
\label{prop:phir}
Avec les notations précédentes, $\pa{x_1, \ldots, x_d}$ est une 
base de $\calM$.

D'autre part, si $x \in \calM \backslash u\calM$, alors 
$\varphi_r \pa x \in \calM \backslash u\calM$.
\end{prop}

\begin{preuve}
Pour le premier énoncé, il suffit de voir que si $x \in
\Fil^r \calM$, alors $\Phi_r \pa x$ s'écrit comme une combinaison linéaire
(à coefficients dans $k\cro{u^p}/u^{ep}$) des $x_i$. 
Comme $\im \phi_r$ engendre $\calM$ comme $k\cro u / u^{ep}$-module, il 
en est de même de la famille $\pa{x_1, \ldots, x_d}$. Comme elle est de 
bon cardinal, elle en est une base.

\medskip

Soit $x \in \calM \backslash u\calM$. On voit en décomposant $x$ sur la 
base des $e_i$, que $\varphi_r \pa x$ s'écrit forcément sous la forme :
$$\varphi_r \pa x = Q_1 (u^p) x_1 + \ldots + Q_d (u^p) x_d$$
où au moins l'un des polynômes $Q_i$ est de valuation nulle. Dans ce 
cas, on a directement $\varphi_r \pa x \in \calM \backslash u\calM$.
\end{preuve}

\medskip

\noindent {\it Remarque}. La deuxième partie de la proposition 
précédente permet de définir correctement les itérés de $\varphi_r$.

\subsubsection*{L'opérateur de monodromie}

Nous allons à présent étudier l'opérateur de monodromie. Pour cela, nous 
notons $\calM_0 = \im \Phi_r$. Par ce qui précède, $\calM_0$ s'identifie 
au $k\cro{u^p} / u^{ep}$-module engendré par les $x_i$. Nous avons alors 
:

\begin{prop}
\label{prop:n}
Pour tout $i$, l'opérateur $c_\pi^i N^i$ induit une application 
$k\cro{u^p} / u^{ep}$-linéaire de $\calM_0$ sur lui-même. Cette 
application est nulle si $i \geq p$.
\end{prop}

\begin{preuve}
Le lemme \ref{lem:phirni} assure que le diagramme suivant :
$$\xymatrix @C=50pt {
\Fil^r \M \ar[r]^{\phi_r} \ar[d]_{u^{ei} N^i} & \M \ar[d]^{c_\pi^i N^i} \\
\Fil^r \M \ar[r]^{\phi_r} & \M}$$
commute. Ainsi si $x$ est dans l'image de $\Phi_r$, alors il 
en est de même de $c_\pi^i N^i\pa x$, et donc que $c_\pi^i N^i$ induit une 
application de $\calM_0$ dans lui-même. D'autre part, $N\pa{u^p x} = 
N\pa{u^p} x + u^p N\pa x = -p u^p x + u^p N\pa x = u^p N\pa x$, ce qui 
prouve bien la linéarité annoncée.

Pour $i \geq p$, l'application $u^{ei} N^i$ est nulle et donc il en est
de même de $c_\pi^i N^i$.
\end{preuve}

\begin{cor}
\label{cor:kern}
Il existe un élément $x \in \calM_0$ non divisible par $u$ tel que $N \pa x 
= 0$.
\end{cor}

\begin{preuve}
Du fait que $c_\pi^p N^p = 0$, il existe $x' \in \calM_0$, $x' \neq 0$, 
tel que $N\pa{x'} = 0$.
Écrivons $x' = u^{pk} x''$ où $x''$ est un élément de $\calM_0$ non 
divisible par $u$ et où $k < e$. On a alors $N\pa{u^{pk}x''} = u^{pk} 
N\pa{x''} = 0$ et donc $N\pa{x''}$ est un multiple de $u^p$.
Notons $n$ le plus petit entier tel que $u^n x'' \in \Fil^r \calM$ de telle 
sorte que l'on ait $\varphi_r \pa {x''} = \Phi_r \pa{u^n x''} = x$. On a :
$$\Phi_r \pa{u^e N\pa{u^n x''}} = c_\pi N \pa{\Phi_r \pa{u^n x''}} = 
c_\pi N \pa x.$$
Mais $N\pa{u^n x''} = -n u^n x'' + u^n N\pa{x''}$. Le premier terme de 
cette somme est dans $\Fil^r \calM$ puisque $u^n x''$ y est. Le second y 
est également puisque $N\pa{x''}$ est un multiple de $u^p$. On en déduit 
que $N\pa{u^n x''} \in \Fil^r \calM$ et donc que $\Phi_r \pa{u^e N\pa{u^n 
x''}} = 0$. Ainsi $N\pa x = 0$.
D'autre part, on a $x = \varphi_r \pa{x''}$ et donc d'après la 
proposition \ref{prop:phir}, $x$ n'est pas divisible par $u$. Ceci 
conclut la preuve du corollaire.
\end{preuve}

\subsubsection*{Description matricielle}

Le but de ce paragraphe est d'écrire sous forme matricielle les 
applications $\Phi_r$ et $N$, explicitations que nous utiliserons dans 
la suite. On fixe 
$\calM$ un objet de $\Mrtilde$ et $\pa{e_1, \ldots, e_d}$ une base 
adaptée de $\calM$, les entiers correspondants étant $n_1, \ldots, n_d$.

\medskip

On note $\Delta$ la matrice diagonale suivante :
$$\Delta = \pa{ \raisebox{0.5\depth}{\xymatrix @R=5pt @C=5pt {u^{n_1} 
\ar@{.}[dr] & \\ & u^{n_d}}} } $$

\begin{deftn}
La \emph{matrice} de $\Phi_r$ dans la base adaptée $\pa{e_1, \ldots, 
e_d}$ est la matrice $G$ définie par l'égalité suivante :
$$\pa{ \begin{array}{c} \Phi_r\pa{u^{n_1} e_1} \\ \vdots \\ 
\Phi_r\pa{u^{n_d} e_d} \end{array} } = \Phi_r \pa{ \Delta \cdot \pa{ 
\begin{array}{c} e_1 \\ \vdots \\ e_d \end{array} } } 
= \t G \cdot \pa{\begin{array}{c} e_1 \\ \vdots \\ e_d \end{array}}$$
\end{deftn}

\noindent {\it Remarque}. Cette définition n'a un sens que si la base 
$\pa{e_1, \ldots, e_d}$ est adaptée. De plus, la présence de la 
transposée sert à rester fidèle à la définition classique 
de la matrice d'une application linéaire.

En gardant les notations du paragraphe précédent, on voit que $G$ est 
simplement la matrice de passage de la base $\pa{e_1, \ldots, e_d}$ à la 
base $\pa{x_1, \ldots, x_d}$. En tant que telle, il s'agit d'une matrice 
inversible.

\medskip

\begin{deftn}
Soit $\pa{a_1, \ldots, a_d}$ une base de $\calM$. La \emph{matrice} de $N$ 
dans la base $\pa{a_1, \ldots, a_d}$ est la matrice $H$ définie par 
l'égalité suivante :
$$\pa{ \begin{array}{c} N\pa{a_1} \\ \vdots \\ N\pa{a_d} \end{array} } = 
N \pa{ \begin{array}{c} a_1 \\ \vdots \\ a_d \end{array} } = \t H \cdot 
\pa{\begin{array}{c} a_1 \\ \vdots \\ a_d \end{array}}$$
\end{deftn}

\medskip

On a une formule de changement de base :

\begin{prop}
\label{prop:basen}
Soient $\A = \pa{a_1, \ldots, a_d}$ et $\B = \pa{b_1, \ldots, b_d}$ deux 
bases de $\calM$, et soit $P$ la matrice de passage de $\A$ à $\B$. On note 
$H_\A$ (resp. $H_\B$) la matrice de $N$ dans la base $\A$ (resp. dans la
base $\B$). On a alors la relation :
$$H_\B = P^{-1} H_\A P + P^{-1} N\pa P$$
\end{prop}

\begin{preuve}
Il s'agit d'un simple calcul. On écrit :
\begin{eqnarray*}
\pa{\begin{array}{c} N\pa{b_1} \\ \vdots \\ N\pa{b_d} \end{array}} & = & 
N \pa { \begin{array}{c} b_1 \\ \vdots \\ b_d \end{array} } \, = \,
N \pa { \t P \cdot \pa{ \begin{array}{c} a_1 \\ \vdots \\ a_d 
\end{array} }} \\
& = & N\pa{\t P} \cdot \pa{\begin{array}{c} a_1 \\ \vdots \\ a_d
\end{array}} + \t P \cdot N \pa { \begin{array}{c} a_1 \\ \vdots 
\\ a_d \end{array} } \\
& = & N\pa{\t P} \t P^{-1} \cdot \pa{\begin{array}{c} b_1 \\ \vdots \\ 
b_d \end{array}} + \t P H_\A \t P^{-1} \pa{\begin{array}{c} b_1 \\ 
\vdots \\ b_d \end{array}}
\end{eqnarray*}
ce qui donne $\t {H_\B} = N\pa{\t P} \t P^{-1} +  \t P H_\A \t P^{-1}$
puis le résultat annoncé en prenant la transposée.
\end{preuve}

\medskip

\noindent {\it Remarque}. Un simple calcul prouve que si $A$ et 
$B$ sont des matrices à coefficients dans $k\cro u / u^{ep}$, 
alors $N\pa{AB} = N\pa A B + A N\pa B$. Ceci a pour conséquence 
l'égalité $P^{-1} N\pa P = - N\pa{P^{-1}} P$ et prouve la
cohérence de la formule lorsque l'on passe d'une base $\A$ à une base 
$\B$ puis que l'on revient à $\A$.

\begin{prop}
\label{prop:nnilp}
Si les $n_i$ sont rangés par ordre croissant, la matrice de $c_\pi N$ dans 
la base $\pa{x_1, \ldots, x_d}$ (où on rappelle que $x_i = \Phi_r 
\pa{u^{n_i} e_i}$) est à coefficients dans $k\cro{u^p}/u^{ep}$ et 
triangulaire inférieure avec des $0$ sur la diagonale.
\end{prop}

\begin{preuve}
La preuve résulte du diagramme commutatif suivant :
$$\xymatrix @C=50pt {
\Fil^r \M \ar[r]^{\phi_r} \ar[d]_{u^e N} & \M \ar[d]^{c_\pi N} \\
\Fil^r \M \ar[r]^{\phi_r} & \M}$$
En effet, fixons un entier $i$ et partons de l'élément $u^{n_i} e_i$. 
Par $\Phi_r$, il s'envoie sur $x_i$ par définition. Puis par $c_\pi N$, il 
s'envoie sur $c_\pi N\pa{x_i}$. Par l'autre chemin, on a d'abord :
$$u^e N\pa{u^{n_i} e_i} = -n_i u^{e+n_i} e_i + u^{e+n_i} N\pa{e_i}$$
Comme $u^{n_i} e_i \in \Fil^r \calM$, le premier terme de la somme 
précédente s'envoie sur $0$ par $\Phi_r$ et $u^{e+n_i} N\pa{e_i} 
\in \Fil^r \calM$.
Pour le second terme, décomposons $N\pa{e_i} = \sum_{j=1}^d a_j 
e_j$ où $a_j \in k\cro{u}/u^{ep}$ est tel que
$u^{e+n_i} N\pa{e_i} \in \Fil^r \calM$ (\emph{i.e.} $u^{e+n_i} a_j \in 
u^{n_j} k\cro u / u^{ep}$). On a alors :
$$\Phi_r\pa { u^e N\pa{u^{n_i} e_i} } = \sum_{j=1}^d \Phi\pa{
u^{e+n_i-n_j} a_j} x_j$$
et les $\Phi\pa{u^{e+n_i-n_j} a_j}$ sont les coefficients de la $j$-ième 
colonne de la matrice de $c_\pi N$. Ils sont donc déjà tous bien dans 
$k\cro{u^p}/u^{ep}$. 

De plus, si $j \leq i$, on a par hypothèse $n_j \leq n_i$ et donc 
$e+n_i-n_j \geq e$. Ainsi $\Phi\pa{u^{e+n_i-n_j} a_j} = 0$ et on a bien 
démontré le résultat annoncé.
\end{preuve}

\medskip

\noindent {\it Remarque}. Cette dernière proposition redémontre en 
particulier, en donnant un résultat plus précis, la proposition 
\ref{prop:n} et le corollaire qui s'ensuit.

\subsection[La catégorie $\widetilde{\text{MF}}^r$]{La catégorie 
$\MFrtilde$}
\label{sec:mfrtilde}
Dans cette partie, nous introduisons des sous-catégories pleines 
$\MFrtilde$ de $\Mrtilde$ qui correspondent aux catégories de 
Fontaine-Laffaille (voir \cite{fontaine-laffaille}) tuées par $p$ pour 
$e=1$.

\bigskip

Commençons par donner une proposition qui caractérise les objets de 
cette sous-catégorie.
Soit $\calM$ un objet de $\Mrtilde_0$. Notons $\calM_0 = \Im \Phi_r$ et plus
généralement $\calM_i = u^i \calM_0$ pour un entier $i$ compris entre $0$ et
$p-1$. Les $\calM_i$ sont des $k\cro{u^p}/u^{ep}$-modules libres et :
$$\calM = \bigoplus_{i=0}^{p-1} \calM_i.$$

\begin{prop}
\label{prop:mfrtilde}
Avec les notations précédentes, les propriétés suivantes sont 
équivalentes :
\begin{itemize}
\item[i)] $\displaystyle \Fil^r \calM = \bigoplus_{i=0}^{p-1} \Fil^r 
\calM \cap \calM_i$ ;
\item[ii)] il existe une base adaptée de $\calM$ formée d'éléments de 
$\calM_0$ ;
\item[iii)] On peut munir $\calM$ d'un opérateur de monodromie nul sur
$\calM_0$ et faisant de $\calM$ un objet de $\Mrtilde$.
\end{itemize}
\end{prop}

\begin{preuve}
La propriété ii) implique de façon presque immédiate les deux autres.
Nous allons montrer que iii) implique i) puis que i) implique ii).

\medskip

Supposons iii). 
Prouvons dans un premier temps que cela implique que $N\pa{\Fil^r \calM}
\subset \Fil^r \calM$. Soit $x \in \Fil^r \calM$. D'après
le diagramme suivant :
$$\xymatrix @C=50pt {
\Fil^r \M \ar[r]^{\phi_r} \ar[d]_{u^e N} & \M \ar[d]^{c_\pi N} \\
\Fil^r \M \ar[r]^{\phi_r} & \M}$$
on a $\Phi_r \pa {u^e N\pa x} = 0$. La proposition \ref{prop:phir} 
implique facilement que $\ker \Phi_r = u^e \Fil^r \calM$ et donc il existe 
$y \in \Fil^r \calM$ tel que $u^e N\pa x = u^e y$. La différence
$N\pa x - y$ est tuée par $u^e$. Elle s'écrit $u^{e\pa{p-1}} z$ pour un 
certain $z \in \calM$. Comme $e\pa{p-1} \geq er$, $u^{e\pa{p-1}} z 
\in \Fil^r \calM$ d'où $N\pa x \in \Fil^r \calM$. Ceci 
prouve la propriété annoncée.

Soit $y \in \Fil^r \calM$. On cherche à construire des $y_i \in \Fil^r 
\calM \cap \calM_i$ tels que $y = y_0 + \ldots + y_{p-1}$.
On peut déjà écrire une égalité de ce type avec $y_i
\in \calM_i$. Appliquons l'opérateur $N$ à cette égalité en remarquant
que puisque $N$ est supposé nul sur $\calM_0$, on a $N\pa{y_i} = -i y_i$. On 
obtient successivement :
\begin{eqnarray*}
y & = & x_0 + y_1 + \ldots + y_{p-1} \\
N\pa y & = & -y_1 + \ldots -\pa{p-1} y_{p-1} \\
N^2\pa y & = & y_1 + \ldots + \pa{p-1}^2 y_{p-1} \\
& \vdots \\
N^{p-1} \pa y & = & y_1 + \ldots + \pa{p-1}^{p-1} y_{p-1}.
\end{eqnarray*}
Les coefficients qui apparaissent forment une matrice de Vandermonde
inversible. Ainsi on peut exprimer les $y_i$ comme combinaisons
linéaires à coefficients dans $\F_p$ des $N^j\pa y$. Par ce qui précède, 
cela entraîne $y_i \in \Fil^r \calM$ et donc bien la propriété voulue.

\medskip

Supposons i). Fixons $\pa{e_1, \ldots, e_d}$ une base
de $\calM_0$ comme $k\cro{u^p}/u^{ep}$-module, et notons $\calM'_0$ le
sous-$k$-espace vectoriel de $\calM_0$ engendré par les $e_i$. Notons
également $\calM'_i = u^i \calM'_0$. Pour tout entier $i$, on a un isomorphisme
$f_i : \calM'_0 \to \calM'_i$ qui est la multiplication par $u^i$. Notons
$F'_i = f_i^{-1}\pa{\Fil^r \calM \cap \calM'_i}$. On obtient 
une filtration croissante par des sous-$k$-espaces vectoriels. Il 
suffit alors pour répondre à la question de considérer une base
$\pa{x_1, \ldots, x_d}$ de $\calM'_0$ compatible à cette filtration.
\end{preuve}

\begin{deftn}
On note $\MFrtilde_0$ la sous-catégorie pleine de $\Mrtilde_0$ formée 
des objets satisfaisant les propriétés de la proposition précédente.
On note $\MFrtilde$ la sous-catégorie pleine de $\Mrtilde$ formée des 
objets dont l'image dans $\Mrtilde_0$ par le foncteur d'oubli
satisfait les propriétés de la proposition précédente.
\end{deftn}

Les lettres $\text{MF}$ font référence à \og modules filtrés \fg\ car
l'on peut donner une nouvelle interprétation de ces objets \emph{via}
des modules filtrés. Avant cela, faisons quelques remarques générales :

\begin{prop}
\label{prop:catmfr}
La catégorie $\MFrtilde$ est une sous-catégorie abélienne de $\Mrtilde$ 
stable par sous-objets et par quotients. De plus, elle est égale à
$\Mrtilde$ si et seulement si $e=r=1$.
\end{prop}

\begin{preuve}
Nous ne savons pas encore à ce stade que $\Mrtilde$ est une catégorie 
abélienne. Nous allons l'admettre momentanément pour prouver la première
partie de la proposition. Il suffit de prouver la stabilité par 
sommes directes, sous-objets et quotients, un noyau étant un 
sous-objet et un conoyau un quotient. Tout cela est immédiat avec la 
caractérisation iii).

\medskip

Traitons le cas $e=r=1$. Soit $\calM$ un objet de $\Mrtilde_0$ et soit $x 
\in \Fil^r \calM$. On peut écrire $x = x_0 + \ldots + x_{p-1}$ avec
$x_i \in \calM_i$. Par hypothèse $u \calM \subset \Fil^r \calM$, et donc
$x_i \in \Fil^r \calM$ pour $i \geq 1$, puis $x_0 \in \Fil^r \calM$. On 
a ainsi vérifié la propriété i).
Réciproquement, considérons $\calM = k\cro u / u^{ep} 
e_1 \oplus k\cro u / u^{ep} e_2$, $\Fil^r \calM = \pa{e_1, u^2 e_2}$,  
$\Phi_r \pa{e_1} = e_2$, $\Phi_r \pa{u^2 e_2} =  e_1 + u e_2$. D'après 
la proposition \ref{prop:nnilp}, un opérateur de monodromie sur $\calM$ 
doit vérifier $N\pa{e_1 + u e_2} = 0$ et $N\pa{e_2} = a \pa{e_1 + u 
e_2}$ pour un certain $a \in k\cro{u^p}/u^{ep}$. On doit avoir $c N 
\circ \phi_r \pa{e_1} = \phi_r \pa{u^e N\pa{e_1}}$, ce qui donne après 
calcul $a = - u^{p\pa{e-1}}$. Il existe donc un unique $N$ valable, et 
il n'est pas nul sur $\im \phi_r$. Cela conclut.
\end{preuve}

\medskip

\begin{prop}
\label{prop:sobj}
Tout objet non nul de $\Mrtilde_0$ (resp. de $\Mrtilde$) admet un 
sous-objet non nul dans $\MFrtilde_0$ (resp. dans $\MFrtilde$ pour
lequel $N$ est nul sur $\Im \Phi_r$).
\end{prop}

\begin{preuve}
La preuve de cette propriété est donnée dans le paragraphe 
\ref{sec:simpn} lors de l'étude des objets simples.
\end{preuve}

\subsubsection*{Les objets de $\MFrtilde$ comme modules filtrés}

Il est possible de décrire la catégorie $\MFrtilde$ avec des 
objets plus proches des objets de Fontaine-Laffaille du cas $e=1$. Soit 
$\calM$ un objet de $\MFrtilde$. Posons $\calM_0 = \Im \Phi_r$ et 
$\calM_i = u^i \calM_0$, pour $0 \leq i \leq p-1$. On a un isomorphisme 
$f_i : \calM_0 \to \calM_i$ qui est la multiplication par $u^i$. 
Définissons $F_{i/e} = f_{er-i}^{-1} \pa{ \Fil^r \calM \cap \calM_i}$. 
On obtient une suite 
décroissante de sous-$k\cro {u^p}/u^{ep}$-modules de $\calM_0$ contenant 
$u^p \calM_0$ telle que $F_0 = \calM_0$ par hypothèse.

L'application $\Phi_r$ induit des applications $\phi_i : F_i \to 
\calM_0$ faisant commuter les diagrammes :
$$\xymatrix @C=50pt {
F_{i+\frac 1 e} \ar[r] \ar[d]_-{\phi_{i+\frac 1 e}} & F_i 
\ar[d]^-{\phi_i} \\
\calM_0 \ar[r]^-{u^p} & \calM_0 }$$

La monodromie, quant à elle, définit une application $c_\pi N : \calM_0 
\to \calM_0$ qui vérifie $\phi(c_\pi) \cro {\pa{c_\pi N} \circ \phi_i} = 
\phi_{i-1} \circ \pa{c_\pi N}$.

Si on remarque pour finir que $k\cro{u^p}/u^{ep}$ est isomorphe en tant 
qu'anneau à $\O_K/p$ (en envoyant $u^p$ sur $\pi$),
on obtient la proposition suivante qui énonce précisément le pont 
entre les catégories $\Mr$ et celles de Fontaine-Laffaille, du moins 
dans le cas modulo $p$ :

\begin{prop}
La catégorie $\MFrtilde$ est équivalente à la catégorie dont les objets 
sont les données suivantes :
\begin{enumerate}
\item un $\O_K / p$-module libre de rang fini $\calM$ ;
\item une filtration décroissante $\pa{F_i}$ de sous-modules 
de $\calM_0$ contenant $\pi \calM$ indexée par les rationnels de dénominateur 
$e$ compris entre $0$ et $r$ telle que $F_0 = \calM$ ;
\item des applications $\Phi$-semi-linéaires $\phi_i : F_i 
\to \calM$ vérifiant ${\phi_i}_{|F_{i+\frac 1 e}} = \pi \phi_{i + 
\frac 1 e}$ ;
\item d'une application linéaire $c_\pi N : \calM \to \calM$ telle que 
$\phi(c_\pi) \cro{ \pa{ c_\pi N} \circ \phi_i } = \phi_{i-1} \circ \pa{
c_\pi N}$ pour tout $i$
\end{enumerate}
et où les flèches sont les morphismes $\O_K / p$-linéaires compatibles à 
toutes les structures.
\end{prop}

\medskip

\noindent
{\it Remarque.} Dans le cas non ramifié (\emph{i.e.} $e=1$ et $c_\pi = 
1$), on retrouve exactement la description des catégories de 
Fontaine-Laffaille modulo $p$ (voir \cite{fontaine-laffaille}). On peut 
étendre cette remarque à toute une sous-catégorie de $\Mr$ comme expliqué 
dans le paragraphe 2.4.1 de \cite{breuil-ens}.

\bigskip

On peut résumer tout ce qui précède par le diagramme suivant :
$$\xymatrix @C=50pt {
\MFrtilde_0 \ar@{^(->}[r] \ar@/_/@{_(->}[d] & \Mrtilde_0 \\
\MFrtilde \ar[u] \ar@{^(->}[r] & \Mrtilde \ar[u] \\
}$$
Les flèches qui montent correspondent aux foncteurs d'oubli évidents. La 
flèche courbe est un foncteur \og $N$ canonique \fg\ qui munit un objet 
de $\MFrtilde_0$ du $N$ (nécessairement unique) donné par le iii) de la 
proposition \ref{prop:mfrtilde}. On pourrait se demander s'il est 
possible de prolonger ce foncteur à tout $\Mrtilde_0$. La réponse est
oui dans le cas $r=1$ (voir le lemme 5.1.2. de \cite{breuil-amer}), et 
non dans le cas général puisqu'il n'est déjà pas vrai que le foncteur 
d'oubli $\Mrtilde \to \Mrtilde_0$ est toujours essentiellement 
surjectif (reprendre l'exemple donné dans la démonstration de la 
proposition \ref{prop:catmfr})

\subsection{Un mot sur le cas $r=1$}
Ce cas est amplement discuté dans \cite{breuil-annals}. Plus exactement, 
Breuil construit là un foncteur contravariant 
entre la catégorie $\widetilde{\underline \M}^1_0$ et la catégorie des 
schémas en groupes finis et plats sur $\O_K$ tués par $p$. Il prouve 
ensuite, en exhibant en quasi-inverse, que ce foncteur est une 
anti-équivalence de catégories.

\medskip

Il étend par la suite ce foncteur à toute la catégorie $\underline 
\M^1_0$ et atteint tous les schémas en groupes sur $\O_K$ tués par une 
puissance de $p$. De cette façon, Breuil retrouve la classification des 
schémas en groupes sur $\O_K$ débutée par Raynaud (\cite{raynaud}) et 
poursuivie par Fontaine (\cite{fontaine-schgrp}) et Conrad 
(\cite{conrad}), et étend même cette classification sans restriction sur 
la ramification.

\subsection{Des catégories abéliennes et artiniennes}
Nous montrons dans ce paragraphe que les catégories $\Mrtilde$ et $\Mr$ 
sont abéliennes. Nous rappelons dans un premier temps que ce résultat
est prouvé dans \cite{breuil-ens} lorsque $e=1$.

\subsubsection*{La catégorie $\Mrtilde$}

Notons $\Mrtilde_{\pa 1}$ la catégorie $\Mr$ obtenue en considérant 
des $k\cro u / u^p$-modules libres à la place de $k \cro u / 
u^{ep}$-modules libres. Un raisonnement rigoureusement 
identique\footnote{Seule la relation de commutation des opérations $N$ et 
$\phi_r$ diffère entre la catégorie $\Mrtilde_{\pa 1}$ et celle 
introduite dans \cite{breuil-ens}, mais l'opérateur $N$ n'intervient pas 
dans la preuve dont il est question.} à celui établi pour prouver le 
corollaire 2.2.3.2 de \cite{breuil-ens} donne :

\begin{theo}
\label{th:abel1}
La catégorie $\Mrtilde_{\pa 1}$ est abélienne et artinienne. Plus 
précisément soit $f : \calX \to \calY$ un morphisme dans $\Mrtilde_{\pa 1}$, 
alors :
\begin{itemize}
\item[i)] $f \pa {\Fil^r \calX} = \Fil^r \calY \cap f\pa \calX$ ;
\item[ii)] Soit $\calK$ le noyau de l'application $k\cro u/u^p$-linéaire
sous-jacente, $\Fil^r \calK = \Fil^r \calX \cap \calK$, $\phi_r : \Fil^r \calK \to \calK$ la 
restriction de $\phi_r : \Fil^r \calX \to \calX$ et $N : \calK \to \calK$ la restriction 
de $N : \calX \to \calX$. Avec ces structures, $\calK$ est un objet de $\Mrtilde_ 
{\pa 1}$ et donne le noyau de $f$ dans $\Mrtilde_{\pa 1}$ ;
\item[iii)] Soit $\calC$ le conoyau de l'application $k\cro u / 
u^p$-linéaire sous-jacente, $\Fil^r \calC$ l'image de $\Fil^r \calY$ dans $\calC$, 
$\phi_r : \Fil^r \calC \to \calC$ l'application qu'induit $\phi_r : \Fil^r \calY \to 
\calY$ et $N : \calC \to \calC$ le quotient de $N : \calY \to \calY$. Avec ces structures, 
$\calC$ est un objet de $\Mrtilde_{\pa 1}$ et donne le conoyau de $f$ dans 
$\Mrtilde_{\pa 1}$.
\end{itemize}
\end{theo}

Nous allons à présent montrer les propriétés analogues pour la catégorie 
$\Mrtilde$.

\medskip

Soit $f : \calX \to \calY$ un morphisme de la catégorie $\Mrtilde$. Notons $\bar 
\calX$ (resp. $\bar \calY$) la réduction de $\calX$ (resp. de $\calY$) modulo $u^p$, et 
$p_\calX : \calX \to \bar \calX$ (resp. $p_\calY : \calY \to \bar \calY$) la projection 
correspondante. On munit $\bar \calX$ et $\bar \calY$ de $\Fil^r$, Frobenius
et opérateurs de monodromie en regardant les structures quotients. On 
obtient des objets de la catégorie $\Mrtilde_{\pa 1}$ et la 
flèche $f$ induit un morphisme $\bar f : \bar \calX \to \bar \calY$ dans cette 
catégorie. Finalement, puisque $er < p-1$, on a :
$$\Fil^r \calX = p_\calX^{-1} \pa{\Fil^r \bar \calX} \quad \text{et} \quad
\Fil^r \calY = p_\calY^{-1} \pa{\Fil^r \bar \calY}$$

\begin{lemme}
L'image (au sens classique) de $f$ est un $k\cro u / u^{ep}$-module 
libre.
\end{lemme}

\begin{preuve}
Comme $f$ commute à $\Phi_r$, elle induit une application $f : \Phi_r 
\pa {\Fil^r \calX} \to \Phi_r \pa {\Fil^r \calY}$ qui est 
$k\cro{u^p}/u^{ep}$-linéaire. En recopiant l'argument de la preuve 
de la proposition \ref{prop:baseadapt}, on prouve qu'il existe des 
éléments $e_1, 
\ldots, e_d, e'_1, \ldots, e'_{d'}$ et des entiers $n_1, \ldots, 
n_d, n'_1, \ldots, n'_{d'}$ tels que $n'_{i'} > 0$ et :
\begin{eqnarray*}
\im f & = & k\cro u / u^{ep} e_1 \oplus \ldots \oplus k\cro u / u^{ep} 
e_d \oplus u^{p n'_1} k\cro u / u^{ep} e'_1 \oplus \ldots \oplus
u^{p n'_{d'}} k\cro u / u^{ep} e'_{d'} \\
\Fil^r \calY \cap \im f & = & u^{n_1} k\cro u / u^{ep} e_1 \oplus 
\ldots \oplus u^{n_d} k\cro u / u^{ep} e_d \oplus u^{p n'_1} k\cro u / 
u^{ep} e'_1 \oplus \ldots \oplus u^{p n'_{d'}} k\cro u / u^{ep} e'_{d'} 
\end{eqnarray*}
De plus, comme $\Phi_r(\Fil^r \calX)$ doit engendrer $\calX$ et que 
$f(\Fil^r \calX) \subset \Fil^r \calY$, on en déduit que $\Phi_r( \Fil^r 
\calY \cap \im f )$ doit au moins engendrer $\im f$. Mais, puisque $er < 
p-1$, on a forcément $\Phi_r (u^p \calY) \subset u^{2p} \calY$ et par
un argument de dimension, le seul moyen de tout concilier est d'avoir 
$d' = 0$, ce qui achève la démonstration.
\end{preuve}

\begin{lemme}
Le noyau et le conoyau de $f$ sont des $k\cro u / u^{ep}$-modules 
libres.
\end{lemme}

\begin{preuve}
En tant que $k\cro u / u^{ep}$-module, l'image de $f$ s'identifie au 
quotient $\calX/\ker f$ et le conoyau au quotient $\calY/ \Im 
f$. Le lemme résulte du fait que si $N \subset M$ est des $k\cro u / 
u^{ep}$-modules de type fini et si deux modules parmi $M$, $N$ et $M/N$ 
sont libres sur $k\cro u / u^{ep}$, alors il en est de même du 
troisième.
\end{preuve}

\bigskip

Notons $\calK$ le noyau de $f$, $\calC$ le conoyau de $f$, $\bar \calK$ 
le noyau de $\bar f$ et $\bar \calC$ le conoyau de $\bar f$ et 
considérons le diagramme suivant :
$$\xymatrix @C=40pt {0 \ar[r] & \ar[r] \ar[d]_{p_\calK} \calK & \ar[r]^f
\ar[d]_{p_\calX} \calX
& \ar[r] \ar[d]_{p_\calY} \calY & \ar[r] \ar[d]_{p_\calC} \calC & 0 \\
0 \ar[r] & \ar[r] \bar \calK & \ar[r]^{\bar f} \bar \calX & \ar[r] \bar 
\calY & \ar[r] \bar \calC
& 0 }$$

\begin{lemme}
\label{lem:noyaumodp}
La flèche $p_\calK$ (resp. $p_\calC$) définie par le diagramme 
précédent est surjective et de noyau $u^p \calK$ (resp. $u^p \calC$). 
Autrement dit $\bar \calK$ s'identifie à $\calK / u^p \calK$ et $\bar 
\calC$ à $\calC / u^p\calC$.
\end{lemme}

\begin{preuve}
Commençons par le noyau et la surjectivité. Soit $\bar x \in \bar 
\calK$. Il se relève en $x \in \calX$ tel que $f\pa x = 0 \pmod {u^p}$. Il
existe donc $y \in \calY$ tel que $f\pa x = u^p y$. On a $u^p y \in \im f$ 
et, puisque $\im f$ est libre sur $k\cro u / u^{ep}$, il existe $y' \in 
\Im f$ tel que $u^p y 
= u^p y'$ et donc $u^p y = f\pa{u^p x'}$ pour un certain $x' \in \calX$.
Mais alors $x-u^p x' \in \calK$ s'envoie sur $\bar x$ par $p_\calK$. Ceci 
prouve la surjectivité.

Soit maintenant $x \in \calK$ tel que $p_\calK\pa x = 0$. On a 
$p_\calX\pa x = 0$ et donc $x$ est un multiple de $u^p$ dans $\calX$. Il 
l'est aussi dans $\calK$ puisque $\calK$ est un $k\cro u/u^{ep}$-module 
libre. Finalement $\ker p_\calK = u^p \calK$.

On utilise des arguments analogues pour le conoyau.
\end{preuve}

\bigskip

On définit $\Fil^r \calK = \calK \cap \Fil^r \calX$, un Frobenius $\Phi_r 
: \Fil^r \calK \to \calK$ et un opérateur de monodromie $N : \calK \to \calK$ déduits 
des opérateurs sur $\calX$. De même, on définit $\Fil^r \calC$ comme l'image de 
$\Fil^r \calY$ par la projection $\calY \to \calC$, un Frobenius et un opérateur de 
monodromie sur $\calC$, les opérateurs sur $\calY$ passant au quotient.

\begin{lemme}
Munis des structures précédentes, les objets $\calK$ et $\calC$ sont des objets 
de la catégorie $\Mrtilde$ et respectivement un noyau et un conoyau de 
l'application $f$.
\end{lemme}

\begin{preuve}
Les conditions de compatibilité et le fait que si les objets sont dans
la catégorie, ils sont noyau ou conoyau est évident. Le seul point 
délicat est la \og surjectivité \fg\ des $\Phi_r$.

Modulo $u^p$, les objets $\calK$ et $\calC$ avec toutes leurs structures 
se réduisent d'après le lemme \ref{lem:noyaumodp} sur $\bar \calK$ et $\bar 
\calC$ et on
sait alors que les $\Phi_r$ définis sur ces objets sont \og surjectifs\fg. 
Notons $\pa{e_1, \ldots, e_d}$ une base adaptée de $\calK$ (qui existe bien) 
et $G$ la matrice de $\Phi_r$ dans cette base. Cette matrice est 
inversible modulo $u^p$ et donc son déterminant est inversible modulo 
$u^p$ puis modulo $u^{ep}$. La matrice $G$ est donc inversible et $\Im 
\Phi_r$ engendre bien tout $\calK$. On raisonne de même pour $\calC$.
\end{preuve}

\begin{cor}
La catégorie $\Mrtilde$ est abélienne et artinienne. Plus 
précisément soit $f : \calX \to \calY$ un morphisme dans $\Mrtilde$, 
alors :
\begin{itemize}
\item[i)] $f \pa {\Fil^r \calX} = \Fil^r \calY \cap f\pa \calX$ ;
\item[ii)] Soit $\calK$ le noyau de l'application $k\cro 
u/u^{ep}$-linéaire sous-jacente, $\Fil^r \calK = \Fil^r \calX \cap 
\calK$, $\phi_r : \Fil^r \calK \to \calK$ la restriction de $\phi_r : 
\Fil^r \calX \to \calX$ et $N : \calK \to \calK$ la restriction de $N : 
\calX \to \calX$. Avec ces structures, $\calK$ est un objet de 
$\Mrtilde$ et donne le noyau de $f$ dans $\Mrtilde$ ;
\item[iii)] Soit $\calC$ le conoyau de l'application $k\cro u / 
u^{ep}$-linéaire sous-jacente, $\Fil^r \calC$ l'image de $\Fil^r \calY$ 
dans $\calC$, $\phi_r : \Fil^r \calC \to \calC$ l'application qu'induit 
$\phi_r : \Fil^r \calY \to \calY$ et $N : \calC \to \calC$ le quotient 
de $N : \calY \to \calY$. Avec ces structures, $\calC$ est un objet de 
$\Mrtilde$ et donne le conoyau de $f$ dans $\Mrtilde$.
\end{itemize}
\end{cor}

\begin{preuve}
On a déjà prouvé ii) et iii). Il ne reste en fait plus qu'à démontrer i) 
car il implique l'isomorphisme entre image et coimage. On a 
évidemment toujours l'inclusion $f \pa {\Fil^r \calX} \subset \Fil^r \calY \cap 
f\pa \calX$.

Soit $y \in \Fil^r \calY \cap f\pa \calX$. La réduction $\bar 
y$ de $y$ modulo $u^p$ est un élément de $\Fil^r \bar \calY \cap \bar f\pa 
{\bar \calX}$ (en gardant les notations précédentes) et d'après le 
théorème \ref{th:abel1}, $\bar y \in \bar f \pa{\Fil^r \bar 
\calX}$. Il existe $\bar x \in \Fil^r \bar \calX$ tel que $\bar y = \bar f 
\pa {\bar x}$. Notons $x$ un relevé de $\bar x$ dans $\Fil^r \calX$. Il 
existe un élément $t \in \calY$ tel que $y = f\pa x + u^p y'$. Les 
éléments $y$ et $f\pa x$ sont dans $\im f$, il en est donc de même de 
$u^p y'$ et puisque $\im f$ est libre sur $k\cro u / u^{ep}$, il existe 
$y'' \in \im f$ tel 
que $u^p y' = u^p y''$. On écrit $y'' = f\pa{x''}$ pour un certain $x'' 
\in \calX$, et il vient $y = f\pa{x+u^p x''}$. Comme $u^{er} \calX \subset 
\Fil^r \calX$ et $er < p-1$, on a $u^p x'' \in \Fil^r \calX$ et donc $x+u^p x'' 
\in \Fil^r \calX$. Finalement $y \in f\pa{\Fil^r \calX}$ et on peut conclure.
\end{preuve}

\subsubsection*{La catégorie $\Mr$}

Si $r=0$, on vérifie facilement (la propriété i) du corollaire précédent 
est alors évidente) que la sous-catégorie pleine de $\Mr$ formé des 
objets tués par $p$ est abélienne et artinienne. On procède ensuite par 
dévissage. La preuve est en tout point analogue à celle déjà connue dans 
le cas $e=1$ et présentée dans le paragraphe 2.3 de \cite{breuil-ens}. 
Les lemmes et les propositions successives gardent un sens dans ce 
contexte plus général, et sont également vrais, les preuves étant encore 
textuellement les mêmes. Nous n'insisterons donc pas davantage et 
laissons le lecteur se reporter à cette référence.

\bigskip

\noindent
{\it Remarque.} De même, on prouve que les catégories $\Mrtilde_0$ et 
$\Mrtilde$ et $\Mr_0$ sont abéliennes et artiniennes.

\section{Classification des objets simples}
\label{sec:simples}
Nous allons donner une classification complète des objets simples de la
catégorie $\Mr$ lorsque le corps résiduel $k$ est algébriquement clos.
Nous essaierons également d'expliquer ce qui se passe lorsque ce n'est 
pas le cas. Pour l'instant, on ne fait aucune hypothèse supplémentaire 
sur $k$.

\bigskip

On considère $\calM$ un objet simple (donc non nul) de $\Mr$. Il est 
obligatoirement tué par $p$. En effet, si ce n'était pas le cas, le 
noyau de la multiplication par $p$ dans $\calM$ fournirait un sous-objet 
strict de $\calM$ (noter que la multiplication par $p$ ne peut pas être 
injective car elle est nilpotente : $\calM$ est supposé être tué par une 
puissance de $p$).
L'objet simple $\calM$ peut être vu dans la catégorie $\Mrtilde$ (du 
moins si $r > 0$, mais dans le cas contraire, le résultat est immédiat
et laissé au lecteur) :
c'est un $k\cro u / u^{ep}$-module muni d'un $\Fil^r$, d'un $\Phi_r$ 
et d'un opérateur de monodromie vérifiant les bonnes propriétés.

\subsection{La monodromie}
\label{sec:simpn}

Si l'on note $\calM_0 = \Im \Phi_r$, l'application de monodromie $N$ 
induit une application linéaire $c_\pi N : \calM_0 \to \calM_0$ 
(voir proposition \ref{prop:n}) et il existe $x_1 \in \calM_0 
\backslash u \calM_0$ tel que $N\pa{x_1} = 0$ (voir corollaire 
\ref{cor:kern}). Notons $x_2 = \varphi_r \pa{x_1}$ (voir définition 
\ref{def:psi}) puis par récurrence $x_{i+1} = \varphi_r \pa{x_i}$, 
ce qui est possible d'après la deuxième partie de la proposition 
\ref{prop:phir}. On a $N\pa {x_i} = 0$ pour tout entier $i$.

\medskip

Notons $\bar x_i$ la réduction modulo $u^p$ de $x_i$. Les $\bar x_i$ 
sont des éléments non nuls de $\calM_0/u^p \calM_0$ qui est un $k$-espace 
vectoriel de dimension finie. Notons 
$n \geq 1$ le plus petit indice tel que $\bar x_{n+1}$ puisse s'écrire 
comme combinaison linéaire des $\bar x_i$ pour $i$ variant de $1$ à 
$n$. Il existe donc $\lambda_i \in k$ tels que :
$$\bar x_{n+1} = \lambda_1 \bar x_1 + \ldots + \lambda_n \bar x_n$$
et on peut supposer $\lambda_1 \neq 0$ quitte à remplacer $x_1$ par 
le plus petit indice $i$ tel que $\lambda_i \neq 0$. Comme $\bar 
x_{n+1} \neq 0$, les $\lambda_i$ ne peuvent être tous simultanément
nuls.

\medskip

Nous allons à présent corriger les $x_i$ pour que cette relation ne soit 
plus vraie seulement modulo $u^p$. On procède par approximations 
successives et on construit une suite indexée par $j$ d'éléments 
$x_i^{\pa j}$ qui sont tels que $x_i^{\pa j} \equiv x_i \pmod 
{u^p}$, $x_{i+1}^{\pa j} = \varphi_r (x_i^{\pa j})$, $N (x_{i+1}^{\pa 
j}) = 0$ et finalement :
$$x_{n+1}^{\pa j} \equiv \lambda_1 x_1^{\pa j} + \ldots + \lambda_n 
x_n^{\pa j} \pmod {u^{jp}}$$
les $\lambda_i$ restant inchangés. On a une solution pour $j=1$. 
Supposons qu'on l'ait pour $j$ et construisons-en une pour $j+1$. On 
cherche un élément $r \in \calM_0$ tel que l'on puisse poser $x_1^{\pa 
{j+1}} = 
x_1^{\pa j} + u^{jp} r$. On définirait alors les $x_i^{\pa {j+1}}$ 
\emph{via} la formule de récurrence $x_{i+1}^{\pa {j+1}} = \varphi_r ( 
x_i^{\pa {j+1}} )$ et il est facile de vérifier que pour tout $i\geq 
2$, on aurait $x_i^{\pa {j+1}} \equiv x_i^{\pa j} \pmod 
{u^{\pa{j+1}p}}$. Au final, il suffit de trouver $r$ tel que :
$$x_{n+1}^{\pa j} \equiv \lambda_1 \pa{x_1^{\pa j} + u^{jp} r} + 
\lambda_2 
x_2^{\pa j} + \ldots + \lambda_n x_n^{\pa j} \pmod {u^{\pa{j+1}p}}$$
mais comme par hypothèse, on a $x_{n+1}^{\pa j} \equiv \lambda_1
x_1^{\pa j} + \ldots + \lambda_n x_n^{\pa j} \pmod {u^{jp}}$, on a bien 
l'existence d'un tel $r$ : il suffit de le prendre tel que 
$u^{jp} \lambda_1 r = x_{n+1}^{\pa j} - \lambda_1 x_1^{\pa j} - 
\ldots - \lambda_n x_n^{\pa j}$. De plus, en appliquant $N$ à cette 
dernière égalité, on voit que $N\pa{u^{jp} \lambda_1 r} = \lambda_1 
u^{jp} N\pa r = 0$ et donc que $N(x_1^{\pa {j+1}}) = 0$
puisque $\lambda_1$ est supposé non nul. Ceci implique la nullité de 
tous les $N(x_i^{\pa {j+1}})$.

\medskip

Pour $j=e$, l'égalité a lieu modulo $u^{ep}$ et donc dans $\calM$. Soit 
$\calK$ le sous-$k\cro u/u^{ep}$-module engendré par les $x_i^{\pa e}$, 
pour $1 \leq i \leq n$. La liberté sur $k$ des $\bar x_i$ assure que 
$\calK$ est un module libre de rang $n$. Notons $\calM_i = u^i \Im \Phi_r$ 
et $\Fil^r \calK = \sum_{i=1}^{p-1} \calM_i \cap \calK$. Par construction, 
$N$ stabilise $\calK$ et $\Phi_r$ envoie $\Fil^r \calK$ sur $\calK$. En 
outre, encore par construction l'image de la restriction de $\Phi_r$ à 
$\Fil^r \calK$ engendre $\calK$ : on voit que l'objet $\calK$ est dans 
la catégorie $\MFrtilde$. Comme $\calM$ est simple, ce sous-objet est tout 
$\calM$. Ainsi on a prouvé la proposition suivante :

\begin{prop}
\label{prop:simpmf}
Soit $\calM$ un objet simple de $\Mr$. Alors $\calM$ est dans la catégorie 
$\MFrtilde$ et l'opérateur de monodromie $N$ est nul sur $\Im \Phi_r$. 
De plus $\calM$ admet une base adaptée de la forme $\pa{x_1, \ldots, x_d}$ 
telle que $N\pa {x_i} = 0$, $x_{i+1} = \varphi_r \pa{x_i}$ (voir 
définition \ref{def:psi}) et :
$$\varphi_r\pa{x_d} = \lambda_1 x_1 + \ldots + \lambda_d x_d$$
où les $\lambda_i$ sont des éléments de $k$ tels que $\lambda_1 \neq 0$.
\end{prop}

\noindent qui admet pour corollaire immédiat la proposition 
\ref{prop:sobj} que l'on vient donc de démontrer.

\subsection{Une base adaptée simple}

Nous allons dans ce paragraphe préciser un peu plus l'énoncé de la 
proposition \ref{prop:simpmf} dans le cas où le corps résiduel $k$ est 
algébriquement clos.

\begin{lemme}
\label{lem:basesimple}
Supposons $k$ algébriquement clos.
Soit $\calM$ un objet simple de $\Mr$. Alors $\calM$ est dans $\Mrtilde$ et il 
existe $\pa{e_1, \ldots, e_d}$ une base adaptée de $\calM$ telle que $N\pa 
{e_i} = 0$, $e_{i+1} = \varphi_r \pa{e_i}$, les indices $i$ étant 
considérés dans $\Z/d\Z$.
\end{lemme}

\begin{preuve}
On sait d'après la proposition précédente qu'il existe une base adaptée 
$\pa{x_1, \ldots, x_d}$ telle que $\varphi_r\pa{x_i} = x_{i+1}$ pour $i$ 
compris entre $1$ et $d-1$ et $\varphi_r\pa{x_d} = \lambda_1 
x_1 + \ldots + \lambda_d x_d$ où $\lambda_i \in k$ et $\lambda_1 \neq 
0$.

Parmi toutes les bases adaptées qui vérifient ces conditions, 
choisissons-en une pour laquelle le nombre de $\lambda_i$ non nuls est
minimal. On écrit alors plutôt :
$$x_{d+1} = \varphi_r\pa{x_d} = \lambda_1 x_{i_1} + \ldots + \lambda_k 
x_{i_k}$$
où tous les $\lambda_i$ sont non nuls et les indices $i_k$ 
sont compris entre $0$ et $d$. De plus, on a $i_1 = 1$.
Notons pour tout $i$, $n_i$ le plus petit entier tel que $u^{n_i} x_i 
\in \Fil^r \calM$. Si tous les $n_{i_j}$ n'étaient pas égaux, 
$\varphi_r \pa{x_{d+1}}$ s'écrirait comme une combinaison linéaire de 
$x_2, \ldots, x_{d+1}$ faisant intervenir strictement moins de $k$ 
termes et la famille $\pa{x_2, \ldots, x_{d+1}}$ fournirait une base 
adaptée de $\calM$ (en reprenant l'étude faite de le paragraphe 
précédent). Mais ceci est en contradiction avec la minimalité 
considérée.

Ainsi tous les $n_{i_j}$ sont égaux et donc égaux à $n_{d+1}$. Par 
récurrence, on prouve que pour tout entier fixé $a$, tous 
les $n_{i_j+a}$ sont égaux, les indices $i_j+a$ étant considérés modulo 
$d$.

Notons $t$ le plus grand commun diviseur de $d$ et de toutes les 
différences $i_j - i_{j'}$. D'après ce qui précède la suite des $n_i$ 
est périodique de période (divisant) $t$. On considère alors le 
sous-$k$-espace vectoriel de $\calM$ engendré par les $x_{tn}$ où $n$ 
parcourt $\Z/d\Z$. L'application $\varphi_r^t$ stabilise 
ce sous-espace et y est $\Phi^t$-semi-linéaire. En particulier, puisque
$k$ est algébriquement clos, il existe un élément $e_1$ de ce 
sous-espace tel que $\varphi_r^t\pa{e_1} = \lambda e_1$ pour un certain 
$\lambda \in k^\star$. Quitte à multiplier $e_1$ par un élément de $k$, 
on peut supposer $\lambda = 1$.

On définit $e_{i+1} = \varphi_r \pa{e_i}$. La famille $\pa{e_1, \ldots, 
e_t}$ engendre un espace stable par $N$ et par $\Phi_r$ qui est par 
construction un sous-objet non nul de $\calM$. C'est donc tout $\calM$. De plus, 
d'après la proposition \ref{prop:simpmf}, $N\pa{e_i} = 0$ pour tout $i$. 
Finalement $\pa{e_1, \ldots, e_t}$ est une base adaptée vérifiant les 
conditions du lemme. Cela conclut.
\end{preuve}

\subsection{Classification proprement dite}
\label{sec:simpalgclos}

\begin{deftn}
Soit $\pa{n_i}$ une suite périodique\footnote{Par \og 
périodique \fg, on entend dans ce papier \og 
périodique dès le début \fg\ et pas \og périodique à partir d'un certain 
rang \fg.} d'entiers compris entre $0$ et $er$. On note $h$ la 
période de cette suite. On définit l'objet $\calM\pa{n_i} \in \Mrtilde$ de 
la façon suivante :
\begin{enumerate}
\item $\displaystyle \calM\pa{n_i} = \bigoplus_{i\in \Z/h\Z} k\cro u / 
u^{ep} e_i$ ;
\item $\displaystyle \Fil^r \calM\pa{n_i} = \bigoplus_{i\in \Z/h\Z} 
u^{n_i}k\cro u / u^{ep} e_i$ ;
\item $\Phi_r \pa {u^{n_i} e_i} = e_{i+1}$ pour tout indice $i$ ;
\item $N\pa{e_i} = 0$ pour tout indice $i$.
\end{enumerate}
\end{deftn}

Il est facile de vérifier que tous ces objets sont bien dans la 
catégorie $\Mrtilde$ et on a le théorème suivant :

\begin{theo}
\label{th:simple}
Supposons $k$ algébriquement clos.
Les objets $\calM\pa{n_i}$ sont des objets simples de la catégorie 
$\Mrtilde$. De plus, si $\calM$ est un objet simple de la catégorie 
$\Mrtilde$, alors il est isomorphe à un certain $\calM\pa{n_i}$.
\end{theo}

\begin{preuve}
Voyons d'abord la simplicité de $\calM\pa{n_i}$. Soit $\calM$ un 
sous-objet non nul de $\calM\pa{n_i}$. L'image de la restriction de 
$\Phi_r$ à $\calM \cap \Fil^r \calM\pa{n_i}$ est supposée engendrer tout $\calM$ ; 
en particulier elle 
n'est pas réduite à $0$ et comprend un élément non divisible par $u^p$, 
disons $x$. On écrit $x = \lambda_1 e_1 + \ldots + \lambda_h e_h$ 
où les $\lambda_i$ sont des polynômes à coefficients dans 
$k\cro{u^p}/u^{ep}$. On peut supposer $\lambda_i \in k$ 
quitte à remplacer $x$ par $\varphi_r \circ \varphi_r \pa x$.

Considérons un $x$ pour lequel le nombre de $\lambda_i$ non nuls est 
minimal et écrivons :
$$x = \lambda_1 e_{i_1} + \ldots + \lambda_k e_{i_k}$$
avec ici tous les $\lambda_i$ non nuls. En appliquant 
$\varphi_r$ éventuellement plusieurs fois, on voit que tous les 
$n_{i_j}$ doivent être égaux car sinon, on obtient un nouvel $x$
qui serait combinaison d'un nombre plus petit de $e_i$. 
On applique alors $\varphi_r$ à l'égalité précédente et comme 
précédemment, on prouve que tous les $n_{i_j+1}$ sont égaux. Par 
récurrence, on voit que pour $a$ fixé tous les $n_{i_j+a}$ sont 
égaux.
Ainsi, pour que la suite $\pa{n_i}$ soit périodique de période 
exactement $h$, il faut que $k=1$, c'est-à-dire que $x$ soit multiple de 
l'un des $e_i$. Mais alors le sous-objet engendré par $x$ est tout 
$\calM\pa{n_i}$ et finalement $\calM = \calM\pa{n_i}$. Ce qui assure la 
simplicité.

\medskip

Voyons la réciproque. On applique le lemme \ref{lem:basesimple} qui
donne une description explicite de l'objet $\calM$. Il reste juste 
à démontrer que la suite $\pa{n_i}$ ne peut-être périodique de période 
divisant strictement $h$. Mais supposons que ce soit le cas et notons
$t$ cette période. On considère le sous-objet engendré par l'élément $x 
= e_{t} + e_{2t} + \ldots + e_{h}$ et on vérifie immédiatement qu'il est 
non nul et strictement inclus dans $\calM$. C'est une contradiction.
\end{preuve}

\bigskip

\noindent
{\it Remarque.} En utilisant la correspondance de \cite{breuil-annals}, 
on retrouve exactement la classification donnée par Raynaud dans 
\cite{raynaud}.

\subsubsection*{Étude des endomorphismes}

On suppose toujours le corps $k$ algébriquement clos.
Soit $\calM$ un objet simple de $\Mrtilde$. Soit $\pa{e_1, \ldots, e_h}$ une 
base adaptée de $\calM$ vérifiant les conditions du théorème
\ref{th:simple}. Nous allons en fait voir que les $e_i$ sont 
presque uniquement déterminés. Plus précisément, on a :

\begin{theo}
Supposons $k$ algébriquement clos.
Si $\lambda \in k$ vérifie $\lambda^{p^h} = \lambda$, alors 
l'application $\psi : \calM \to \calM$ définie par $\psi\pa{e_i} = 
\lambda^{p^i} e_i$ est un endomorphisme de $\calM$. Ce sont les seuls.
\end{theo}

\begin{preuve}
Déjà il est facile de vérifier que les applications définies dans 
l'énoncé du théorème sont bien compatibles au $\Fil^r$, au Frobenius 
et à l'opérateur de monodromie.

\medskip

Pour la réciproque, nous allons anticiper sur des résultats ultérieurs
donnant une application non nulle $\End\pa \calM \to \End\pa{\Tst\pa \calM}$
(la non-nullité se déduit de la fidélité du foncteur $\Tst$, voir 
corollaire \ref{cor:fidele}). D'autre part $\End\pa \calM$ est un corps 
\emph{a priori} non commutatif et $\End\pa{\Tst\pa \calM}$ est un corps fini 
à $p^h$ éléments (cela se déduit de théorème \ref{th:tstsimple}). On en 
déduit facilement que la flèche $\End\pa \calM \to \End\pa{\Tst\pa \calM}$ est 
bijective, ce qui prouve le théorème.
\end{preuve}

\begin{cor}
Les objets simples $\calM\pa{n_i}$ et $\calM\pa{m_i}$ sont isomorphes si et 
seulement si la suite $\pa{m_i}$ se déduit de la suite $\pa{n_i}$ par 
translation.
\end{cor}

\begin{preuve}
Si deux objets $\calM\pa{n_i}$ et $\calM\pa{m_i}$ sont isomorphes, on peut 
transporter une base adaptée de $\calM\pa{n_i}$ à $\calM\pa{m_i}$ et le 
théorème précédent entraîne la conclusion voulue.
\end{preuve}

\subsubsection*{Un autre point de vue}

Donnons finalement un point de vue différent sur cette classification, 
peut-être plus agréable à retenir.

\medskip

Soit $\pa{n_i}_{i \in \N}$ une suite quelconque d'entiers compris entre 
$0$ et $p-1$ et soit $t$ le rationnel dont le développement \og 
décimal \fg\ en base $p$ est :
$$t = 0 \, , \, n_1 \; n_2 \; n_3 \; n_4 \; \ldots$$

\medskip

On a une propriété classique :

\begin{prp}
Avec les notations précédentes, les suites périodiques sont 
exactement celles qui correspondent aux rationnels de $\Z_{\pa p} \cap 
\cro{0,1}$ où $\Z_{\pa p}$ désigne le localisé de $\Z$ en $p$.
\end{prp}

On peut alors poser la définition suivante :

\begin{deftn}
Soit $\rat$ l'ensemble des classes d'équivalence d'éléments de $\Z_{\pa 
p}$ pour la relation d'équivalence suivante : $a \sim b$ si et 
seulement s'il existe un entier $n$ tel que $a \equiv p^n b \pmod \Z$.
\end{deftn}

La dernière relation d'équivalence n'est pas mystérieuse : elle 
correspond simplement à un décalage des décimales du nombre. En 
particulier, à cause de la périodicité, les classes d'équivalence sont 
toutes finies.

\bigskip

Dans ces conditions, $\rat$ classifie exactement les objets 
simples de la catégorie $\Mr$ (\emph{via} la correspondance que l'on a 
décrite précédemment).

\medskip

Nous verrons par la suite que le \og\emph{rationnel classifiant}\fg\ 
va réapparaître de façon naturelle.

\section{Étude du foncteur $\Tst$}

\subsection{Un système préliminaire}
\label{sec:systeme}
Ce paragraphe présente une version légèrement différente de résultats 
classiques et par exemple déjà discutés dans \cite{wach} ou dans le 
paragraphe 3.3.2 de \cite{breuil-ens}. On suppose
dans ce paragraphe que le corps résiduel $k$ est algébriquement clos.

\bigskip

On considère un entier $h$ strictement postitif. On fixe $\eta^{\pa h}$ 
une racine $\pa{p^h-1}$-ième de l'uniformisante $\pi$ de $K$ et on 
appelle $K^{\pa h}$ l'extension de $K$ engendrée par cette racine. On 
rappelle que $K^{\pa h}/K$ est totalement et modérément ramifiée de 
degré $p^h - 1$. On rappelle également que la limite inductive de toutes 
ces extensions est l'extension maximale modérément ramifiée 
de $K$. Par la suite, lorsqu'il n'y aura pas de risque d'ambiguité, on 
notera $\eta$ à la place de $\eta^{\pa h}$. On rappelle enfin que 
$\pi_1$ désigne une racine $p$-ième de $\pi$.

\medskip

On s'intéresse au système d'équations suivant :
$$\pa S : \left\{
\begin{array}{lcr}
\frac {\pa{\pi_1^{n_1} \hat x_1 + \hat c_1}^p}{\pi^{er}} & = & \varpi 
\hat x_2 + \hat r_1 \\
\frac {\pa{\pi_1^{n_2} \hat x_2 + \hat c_2}^p}{\pi^{er}} & = & \varpi 
\hat x_3 + \hat r_2 \\
& \vdots \\
\frac {\pa{\pi_1^{n_h} \hat x_h + \hat c_h}^p}{\pi^{er}} & = & \varpi 
\hat x_1 + \hat r_h \\
\end{array} \right.$$
où $\varpi \in \acco{-1,1}$ est un signe, les $n_i$ sont des entiers 
fixés tous compris entre $0$ et $er$, et où les $\hat r_i$ et les $\hat 
c_i$ sont des éléments de $\O_{\bar K}$. Les inconnues sont les $\hat 
x_i$ que l'on cherche également dans $\O_{\bar K}$. On pose dans la 
suite $m_i = er - n_i$.

\subsubsection*{Sans coefficient constant}

On s'intéresse tout d'abord au cas où toutes les constantes $\hat r_i$ 
et $\hat c_i$ sont nulles. Il est alors possible de résoudre directement 
le système dans $\O_{\bar K}$. En effet, le système se réécrit 
simplement :
$$\hat x_i ^p = \varpi \pi^{m_i} \hat x_{i+1}.$$
Par des manipulations simples, on voit que $\hat x_1$ doit être solution 
de l'équation :
$$\hat x_1 ^{p^h} = \varpi^h \pi^{s_1} \hat x_1$$
où $s_1$ est défini par la formule :
$$s_1 = m_1 p^{h-1} + m_2 p^{h-2} + \ldots + m_{h-1} p + m_h.$$
Cette équation admet $p^h$ solutions qui sont $0$ et toutes les racines 
$\pa{p^h-1}$-ièmes de $\varpi^h \pi^{s_1}$. À partir de $\hat x_1$, on 
reconstruit les autres $\hat x_i$ et on vérifie qu'ils forment bien une 
solution du système.

\medskip

On peut présenter les choses de façon plus homogène en procédant comme 
suit. On pose pour tout $i \in \Z/h\Z$ :
$$s_i = m_i p^{h-1} + m_{i+1} p^{h-2} + \ldots + m_{i+h-2} p + 
m_{i+h-1}.$$
Si $\epsilon$ est une racine $\pa{p^h-1}$-ième de $\varpi$ (qui est 
déjà dans $K$), la famille des $\hat x_i = \varpi^i \epsilon^{p^i} 
\eta^{s_i}$ est une solution de $\pa S$. Toutes les solutions 
s'obtiennent ainsi à l'exception de la solution nulle $\hat x_1 = \ldots 
= \hat x_h = 0$.

\subsubsection*{Un lemme à la Hensel}

On ne suppose plus que les constantes sont nulles et on cherche un lien 
entre les solutions de $\pa S$ modulo $p$ et les solutions de $\pa S$ 
dans $\O_{\bar K}$ :

\begin{lemme}
\label{hensel}
Avec les notations précédentes, si le système $\pa S$ admet une 
solution $\pa{x_1, \ldots, x_h}$ modulo $p$, alors cette solution se 
relève dans $\O_{\bar K}$ en une solution $\pa{\hat x_1, \ldots, \hat 
x_h}$.
\end{lemme}

\begin{preuve}
On construit cette solution par approximations successives. Fixons tout 
d'abord une extension finie $L$ de $K^{\pa h}$ suffisamment grande pour 
contenir tous les $\hat r_i$, les $\hat c_i$ et pour que tous les $x_i$ 
puissent s'y relever. L'extension $L/K^{\pa h}$ est totalement ramifiée 
(puisque $k$ est supposé algébriquement clos), disons de degré $d$.
Notons $\O_L$ l'anneau des entiers de $L$.

On va construire une suite de $(x_1^{\pa n}, x_2^{\pa n}, \ldots, 
x_h^{\pa n})$ de solutions compatibles du système $\pa S$ modulo 
$\eta^n$ dans $\O_L$. Il suffira par la suite de prendre la limite de 
cette suite pour avoir une solution du système dans $\O_L$ et donc dans 
$\O_{\bar K}$.

On a déjà, par hypothèse, un $h$-uplet pour $n = e\pa{p^h-1}$. 
Les suivants se construisent par récurrence. On part d'un entier $n \geq 
e\pa{p^h-1}$ et d'éléments $x_1^{\pa n}, \ldots, x_h^{\pa n}$ vérifiant :
$$\frac {\pa{\pi_1^{n_i} x_i^{\pa n}+\hat c_i}^p}{\pi^{er}} \equiv 
\varpi x_{i+1}^{\pa n} + \hat r_i \pmod {\eta^n}$$
pour tout indice $i$ pris dans $\Z/h\Z$ et on cherche à construite $y_1, 
\ldots, y_h$, tels que :
$$\frac {\pa{\pi_1^{n_i} x_i^{\pa n} + \hat c_i + \pi_1^{n_i} \eta^n
y_i}^p}{\pi^{m_i}} \equiv \varpi x_{i+1}^{\pa n} + \eta^n y_{i+1} + \hat 
r_i \pmod {\eta^{n+1}}$$
Un calcul donne :
$$\frac {\pa{\pi_1^{n_i} x_i^{\pa n} + \hat c_i + \pi_1^{n_i} \eta^n 
y_i}^p}{\pi^{er}} = \frac {\pa{\pi_1^{n_i} x_i^{\pa n} + \hat c_i}^p} 
{\pi^{er}} + \sum_{k=1}^p \ncr p k \frac {\pi_1^{k n_i} 
\eta^{kn}}{\pi^{er}} y_i^k \pa{\pi_1^{n_i} x_i^{\pa n} + \hat 
c_i}^{p-k}$$
Soit un entier $k \geq 1$. On a $v\pa{\frac {\pi_1^{kn_i} 
\eta^{kn}}{\pi^{er}}} = \frac {kn_i} p + \frac 
{kn}{p^h-1} - er$. D'autre part, $\frac {(\pi_1^{n_i} x_i^{\pa 
n}+\hat c_i)^p} {\pi^{er}}$ est un entier, donc de valuation positive et 
on en déduit que $v(\pi_1^{n_i} x_i^{\pa n} + \hat c_i) \geq \frac {er} 
p$. On obtient :
\begin{eqnarray*}
v_i = v\pa {\frac {\pi_1^{kn_i} \eta^{kn}}{\pi^{er}} \pa{\pi_1^{n_i} 
x_i^{\pa n} + \hat c_i}^{p-k}} & \geq & 
\frac{kn}{p^h-1} + \frac {k n_i} p - er + \pa{p-k} \frac {er} p \\
& = & \frac n {p^h-1} + \pa{k-1}\frac n{ p^h-1} - \frac {k m_i} p 
\end{eqnarray*}
Comme par hypothèse $n \geq e\pa{p^h-1}$ et $m_i \leq er$, il vient :
$$v_i \geq \frac n {p^h-1} + e\pa{k-1} - \frac {ker} p = \frac n{p^h-1} 
- e + ek \pa{1 - \frac r p} $$
Maintenant si $k < p$, le coefficient binomial $\ncr p k$ est multiple 
de $p$ et donc :
$$v\pa {\ncr p k \frac {\pi_1^{kn_i} \eta^{nk}}{\pi^{er}} 
\pa{\pi_1^{n_i} x_i^{\pa n} + \hat c_i}^{p-k}} \geq e + v_i 
\geq \frac n{p^h-1} + ek \pa{1 - \frac r p} \geq \frac {n+1}{p^h-1}$$
la dernière inégalité résultant du fait que $r \leq er \leq p-2$.

On en déduit que tous les termes de la somme pour $k$ compris 
strictement entre $0$ et $p$ sont nuls modulo $\eta^{n+1}$. En fait, 
c'est aussi le cas pour $k = p$. En reprenant les égalités précédentes, 
on voit que :
$$v_p \geq \frac n {p^h-1} + e\pa{p-1-r}$$
mais $p-1-r \geq 1$ et donc on a également $v_p \geq \frac{n+1} 
{p^h-1}$. Finalement le système que l'on a à résoudre se réduit à :
$$\frac {\pa{\pi_1^{m_i} x_i^{\pa n} + \hat c_i}^p} {\pi^{er}}
\equiv \varpi x_{i+1}^{\pa n} + \eta^n y_i + \hat r_i \pmod 
{\eta^{n+1}}$$
mais on sait que la différence $\frac {(\pi_1^{m_i} x_i^{\pa n} + \hat 
c_i)^p} {\pi^{er}} - \varpi x_{i+1}^{\pa n} - \hat r_i$ est un multiple 
de $\eta^n$, et donc s'écrit $\eta^n q_i$. Il suffit ensuite de choisir 
$y_i = q_i$ pour avoir la solution que l'on cherchait.
\end{preuve}

\subsubsection*{Résolution du système}

Une première conséquence du lemme que l'on vient de prouver est la 
résolution du système $\pa S$ modulo $p$ lorsque les constantes $\hat 
r_i$ et $\hat c_i$ sont toutes nulles :

\begin{lemme}
\label{lem:systnul}
Supposons que les constantes $\hat r_i$ et $\hat c_i$ soient nulles.
Mise à part la solution nulle, les solutions de $\pa S$ dans $\O_{\bar 
K}/p$ s'écrivent $x_i = \varpi^i \epsilon^{p^i} \bar \eta^{s_i}$ où 
$\epsilon \in \O_K$ est une racine $\pa{p^h-1}$-ième de $\varpi^h$, 
$\bar \eta$ 
est la réduction de $\eta$ dans $\O_{\bar K}/p$ et :
$$s_i = m_i p^{h-1} + m_{i+1} p^{h-2} + \ldots + m_{i+h-2} p + 
m_{i+h-1}.$$
\end{lemme}

\begin{preuve}
Si $\epsilon$ est une racine $\pa{p^h-1}$-ième de $\varpi^h$ et si 
$x_i$ désigne la réduction modulo $p$ de $\varpi^i \epsilon^{p^i} 
\eta^{s_i}$, le uplet $\pa{x_1, \ldots, x_h}$ est
solution du système dans $\O_{\bar K}/p$. De plus, si $\epsilon$ et 
$\epsilon'$ sont deux racines 
$\pa{p^h-1}$-ièmes de $\varpi^h$ distinctes, on a pour tout entier $i$, 
$\epsilon^{p^i} \neq \epsilon'^{p^i}$ dans le corps résiduel et donc 
$\epsilon^{p^i} - \epsilon'^{p^i}$ est de valuation nulle. On en déduit,
puisque $v(\eta^{s_i}) = \frac{s_i}{p^h-1}<1$, que $x_i = 
\varpi^i \epsilon^{p^i} \eta^{s_i}$ et $x'_i = \varpi^i \epsilon'^{p^i} 
\eta^{s_i}$ sont distincts dans $\O_{\bar K}/p$.

\medskip

On a ainsi trouvé $p^h$ solutions à $\pa S$ modulo $p$. Le lemme 
\ref{hensel} assure qu'il y en a au moins autant dans 
$\O_{\bar K}$. Mais on a vu qu'il y en a exactement $p^h$ dans $\O_{\bar 
K}$, on les a donc toutes.
\end{preuve}

\bigskip

Passons au cas général. On reprend le système $\pa S$ mais on ne suppose 
plus la nullité de $\hat r_i$ et de $\hat c_i$.

\begin{theo}
\label{th:systeme}
Supposons que le système $\pa S$ admette une solution dans $\O_{\bar 
K}$, alors il admet toujours $p^h$ solutions dans $\O_{\bar K}$ et 
$p^h$ solutions dans $\O_{\bar K}/p$. De plus l'application 
de réduction modulo $p$ définit une bijection entre ces ensembles de 
solutions.

En outre si $\pa{x_1, \ldots, x_h}$ et $\pa{y_1, \ldots, y_h}$ sont deux 
solutions distinctes dans $\O_{\bar K}/p$, alors il existe $\epsilon$ 
une racine $\pa{p^h-1}$-ième de $\varpi^h$ telle que $y_i = x_i + 
\varpi^i \epsilon^{p^i} \eta^{s_i}$ pour tout indice $i \in \Z/h\Z$ où 
$s_i$ est défini par la formule :
$$s_i = m_i p^{h-1} + m_{i+1} p^{h-2} + \ldots + m_{i+h-2} p + 
m_{i+h-1}$$
\end{theo}

\begin{preuve}
Soit $(\hat x_1, \ldots, \hat x_h)$ une solution de $\pa S$ dans 
$\O_{\bar K}$. Si l'on note $x_i \in \O_{\bar K}/p$ la réduction modulo 
$p$ de $\hat x_i$, le uplet $\pa{x_1, \ldots, x_h}$ est solution de 
$\pa S$ modulo $p$. Prenons $\epsilon$ une racine $\pa{p^h-1}$-ième de 
$\varpi^h$ et posons $y_i = x_i + \varpi^i \epsilon^{p^i} \eta^{s_i}$. Un 
calcul donne :
$$\pa{\pi_1^{n_i} \hat x_i + \pi_1^{n_i} \varpi^i 
\epsilon^{p^i}\eta^{s_i} + 
\hat c_i}^p = \pa{\pi_1^{n_i} \hat x_i + \hat c_i}^p + \pi^{n_i}
\varpi^i \epsilon^{p^{i+1}} \eta^{p s_i} + \sum_{k=1}^{p-1} \ncr p k 
\pi_1^{k n_i} \varpi^{ki} \epsilon^{k p^i}\eta^{k s_i} \pa{\pi_1^{n_i} 
\hat x_i + \hat c_i}^{p-k}$$
Or les $\hat x_i$ forment une solution de $\pa S$ et donc on a $v 
(\pi_1^{n_i} 
\hat x_i + \hat c_i) \geq \frac {er} p$. Également, on a $v \pa 
{\eta^{s_i}} = \frac {s_i}{p^h-1} \geq \frac {m_i} p$. Finalement, on 
obtient :
$$v \pa { \pi_1^{k n_i} \varpi^{ki} \epsilon^{k p^i}\eta^{k s_i} 
\pa{\pi_1^{n_i} \hat x_i + \hat c_i}^{p-k} } \geq k \frac {n_i} p + k 
\frac {m_i} p + \pa{p-k} \frac{er} p = er$$
Ainsi tous les termes de la somme sont des multiples de $p \pi^{er}$. 
Modulo $p$, il reste :
\begin{eqnarray*}
\frac 1 {\pi^{er}} \pa{\pi_1^{n_i} \hat x_i + \pi_1^{n_i} 
\varpi^i \epsilon^{p^i}\eta^{s_i} + \hat c_i}^p & \equiv & \frac 
{\pa{\pi_1^{n_i} \hat x_i + \hat c_i}^p} {\pi^{er}} + \frac {\pi^{n_i} 
\varpi^i \epsilon^{p^{i+1}} \eta^{p s_i}}{\pi^{er}} \pmod p \\
& \equiv & \varpi \hat x_{i+1} + \hat r_i + \varpi \frac {\varpi^{i+1} 
\epsilon^{p^{i+1}} \eta^{p s_i}}{\pi^{m_i }} \pmod p
\end{eqnarray*}
On remarque que $p s_i = s_{i+1} + m_i\pa{p^h-1}$, puis que $\pa{y_1, 
\ldots, y_h}$ est solution de $\pa S$.
On conclut en reprenant la démonstration du lemme \ref{lem:systnul}.
\end{preuve}

\bigskip

Voici un dernier corollaire qui nous sera utile par la suite :

\begin{cor}
\label{cor:systeme}
Soit $g$ un élément du groupe de Galois $G_K$ qui fixe tous les $\hat 
r_i$ et tous les $\hat c_i$. Soit $\pa{\hat x_1, \ldots, \hat x_h}$ une 
solution de 
$\pa S$ dans $\O_{\bar K}$. On note $x_i$ la réduction modulo $p$ de 
$\hat x_i$. Alors, pour tout $i \in \Z/h\Z$, $g$ fixe $\hat x_i$ si et 
seulement si $g$ fixe $x_i$.
\end{cor}

\begin{preuve}
Il suffit de montrer que si $g$ fixe les $x_i$ alors $\pa{g 
\hat x_1, \ldots, g \hat x_h}$ est aussi solution de $\pa S$. En 
effet, d'après le théorème précédent, si ces deux solutions sont 
distinctes dans $\O_{\bar K}$, elles le sont aussi modulo $p$. 
Le théorème \ref{th:systeme} donne ceci : dans le cas où les deux 
solutions sont distinctes, dans $\O_{\bar K}$ comme dans $\O_{\bar 
K}/p$, toutes les \og coordonnées \fg\ des $h$-uplets sont distinctes. 
Le corollaire en découle directement.
\end{preuve}

\bigskip

\noindent
{\it Remarque.} Dans le cas où tous les $\hat c_i$ sont nuls, il existe 
toujours une solution au système dans $\O_{\bar K}$. En effet, en 
combinant les équations, on aboutit à une unique équation polynomiale 
à coefficients entiers que doit vérifier $\hat x_1$. Comme $\bar K$ est 
algébriquement clos, cette équation admet une solution.

\subsection{Calcul sur les objets simples}
\label{sec:tstsimple}
Dans ce paragraphe uniquement, on suppose le corps 
résiduel $k$ algébriquement clos. On suppose également que $\pi$ est 
choisi tel que $\pi^e = p$, ce qui est toujours possible si $r > 0$ (si 
$r=0$, les résultats se démontrent indépendemment et facilement). Ainsi 
$E\pa u = u^e - p$ et $c_\pi = -1$. Soit $\calM$ un objet simple de $\Mr$. 
Le théorème \ref{th:simple} affirme que $\calM$ est de la forme $\calM 
\pa{n_i}$ pour une certaine suite périodique $\pa{n_i}$. Notons $h$ sa 
période.

\medskip

L'image de $\calM$ par le foncteur $\Tst$ s'identifie, comme 
le prouve le lemme \ref{lem:tsttilde}, à l'ensemble 
$\hom (\calM, \hat A)$. Se donner un tel morphisme revient à se 
donner pour tout $i$, un élément $x_i \in \hat A$, image de $e_i$, ces 
éléments $x_i$ étant soumis à certaines relations que nous allons 
expliciter. On rappelle que d'après le lemme \ref{lem:hata}, l'anneau 
$\hat A$ s'identifie à $\pa{\O_{\bar K} \brac X} / p$.

\begin{lemme}
\label{lem:nzero}
L'ensemble des $x \in \hat A$ tels que $N\pa x = 0$ est $\O_{\bar K}/p$.
\end{lemme}

\begin{preuve}
Le lemme résulte directement du fait que $N\pa X$ est une unité de $\hat 
A$.
\end{preuve}

\bigskip

De $N\pa {e_i} = 0$, on déduit $N \pa{x_i} = 0$ et donc d'après le lemme 
précédent $x_i \in \O_{\bar K} / p$. Intéressons-nous maintenant à la 
condition imposée par le Frobenius. Sur l'objet $\calM$, $\phi_r$ est défini 
par $\phi_r\pa{u^{n_i} e_i} = e_{i+1}$. Cela impose donc deux choses : 
l'élément $u^{n_i} x_i$ appartient à $\Fil^r \pa{\O_{\bar K} \brac 
X}/p$ et on a l'égalité $\phi_r \pa{ u^{n_i} x_i } = x_{i+1}$.

\medskip

On rappelle que l'on avait appelé $p_1$ (resp. $\pi_1$) une racine 
$p$-ième de $p$ (resp. de $\pi$) et comme $\pi^e = p$, on peut supposer
en outre que $\pi_1^e = p_1$. D'autre part, si $x = 
\sum_{i \geq 0} a_i \frac{X^i}{i!} \in \hat A$ ($a_i \in \O_{\bar K} / 
p$), alors $x \in \Fil^r \hat A$ si et seulement si $a_i$ est un 
multiple de $\bar p_1^{r-i}$ pour tout entier 
$i$ compris entre $0$ et $r$. Comme, dans $\hat A$, $u = \frac 
{\pi_1}{1+X}$, on a 
$u^{n_i} x_i \in \Fil^r \hat A$ si et seulement si $\pi_1^{n_i} x_i \in 
\Fil^r \hat A$, c'est-à-dire $\pi_1^{er-n_i}$ divise $x_i$ pour tout $i 
\in \Z/h\Z$.

\medskip

Soit $\hat x_i$ un relevé de $x_i$ dans $\O_{\bar K}$ qui est un 
multiple de $\pi_1^{er-n_i}$. Par définition $\phi_r \pa{ u^{n_i} x_i }$ 
est la réduction modulo $p$ de :
$$\frac 1 {p^r} \cdot \phi\pa{\frac{\pi_1^{n_i}}{\pa{1+X}^{n_i}} \hat 
x_i} = \frac {\pa{-1}^r} {p^r} \cdot \frac {\pi_1^{p n_i}}{\pa{1+X}^{p 
n_1}} \hat x_i ^p = \pa{-1}^r \frac 1 {\pa{1+X}^{p n_1}} \cdot \frac 1 
{\pi^{er-n_i}} \hat x_i ^p.$$
Or modulo $p$, $\pa{1+X}^p = 1$ et finalement $\phi_r\pa{ 
u^{n_i} x_i } = \pa{-1}^r \frac {x_i ^p} {\pi ^{er-n_i}}$.

Ces équations fournissent un système qui est exactement celui étudié 
dans le paragraphe \ref{sec:systeme} avec $\varpi = \pa{-1}^r$ et $\hat 
c_i = \hat r_i = 0$. En particulier, le lemme \ref{lem:systnul} nous 
fournit directement les solutions.

\medskip

On vient de prouver le théorème \ref{intro:tstsimple} dont nous 
rappelons l'énoncé :

\begin{theo}
\label{th:tstsimple}
Supposons $k$ algébriquement clos et $er < p-1$. Si l'objet simple $\calM$ 
s'identifie à $\calM\pa{n_i}$ pour une suite $\pa{n_i}$ périodique de 
période $h$ (voir théorème \ref{th:simple}), alors la représentation 
galoisienne $\Tst\pa \calM$ est isomorphe à :
$$\theta_1^{m_1} \theta_2^{m_2} \ldots \theta_h^{m_h}$$
où $m_i$ est défini par $n_i + m_i = er$ et où les $\theta_i$ sont les
caractères fondamentaux de niveau $h$.

En particulier, pour tout objet $\calM$ de $\Mr$ tué par $p$, les exposants
qui décrivent l'action de l'inertie modérée sur la semi-simplifiée
modulo $p$ de $\Tst\pa \calM$ sont tous compris entre $0$ et $er$.
\end{theo}

\subsection{Exactitude et fidélité}
\subsubsection*{Exactitude}

\begin{theo}
Le foncteur $\Tst$ de la catégorie $\Mr$ dans la catégorie des 
$\Z_p$-représentations galoisiennes de torsion est exact.
\end{theo}

\begin{preuve}
La preuve est en tout point semblable à celle donnée 
dans le paragraphe 3.2.1. de \cite{breuil-ens}, et dans le paragraphe
2.3.1. de \cite{breuil-invent}.
\end{preuve}

\subsubsection*{Fidélité}

Commençons par le lemme suivant :

\begin{lemme}
\label{lem:imsimple-clos}
Supposons $k$ algébriquement clos.
L'image par le foncteur $\Tst$ d'un objet simple de $\Mr$ est une 
représentation irréductible.
\end{lemme}

\begin{preuve}
Par le théorème \ref{th:tstsimple}, on connaît l'image d'un objet simple 
par le foncteur $\Tst$. On vérifie directement que cette image est une 
représentation galoisienne irréductible.
\end{preuve}

\begin{cor}
\label{cor:lgalgclos}
Supposons $k$ algébriquement clos.
Si $\calM$ est un objet de $\Mr$, on a :
$$\lg\pa{\calM} = \lg\pa{\Tst\pa \calM}$$
\end{cor}

\begin{preuve}
Cela découle directement du lemme précédent et de l'exactitude.
\end{preuve}

\bigskip

\noindent
{\it Remarque.} Ces deux derniers résultats restent vrais si $k$ n'est 
pas algébriquement clos (voir théorème \ref{th:longueur}).

\begin{cor}
\label{cor:fidele}
Le foncteur $\Tst$ de la catégorie $\Mr$ dans la catégorie des 
$\Z_p$-représentations galoisiennes de torsion est fidèle.
\end{cor}

\begin{preuve}
Supposons dans un premier temps $k$ algébriquement clos.
Soit $f : \calX \to \calY$ un morphisme dans la catégorie $\Mr$ tel que 
$\Tst\pa f = 0$. On a la suite exacte dans $\Mr$ :
$$\xymatrix{
0 \ar[r] & \Ker f \ar[r] & \calX \ar[r]^-{\tilde f} & \Im f \ar[r] & 0 }$$

En outre l'application $\Im f \to \calY$ est injective et donc la flèche 
déduite $\Tst\pa {\calY} \to \Tst\pa{\Im f}$ est surjective. On en déduit 
que $\Tst(\tilde f) = 0$. En appliquant le foncteur exact $\Tst$ à la 
suite exacte écrite précédemment, on voit que $\Tst\pa {\im f} = 0$.
D'après le corollaire précédent, $\im f = 0$, puis $f=0$.

\medskip

Pour le cas général, notons $K^\nr$ le complété $p$-adique de 
l'extension maximale non ramifiée de $K$. Son corps résiduel 
s'identifie à une clôture algébrique 
$\bar k$ de $k$. Désignons par $S_\nr$ l'anneau $S$ construit à 
partir de $K^\nr$ et par $\Mr_\nr$ la catégorie de modules sur $S_\nr$.

Si $\calM$ est un objet de $\Mr$, alors $\calM_\nr = S_\nr \otimes_S 
\calM$ est un objet de $\Mr_\nr$ et l'application :
$$\begin{array}{rcl}
\Tst\pa \calM & \to & \Tst(\calM_\nr) \\
f & \mapsto & \cro{s \otimes x \mapsto s f \pa x}
\end{array}$$
est un isomorphisme commutant à l'action de $\gal\pa{\bar K / K^\nr}$. 
De plus, le morphisme $\iota_M : \calM \to \calM_\nr$, $x \mapsto 1 
\otimes x$ est injectif.

Soit $f : \calX \to \calY$ un morphisme dans la catégorie $\Mr$ tel que 
$\Tst\pa f = 0$. Il induit un morphisme $f_\nr : \calX_\nr \to  
\calY_\nr$ de la catégorie $\Mr_\nr$ et on a $\Tst(f_\nr) = 0$. Par 
la fidélité dans le cas algébriquement clos, il vient $f_\nr = 0$. La 
composée $\iota_{\calY} \circ f : \calX \to \calY_\nr$ est nulle et 
comme $\iota_{\calY}$ est injectif, $f$ est nulle. Ceci démontre la 
fidélité.
\end{preuve}

\section{Pleine fidélité du foncteur $\Tst$}

Dans cette partie, on suppose à nouveau dans un premier temps que le 
corps résiduel $k$ est algébriquement clos. La propriété de pleine 
fidélité reste valable sans cette hypothèse et nous verrons dans le 
dernier paragraphe comment le cas général se déduit 
simplement du cas \og algébriquement clos \fg.

\medskip

Par un argument classique (voir \cite{fontaine-laffaille}), on se ramène 
à prouver le lemme suivant :

\renewcommand{\N}{\mathcal N}

\begin{lemme}
\label{lem:extinj}
Soient $\M$ et $\N$ deux objets simples de $\Mr$. Alors la flèche
canonique $\ext^1\pa{\M,\N} \to \ext^1\pa{\Tst\pa\N, \Tst\pa\M}$ est 
injective.
\end{lemme}

\subsection{Le module $\Ass$}
\label{sec:defass}
Pour prouver le lemme \ref{lem:extinj}, on considère $\M$ et $\N$ deux 
objets simples, $\calX$ une extension dans la catégorie $\Mr$ de 
ces deux objets telle que $\Tst\pa \calX$ soit isomorphe au produit 
direct $\Tst\pa\M \times \Tst\pa\N$. Il nous faut montrer que 
$\calX$ est isomorphe à $\M \times \N$.

\medskip

Les hypothèses impliquent que $\calX$ est tué par $p$. En effet, 
$\Tst\pa \calX$ est tué par $p$, ce qui signifie que la multiplication 
par $p$ sur $\Tst\pa \calX$ est l'application nulle. Par fidélité, on 
en déduit que la multiplication par $p$ sur $\calX$ est également 
l'application nulle. Ainsi on peut travailler dans les catégories 
$\Mrtilde$. D'autre part, si $r = 0$, il y a un unique objet
simple à isomorphisme près, ce qui règle rapidement ce cas. Ainsi on
peut supposer $r > 0$ et supposer à nouveau $\pi^e = p$.

\bigskip

Commençons par donner une caractérisation, faisant intervenir 
explicitement le foncteur $\Tst$, des objets de $\Mrtilde$ qui 
sont semi-simples.

\medskip

On construit un sous-module $\Ass$ de $\hat A$ ($\text{ss}$ pour 
\emph{semi-simple}). Pour cela, comme dans le 
paragraphe \ref{sec:systeme}, on fixe, pour tout entier 
$h$, $\eta^{\pa h}$ une racine $\pa{p^h-1}$-ième de l'uniformisante 
$\pi$. On impose en outre une condition de compatibilité : on demande 
que lorsque $h'$ divise $h$, on ait :
$$\pa{\eta^{\pa {h'}}}^{\frac{p^h-1}{p^{h'}-1}} = \eta^{\pa h}$$
De cette façon, si $s \in \Z_{\pa p}$ (le localisé de 
$\Z$ en $p$) on pourra sans ambiguité parler de $\pi^s$. En effet, comme 
tout nombre premier à $p$ admet un multiple de la forme $p^h-1$, on
peut toujours écrire $s=\frac a {p^h-1}$, et poser :
$$\pi^s = \pa{\eta^{\pa h}}^a.$$
La condition de compatibilité dit précisément que le résultat ne dépend 
pas de la fraction choisie pour représenter $s$. En outre, on a les 
formules évidentes $\pi^s \times \pi^{s'} = \pi^{s+s'}$ et $\pi^{ns} = 
\pa{\pi^s}^n$ si $s$ et $s'$ sont dans $\Z_{\pa p}$ et si $n$ est un 
entier.

\medskip

Reprenons la description donnée tout à la fin du paragraphe 
\ref{sec:simpalgclos}. Choisissons un élément $t \in \rat$ classifiant un 
certain objet simple $\M$ de la catégorie $\Mrtilde$. Appelons $t_1, 
\ldots, t_h$ les rationnels de $\Z_{\pa p} \cap \left[0,1\right[$
correspondant à $t$. Précisément si c'est la suite $\pa{n_i}_{i \in 
\Z/h\Z}$ qui classifie $\M$, on aura :
$$t_i = 0 \, , \, n_i n_{i+1} {\ldots} n_{i+h-1} \; 
\overline{n_i n_{i+1} {\ldots} n_{i+h-1}} \; \ldots$$
Si l'on pose $v_i = \frac {er}{p-1} - t_i$, on voit d'après 
le calcul fait dans le paragraphe \ref{sec:tstsimple} que tout élément 
de $\Tst\pa \M$ tombe dans le sous-$k\cro u / u^{ep}$-module de $\hat A$ 
engendré par les $\pi^{v_i}$.
On pose, pour tout $t \in \rat$ :
$$\Ass_t = k\cro u / u^{ep} \cdot \pi^{v_1} + k\cro u / u^{ep} 
\cdot \pi^{v_2} + \ldots + k\cro u / u^{ep} \cdot \pi^{v_h}$$
où l'entier $h$ dépend de $t$. La somme précédente est 
directe (voir lemme \ref{lem:direct}). Il faut faire attention au fait 
que les modules $k\cro u /  u^{ep} \cdot \pi^{v_i}$ ne sont pas libres, 
car par exemple on a toujours $u^{ep-1} \pi^{v_i} = 0$, sauf dans le cas 
très particulier où $h=1$ et $n_1 = er$. En particulier $\Ass_t$ n'est 
\emph{pas} isomorphe à $\M$.

\medskip

\begin{deftn}
On pose :
$$\Ass = \sum_{t \in \rat} \Ass_t \subset \hat A.$$
Autrement dit, $\Ass$ est le sous-$k\cro u/u^{ep}$-module engendré par 
les $\pi^{t'}$ où $t'$ parcourt l'ensemble des rationnels compris 
strictement entre $0$ et $1$ et dont l'écriture \og décimale \fg\ en 
base $p$ ne comporte que des chiffres compris entre $0$ et $er$.
\end{deftn}

\begin{lemme}
\label{lem:direct}
Le morphisme évident :
$$\bigoplus_{t'} k\cro u / u^{ep} \cdot \pi^{t'} \to \Ass$$
est un isomorphisme (où la somme est à nouveau étendue aux $t'$ 
rationnels compris strictement entre $0$ et $1$ et dont l'écriture \og 
décimale \fg\ en base $p$ ne comporte que des chiffres compris entre $0$ 
et $er$).
\end{lemme}

\noindent
Avant de faire la démonstration, insistons sur le fait 
que la notation est trompeuse : le module $k\pa u / u^{ep} \cdot 
\pi^{t'}$ n'est pas libre, il doit être vu comme un sous-module de $\hat 
A$. Le lemme dit donc que la somme dans $\hat A$ de tous ces 
sous-modules est directe.

\bigskip

\begin{preuve}
La surjectivité est une conséquence immédiate de la définition de 
$\Ass$. Passons à l'injectivité. Considérons une relation de la forme :
$$P_1 \pa u \pi^{v_1} + \ldots + P_n\pa u \pi^{v_n} = 0$$
où les $v_i$ sont deux à deux distincts et où on peut supposer que tous 
les polynômes $P_i \in k\cro u / u^{ep}$ sont non nuls. Il faut alors 
montrer que tous les termes de la somme $P_i \pa u \pi^{v_i}$ sont nuls, 
et ceci va résulter d'un simple calcul de valuation.

On écrit $u = \pi_1 X'^{p-1}$ où l'on rappelle que $X' = 1+X$ vérifie la 
relation $X'^p = 1$. En identifiant les coefficients en $X'$, on obtient 
pour tout $j$ compris entre $0$ et $p-1$ des égalités de la forme :
$$P_1^{\pa j} \pa \pi \pi^{v_1} + \ldots + P_n^{\pa j}\pa \pi \pi^{v_n} 
= 0$$
où les $P_i^{\pa j}$ sont des polynômes à coefficients dans $\O_{\bar 
K} / p$ . On rappelle que l'on 
dispose d'une valuation sur $\O_{\bar K}/p$ et que le fait d'être nul 
signifie simplement d'être de valuation supérieure à $e$. La valuation 
de $P_i^{\pa j} \pa \pi$ est un entier. Comme $v_i \in \Z_{\pa p} \cap 
[0,1[$, et que tous les $v_i$ sont deux à deux distincts, les valuations 
de $P_i^{\pa j}\pa \pi \pi^{v_i}$ sont aussi deux à deux distinctes, et 
on a :
$$v\pa{P_1^{\pa j} \pa \pi \pi^{v_1} + \ldots + P_n^{\pa j}\pa \pi 
\pi^{v_n}} = \min_i v\pa{P_i^{\pa j} \pa \pi \pi^{v_i}}$$
En particulier, la somme est nulle si et seulement si tous les termes 
sont nuls, ce qui est bien ce que l'on voulait prouver.
\end{preuve}

\bigskip

Soit $\calX$ un objet de la catégorie $\Mrtilde$. L'injection $\Ass \to \hat 
A$ fournit une flèche injective $\hom\pa{\calX, \Ass} \to \Tst \pa \calX$.

\begin{lemme}
\label{lem:caracass}
L'objet $\calX$ est semi-simple si et seulement si la flèche 
précédente est surjective (et donc un isomorphisme).
\end{lemme}

\begin{preuve}
Le sens direct est facile : si $\calX$ est semi-simple et s'écrit donc comme 
la somme $\calX = \M_1 \oplus \ldots \oplus \M_n$ pour certains objets 
simples $\M_i$, alors $\Tst\pa \calX$ se décompose lui aussi comme la somme 
directe :
$$\Tst\pa \calX = \Tst\pa {\M_1} \oplus \ldots \oplus \Tst\pa {\M_n}$$
et on a déjà vu que $\Tst \pa {\M_i} = \hom\pa{\M_i, \Ass}$.

\medskip

Faisons la réciproque. Le lemme 2.3.1.2 de \cite{breuil-invent} affirme  
que le cardinal de $\Tst\pa \calX$ est $p^{\rg \calX}$ où $\rg \calX$ désigne le 
rang de $\calX$ en tant que $k \cro u / u^{ep}$-module. On prouve par 
récurrence sur la longueur de l'objet $\calX$ que $\card \hom\pa{\calX, \Ass} 
\leq p^{\rg \calX}$ et qu'il y a égalité si et seulement si $\calX$ est 
semi-simple. Cela entraînera bien le résultat annoncé dans le lemme.

Le résultat est évident si $\calX$ est simple (de longueur $1$). Prenons 
un objet $\calX$ de longueur $n+1$. Il existe une suite exacte courte de 
la forme :
$$\xymatrix{0 \ar[r] & \M \ar[r] & \calX \ar[r] & \N \ar[r] & 0 }$$
où $\M$ est un objet simple et $\N$ est un objet de $\Mrtilde$ de 
longueur $n$. Par application du foncteur contravariant $\hom \pa 
{\cdot, \Ass}$, on en déduit une suite exacte à gauche :
$$\xymatrix{0 \ar[r] & \hom\pa{\N, \Ass} \ar[r] & \hom\pa{\calX, \Ass} 
\ar[r] & \hom\pa{\M, \Ass}}$$
d'où :
$$\card \hom\pa{\calX, \Ass} \leq \card \hom\pa{\N, \Ass} \cdot \card 
\hom\pa{\M, \Ass} \leq p^{\rg \N} \cdot p^{\rg \M} = p^{\rg \calX}$$
Pour que les deux inégalités précédentes soient des égalités, il faut 
que la flèche $\hom\pa{\calX, \Ass} \to \hom\pa{\M, \Ass}$ soit 
surjective et que $\card \hom\pa{\N, \Ass} = p^{\rg \N}$. D'après 
l'hypothèse de récurrence, cette dernière condition implique que $\N$ 
est semi-simple.

Exploitons la première condition. Soit $\psi \in \hom\pa{\M,
\Ass}$, $\psi \neq 0$. Si $t$ désigne le \og rationnel 
classifiant \fg\ de $\M$, $\psi$ tombe dans un $\Ass_t$ qui est un 
facteur direct de $\Ass$. Par hypothèse, $\psi$ se prolonge à tout 
$\calX$. On s'intéresse à la composée $s : \calX \to \Ass \to \Ass_t$ 
où la première flèche est $\psi$ ainsi prolongée et la seconde flèche 
est la projection canonique.

Notons $\pa{e_1, \ldots, e_d}$ une base adaptée de $\calX$ pour les entiers 
$n_1, \ldots, n_d$ et notons pour tout $i$, $f_i$ un relevé de 
$s\pa{e_i}$ dans $\M$, qui existe puisque tous les morphismes non nuls 
$\M \to \Ass_t$ sont surjectifs. Nous allons corriger les $f_i$ pour que 
la flèche $s : \calX \to \M, e_i \mapsto f_i$ définisse un scindage de :
$$\xymatrix {
0 \ar[r] & \M \ar[r] & \calX \ar[r] & \N \ar[r] & 0 }.$$
Les $f_i$ sont uniques modulo $u^e \Fil^r \calX$ (on peut faire beaucoup 
mieux en fait, mais ce ne sera pas utile). En particulier, quelle que 
soit la façon de les choisir, la flèche $s$ obtenue respecte $\Fil^r$. 
D'autre part, on a :
$$\pa{ \begin{array}{c}
\Phi_r\pa{u^{n_i} f_i} \\ \vdots \\ \Phi_r\pa{u^{n_d} f_d}
\end{array} } = \t G
\pa{ \begin{array}{c} f_i \\ \vdots \\ f_d \end{array}  }
+ \pa{ \begin{array}{c} r_i \\ \vdots \\ r_d \end{array} } $$
où $G$ désigne la matrice de $\Phi_r$ dans la base adaptée $\pa{e_1, 
\ldots, e_d}$ et où les $r_i$ sont des éléments de $u^e \Fil^r \calX$. On 
voit donc que si l'on remplace le vecteur $\pa {\begin{array}{c} f_i \\ 
\vdots \\ f_d \end{array} }$ par le vecteur $\pa {\begin{array}{c} f_i 
\\ \vdots \\ f_d \end{array} } + \t G^{-1} \pa{ \begin{array}{c} r_i \\ 
\vdots \\ r_d \end{array} }$, on obtient une flèche compatible à 
$\Fil^r$ et à $\phi_r$.

\medskip

Pour prouver que cette rétraction est également compatible à $N$, on 
considère le diagramme commutatif suivant :
$$\xymatrix {
& \Fil^r \M \ar[rr]^{\Phi_r} \ar[dd] & & \M \ar[dd]^{c N} \\
\Fil^r \calX \ar[rr] \ar[dd]_{u^e N} \ar[ur]^s & \ar@{ }[d]_{u^e N} \ar@{ 
}[r]_{\Phi_r} & \calX 
 \ar[dd] \ar[ur]^s \\
& \Fil^r \M \ar[rr] & \ar@{ }[r]^{\Phi_r} \ar@{ }[d]_{c N} & 
\M \\
\Fil^r \calX \ar[rr]^{\Phi_r} \ar[ur]^s & & \calX \ar[ur]^s }$$
Les faces du cube situées devant, derrière, au-dessus et au-dessous 
commutent. La face de gauche commute modulo $u^e \Fil^r \M$ et donc
$\phi_r \circ \pa{u^e N} \circ s = \phi_r \circ s \circ \pa{u^e N}$. Une 
chasse au diagramme permet d'obtenir $\pa{cN} \circ s \circ \phi_r = s 
\circ \pa{cN}\circ \phi_r$, ce qui permet de conclure puisque $\phi_r 
\pa{\Fil^r \calX}$ engendre tout $\calX$.
\end{preuve}

\subsection{Le calcul de $\hom(\N,\hat A / \Ass)$}
\label{sec:homass}
Rappelons que notre objectif est de prouver le 
lemme \ref{lem:extinj}. On considère donc $\calX$, objet de $\Mrtilde$ et 
extension de deux objets simples $\M$ et $\N$. On suppose que $\Tst\pa 
\calX \simeq \Tst\pa \M \times \Tst\pa \N$ et on veut montrer que $\calX$ 
est semi-simple. Pour cela 
d'après le lemme \ref{lem:caracass}, il suffit de prouver que tout 
élément de $\Tst \pa \calX$ définit un morphisme qui tombe dans $\Ass$. 
Soit $\psi \in \Tst\pa \calX$. On peut dessiner le diagramme suivant :
$$\xymatrix @R=15pt {
0 \ar[r] & \M \ar[r]^-{f} \ar[rdd]_-{0} & \calX \ar[r] \ar[d]^-{\psi} & 
\N \ar[r] \ar@{.>}[ldd]^-{\tilde \psi} & 0 \\ & & \hat A 
\ar[d]^-{\pr} \\ & & \hat A/\Ass }$$
La composée $\psi \circ f$ est un morphisme de $\M$ dans $\hat A$, 
qui tombe dans $\Ass$ par simplicité de $\M$ et devient 
nulle lorsqu'elle est composée avec la projection canonique. Il existe 
donc une flèche $\tilde \psi : \N \to \hat A / \Ass$ faisant 
commuter le diagramme. L'objectif de ce paragraphe est d'étudier plus en 
détail cette flèche.

\medskip

Notons d'abord que le quotient $\hat A/\Ass$ hérite d'une 
filtration, d'un Frobenius et d'un opérateur de monodromie : on définit 
$\Fil^i(\hat A/\Ass) = \pr (\Fil^i \hat A)$ et on vérifie que 
$N\pa{\Ass} \subset \Ass$ et que $\Phi_i (\Ass \cap \Fil^i \hat A) 
\subset \Ass$. Cela suffit pour transporter les structures.

\bigskip

Comme $\M$ est un objet simple, on sait le décrire précisément :
par le théorème \ref{th:simple}, il existe un entier $h$, des éléments 
$e_1, \ldots, 
e_{h}$ qui forment une base de $\M$ et des entiers $n_1, \ldots, n_{h}$  
le tout tel que $\Fil^r \M = u^{n_1} e_1 + \ldots +
u^{n_{h}} e_{h}$, $\Phi_r\pa{u^{n_{i}}e_{i}} = e_{i+1}$ et $N \pa{e_{i}} 
= 0$, les indices étant considérés dans $\Z/h\Z$. De même, on a une 
description de $\N$ : il existe un entier $h'$, des éléments $e'_1, 
\ldots, e'_{h'}$ et des entiers $n'_1, \ldots, n'_{h'}$ le tout 
vérifiant des conditions analogues.

Dans un premier temps, comme $\psi$ commute à $N$, on a
$N (\tilde \psi \pa{e'_{i'}}) = 0$ pour tout indice $i'$. On 
cherche donc les éléments de $\hat A$ dont l'image par $N$ tombe dans 
$\Ass$. C'est l'objet du lemme suivant. On rappelle que, par le
lemme \ref{lem:descy} :
$$\hat A \simeq \pa{\O _{\bar K} \cro{X'} \brac Y}/\pa{X'^p-1, p}$$
l'isomorphisme consistant à faire correspondre $X'$ à $1+X$ et $\frac 
{Y^i}{i!}$ à $\frac 1 {i!} \pa{\frac {\pa{1+X}^p-1} p}^i$.

\begin{lemme}
\label{lem:nass}
Avec les notations précédentes, l'ensemble des $x \in \hat A$ tels que 
$N\pa x \in \Ass$ est $\Ass + \O_{\bar K}/p + \pa{\Ass \cap \O_{\bar 
K}/p} Y$.
\end{lemme}

\begin{preuve}
Soit $x \in \hat A$ tel que $N\pa x \in \Ass$. Il s'écrit :
$$x = \sum_{j\geq0} P_j\pa{X'} \frac {Y^j}{j!}$$
les $P_j$ étant des polynômes de degré inférieur à $p-1$ à coefficients 
dans $\O_{\bar K}/p$ nuls pour $j \gg 0$. On a :
$$N\pa x = \sum_{j\geq 0} \pa{X' P'_j\pa{X'} + P_{j+1}\pa {X'}} \frac 
{Y^j}{j!}.$$
On remarque que \emph{via} les identifications faites, $\Ass$ est 
entièrement inclus dans $\O _{\bar K}/p \cro{X'}/(X'^p-1)$ et donc 
il suffit de vérifier les conditions :
\begin{enumerate}
\item $X' P'_0\pa{X'} + P_1\pa{X'} \in \Ass$
\item $X' P'_j \pa{X'} + P_{j+1}\pa{X'} = 0$ pour tout $j \geq 1$
\end{enumerate}
La deuxième condition entraîne $P_1\pa{X'} = b$ pour un certain $b \in 
\O_{\bar K} / p$ et $P_j\pa{X'} = 0$ pour tout $j \geq 2$.

Exploitons maintenant la première condition. Écrivons $P_0\pa{X'} = a_0 
+ a_1 X' + \ldots + a_{p-1} X'^{p-1}$ où $a_i \in \O_{\bar K} / p$. On 
obtient :
$$b + a_1 X' + 2 a_2 X'^2 + \ldots + \pa{p-1} a_{p-1} X'^{p-1} \in 
\Ass.$$
Par définition de $\Ass$ et en remarquant que $u \in \Ast$ 
correspond à $\pi_1 X'^{p-1} \in \hat A$, on voit que tous les
termes de la somme précédente sont éléments de $\Ass$.
En particulier on a $b \in \Ass$. D'autre part, les entiers 
$2, \ldots, p-1$ sont inversibles dans $\O_{\bar K} / p$ et donc tous 
les $a_i X'^i$, pour $i \geq 1$, sont aussi éléments de $\Ass$. Cela 
prouve finalement que $P_0\pa{X'} \in \Ass + \O_{\bar K} / p$ puis
la conclusion annoncée.

Il reste à faire la réciproque, mais elle est immédiate au vu du calcul 
précédent.
\end{preuve}

\subsection{Fin de la preuve}
Choisissons des relevés de $e'_{i'}$ dans $\calX$, relevés que l'on 
appelle encore $e'_{i'}$. D'après le lemme \ref{lem:nass}, l'application 
$\psi$ a la forme suivante :
\begin{eqnarray*}
\psi\pa{e_i} & = & x_i \\
\psi\pa{e'_{i'}} & = & a_{i'} + a'_{i'} + b_{i'} Y
\end{eqnarray*}
où $a_{i'} \in \O_{\bar K}/p$, $a'_{i'} \in \Ass$, $b_{i'} \in \Ass \cap 
\O_{\bar K}/p$, et où on connaît précisément la forme des $x_i$ d'après 
le calcul fait dans le paragraphe \ref{sec:tstsimple} : si $x_i \neq 0$, 
si l'on note comme dans le paragraphe \ref{sec:defass} :
\begin{eqnarray*}
t_i & = & 0 \, , \, n_i n_{i+1} {\ldots} n_{i+h-1} \;
\overline{n_i n_{i+1} {\ldots} n_{i+h-1}} \; \ldots \\
t'_{i'} & = & 0 \, , \, n'_{i'} n'_{i'+1} {\ldots} n'_{i'+h'-1} \;
\overline{n'_{i'} n'_{i'+1} {\ldots} n'_{i'+h'-1}} \; \ldots
\end{eqnarray*}
et si on pose $v_i = \frac {er}{p-1} - t_i$ et $v'_{i'} = \frac
{e\pa{r-1}}{p-1} - t'_{i'}$, il existe deux racines $\pa{p^h-1}$-ièmes
de l'unité, $\epsilon$ et $\epsilon'$, telles que $\hat x_i =
\pa{-1}^{ri} \epsilon^{p^i} \pi^{v_i}$ et où $x_i$ est la réduction modulo 
$p$ de $\hat x_i$ (resp $\hat b_{i'}$).

\medskip

De plus, en remarquant qu'il existe $z \in \calM$ tel que 
$u^{n'_{i'}} e'_{i'} + z \in \Fil^r \calX$, on obtient des relations 
de la forme :
\begin{equation}
\label{eq:phir}
\phi_r\pa{ u^{n'_{i'}} \pa{a_{i'} + b_{i'} Y} + c_{i'} } = 
a_{i'+1} + b_{i'+1} Y + r_{i'+1}
\end{equation}
où $c_{i'}$ et $r_{i'}$ sont des éléments de $\Ass$. Écrivons :
$$c_{i'} = c_{i'}^{\pa 0} + c_{i'}^{\pa 1} u + \ldots + c_{i'}^{\pa 
{p-1}} u^{p-1}$$
avec $c_{i'}^{\pa j} \in \Ass \cap \O_{\bar K} / p$. On peut 
bien arrêter la somme à $p-1$ car $u^p = \pi$ et si $x \in \O_{\bar K} / 
p$, on a bien $\pi x \in \Ass$ si et seulement si $x \in \Ass$. On peut 
également supposer $c_{n'_{i'}} = 0$ quitte à modifier $a_{i'}$.
Décomposons $u_{n'_{i'}} \pa{a_{i'} + b_{i'} Y} + c_{i'}$ de la
façon suivante :
\begin{eqnarray*}
u^{n'_{i'}} \pa{a_{i'} + b_{i'} Y} + c_{i'} & \equiv &
    \pi_1^{n'_{i'}} a_{i'} + c_{i'}^{\pa 0} + c_{i'}^{\pa 1} \pi_1 + 
    \ldots + c_{i'}^{\pa {p-1}} \pi_1^{p-1} \\
& & - \cro{n'_{i'} \pi_1 ^{n'_{i'}} a_{i'} + 
    c_{i'}^{\pa 1} \pi_1 + \ldots + \pa{p-1} c_{i'}^{\pa {p-1}} 
    \pi_1^{p-1} + \pi_1 ^{n'_{i'}} b_{i'} } X \\
& \equiv & U - V X \pmod {\frac{X^i}{i!}, i \geq 2}
\end{eqnarray*}
Cette quantité doit appartenir à $\Fil^r \hat A$. On en déduit que
$\pi_1^{er}$ divise $U$ et $\pi_1^{er-e}$ divise $V$. De plus, en 
identifiant les termes constants en $Y$ dans \ref{eq:phir}, on obtient
les relations :
$$\phi_r \pa U = a_{i'+1} + r_{i'+1}$$
qui impliquent $r_{i'} \in \Ass \cap \O_{\bar K} / p$.
Notons $\hat a_{i'} \in \O_{\bar K}$ un relevé de $a_{i'}$ et $\hat 
c_{i'}^{\pa j} \in \O_{K^\mr}$ ($K^\mr$ désigne l'extension maximale 
modérément ramifiée de $K$) des relevés de $c_{i'}^{\pa j}$ et 
posons :
$$\hat c_{i'} = \hat c_{i'}^{\pa 0} + \hat c_{i'}^{\pa 1} \pi_1 +
\ldots + \hat c_{i'}^{\pa {p-1}} \pi_1^{p-1} \in K^\mr \cro{\pi_1}.$$
Notons finalement $\hat r_{i'} \in \O_{K^\mr}$ un relevé de $r_{i'}$.
Intéressons-nous au système (dont les inconnus sont les $\hat x_{i'}$) 
donné par les équations :
$$\frac {\pa{\pi_1^{n'_{i'}} \hat x_{i'} + \hat c_{i'}}^p}{\pi^{er}} =
\pa{-1}^r \pa {\hat x_{i'+1} + \hat r_{i'}}.$$
On vient de voir que les $a_{i'}$ forment une solution modulo $p$, qui
se remonte d'après le lemme \ref{hensel} en une solution dans
$\O_{\bar K}$ que l'on note $\hat a_{i'}$. Le corollaire
\ref{cor:systeme} s'applique : un élément du groupe de Galois $G_K$
qui fixe les $\hat c_{i'}$ et les $\hat r_{i'}$, fixe $\hat a_{i'}$ si 
et seulement s'il fixe $a_{i'}$.

\bigskip

D'autre part, on rappelle que par hypothèse la suite :
$$\xymatrix{0 \ar[r] & \Tst \pa \N \ar[r] & \Tst \pa \calX \ar[r] & \Tst 
\pa \M \ar[r] & 0 }$$
est exacte et que l'on dispose d'une section $s : \Tst \pa \M \to \Tst 
\pa \calX$ qui commute à l'action de Galois.

Soit $\psi \in \Tst \pa \M$ décrit comme on vient de le voir. Le 
morphisme $s \pa \psi \in \Tst \pa \calX$ prolonge $\psi$, on 
l'appellera simplement $\psi$ par la suite. Comme $s$ est 
compatible à Galois, pour tout élément $\sigma$ stabilisant les $x_i$, 
on a :
$$\sigma\pa{a_{i'} + c_{i'} + b_{i'} Y} = \sigma a_{i'} + 
\sigma c_{i'} + \sigma b_{i'} t\pa \sigma + \sigma b_{i'} Y = 
a_{i'} + c_{i'} + b_{i'} Y$$
où on rappelle que $\sigma\pa Y = Y + t\pa\sigma$. On a vu que $t \pa 
\sigma \in \O_{\bar K}/p$ (voir lemme \ref{lem:tsigma0}) et donc de 
l'égalité précédente, on déduit en particulier :
$$\sigma a_{i'} + \sigma c_{i'} + \sigma b_{i'} t \pa \sigma = 
a_{i'} + c_{i'}.$$
Si $t\pa \sigma = 0$ et si $\sigma$ fixe $c_{i'}$, on obtient 
$\sigma a_{i'} = a_{i'}$. En particulier, tout $\sigma \in \gal\pa{\bar 
K / K^\mr\cro{\pi_1}}$ vérifie $\sigma a_{i'} = a_{i'}$. Comme de plus
tout tel $\sigma$ fixe $\hat c_{i'}$ et $\hat r_{i'}$, il vient $\sigma 
\hat a_{i'} = \hat a_{i'}$ puis $\hat a_{i'} \in K^\mr\cro{\pi_1}$.

Dès lors, la relation $\phi_r \big(\pi_1^{n'_{i'}} a_{i'} + c_{i'}\big) 
= 
a_{i'+1} + r_{i'+1}$ entraîne que $a_{i'+1} \in \O_{K^\mr}$, et ce
bien sûr pout tout $i$. Ainsi, il existe un entier $d$ tel que l'on 
puisse écrire :
$$a_{i'} = \sum_{v \in I} \lambda_v \pi^v$$
où $I$ désigne l'ensemble des rationnels dans $[0,1[$ et 
ayant pour démoninateur $\pa{p^d-1}$ et où les $\lambda_v$ sont 
des éléments de $\O_K / p$. Soit $I_{\text{ss}}$ l'ensemble des 
rationnels appartenant à $I$ dont le développement \og décimal \fg\ en 
base $p$ ne fait intervenir que des chiffres compris entre $0$ et $er$.
Soit $I_{\overline{\text{ss}}} = I \backslash I_{\text{ss}}$. On pose :
$$a_{\text{ss},i'} = \sum_{v \in I_{\text{ss}}} \lambda_v \pi^v \quad 
\text{et} \quad a_{\overline{\text{ss}},i'} = \sum_{v \in I_{\overline 
{\text{ss}}}} \lambda_v \pi^v.$$
Alors $a_{i'} = a_{\text{ss},i'} + a_{\overline{\text{ss}},i'}$, 
$a_{\text{ss},i'} \in \Ass$ et on vérifie que :
$$\Phi_r \pa{\pi_1^{n'_{i'}} a_{\overline{\text{ss}},i'} } = 
a_{\overline{\text{ss}},{i'+1}}.$$
On sait résoudre cette équation et ses solutions sont dans $\Ass$. Cela
entraîne $a_{\overline{\text{ss}},i'}=0$ pour tout indice $i'$. Ainsi 
$a_{i'} \in \Ass$.

\bigskip

Reprenons à présent l'élément $U = \pi_1^{n'_{i'}} a_{i'} + c_{i'}^{\pa 
0} + c_{i'}^{\pa 1} \pi_1 + \ldots + c_{i'}^{\pa {p-1}} \pi_1^{p-1}$.
Comme on a vu que $a_{i'} \in \Ass$, sa valuation est un élément de 
$\Z_{\pa p}$. Ainsi les valuations de termes non nuls intervenant dans 
$U$ sont deux à deux distinctes. Puisque $\pi^{er}$ divise $U$, on en 
déduit qu'il divise chacun de ces termes. En particulier, cela implique :
$$\phi_{r-1} \pa V = \phi_{r-1} \pa{ \pi^{n'_{i'}} b_{i'}}$$
et en regardant la composante sur $Y$, la relation (\ref{eq:phir}) 
implique :
$$\phi_{r-1} \pa{ \pi^{n'_{i'}} b_{i'}} = b_{i'+1}$$
On a déjà résolu plusieurs fois ce système. En particulier (c'est tout
ce dont on aura besoin), $b_{i'} \in \Ass$ et si $b_{i'} \neq 0$, on a :
$$v\pa{b_{i'}} \leq \frac {e\pa{r-1}} p$$
Supposons par l'absurde qu'il existe un indice $i'$ tel que $b_{i'} \neq 
0$. 
Soit $\sigma \in \gal\pa{\bar K/K^\mr}$ ne fixant pas $\pi_1$. On a 
$\sigma b_{i'} = b_{i'}$ et $t \pa \sigma \neq 0$. On a démontré
dans le lemme \ref{lem:tsigma} que $t \pa \sigma$ était congru à une
racine $\pa{p-1}$-ième de $\pa{-p}$. En particulier, il est de valuation 
$\frac e {p-1}$. On en déduit :
$$v\pa{b_{i'} t\pa \sigma } \leq \frac {e\pa{r-1}}{p-1} + \frac e {p-1} 
= \frac {er}{p-1} < e$$
et donc $b_{i'} t\pa \sigma$ est non nul dans $\O_{\bar K}/p$. Mais on a 
l'égalité :
$$\sigma a_{i'} + \sigma c_{i'} + \sigma b_{i'} t \pa \sigma = 
a_{i'} + c_{i'}$$
qui se simplifie ici en $b_{i'} t\pa\sigma = 0$. C'est une 
contradiction. Ainsi $b_{i'} = 0$ pour tout $i$.

\medskip

En conclusion, l'application $\psi$ prend la forme suivante :
\begin{eqnarray*}
\psi\pa{e_i} & = & x_i \\
\psi\pa{e'_{i'}} & = & a_{i'} + a'_{i'}
\end{eqnarray*}
avec $x_i$, $a_{i'}$ et $a'_{i'}$ éléments de $\Ass$. On en déduit que
$\psi$ tombe dans $\Ass$.

\bigskip

Maintenant, tout élément de $\hom(\calX, \hat A)$ s'écrit comme une 
somme d'un élément de $\hom(\calM, \hat A)$ et de l'image par $s$
d'un élément de $\hom(\calN, \hat A)$. Ainsi on a bien prouvé que 
$\hom(\calX, \hat A) = \hom(\calX, \Ass)$ et par suite que le
foncteur $\Tst$ est pleinement fidèle, du moins dans le cas où $k$ est
algébriquement clos.

\subsection{Récapitulatif et conclusion}
Récapitulons tout ce que l'on vient de voir. On a prouvé sans 
hypothèse sur le corps résiduel $k$ que le foncteur $\Tst$ est 
toujours exact et fidèle. On a également prouvé, pour l'instant, que si
ce corps résiduel est algébriquement clos, alors le foncteur $\Tst$ 
était également plein. En procédant comme dans le paragraphe 6.2 de
\cite{fontaine-laffaille}, on peut déduire le résultat pour $k$ 
quelconque du résultat pour $k$ algébriquement clos :

\begin{theo}
Le foncteur $\Tst$ de la catégorie $\Mr$ dans la catégorie des 
représentations $\Z_p$-linéaires de torsion du groupe de Galois $G_K$
est exact et pleinement fidèle.
\end{theo}

\noindent
{\it Remarque.} L'image essentielle du foncteur $\Tst$ est incluse dans 
la catégore des $\Z_p$-représenta\-tions de longueur finie de $G_K$ 
comme le montre le théorème \ref{th:longueur} que nous prouvons par la 
suite.

\bigskip

Nous pouvons finalement répondre complètement à la conjecture A.2 
formulée dans \cite{breuil-invent}. Mais avant cela, nous allons
énoncer et prouver une propriété formelle :

\begin{prp}
Soient $A$ et $B$ deux catégories abéliennes et artiniennes. Soit
$F : A \to B$ un foncteur additif, exact et pleinement fidèle qui
est tel que l'image de tout objet simple de $A$ est encore simple
dans $B$. Alors l'image essentielle de $F$ est stable par sous-objets
et par quotients.
\end{prp}

\begin{preuve}
On se ramène directement au cas où $A$ est une sous-catégorie pleine
de $B$. L'hypothèse dit que les objets simples de $A$ restent simples
dans $B$. En particulier si $M$ est un objet de $A$ et si :
$$0 = M_0 \subset M_1 \subset \ldots \subset M_m = M$$
est une suite de Jordan-Hölder dans $A$, elle restera une suite de 
Jordan-Hölder dans $B$. Il s'agit de prouver que la catégorie $A$ est
stable par sous-objets et par quotients.

Introduisons pour cela $A'$ la sous-catégorie pleine de $A$ formée
des objets dont tous les quotients de Jordan-Hölder sont dans $B$. C'est
une sous-catégorie abélienne de $A$ qui est stable par sous-objets et 
par quotients. \'Evidemment $A$ est une sous-catégorie de $A'$, on peut
donc supposer que $A' = B$ ou si l'on préfère que les objets simples de 
$A'$ et ceux de $B$ sont les mêmes.

Soit $M$ un objet de $A$ et $N$ un sous-objet de $M$. En
considérant des suites de Jordan-Hölder de $N$ et de $M/N$, on voit que
l'on peut écrire une suite de Jordan-Hölder de la forme suivante :
$$0 = M_0 \subset M_1 \ldots \subset M_n = N \subset M_{n+1} \subset
\ldots \subset M_m = M$$
Le quotient $M_m/M_{m-1}$ est un objet simple et donc un objet de $A$.
Par suite le noyau de la projection $M_m \to M_m/M_{m-1}$ qui 
s'identifie à $M_{m-1}$ est également objet de $A$. Par récurrence, on 
montre que tous les $M_i$ sont objets de $A$ et donc qu'il en est 
de même de $N$. Ceci prouve la stabilité par sous-objets, la stabilité 
par quotients se traite de façon totalement identique.
\end{preuve}

\bigskip

\noindent {\it Remarque.}
Cette propriété redémontre en particulier le fait que la 
sous-catégorie $\MFrtilde$ de $\Mrtilde$ est stable par sous-objets et 
par quotients, puisque l'on a vu dans la proposition \ref{prop:sobj} 
que tous les objets simples de $\Mrtilde$ étaient dans $\MFrtilde$.

\medskip

On peut désormais énoncer le théorème qui résout la conjecture
mentionnée précédemment :

\begin{theo}
\label{theo:tststable}
L'image essentielle du foncteur $\Tst$ est stable par sous-objets et par 
quotients et indépendante du choix de l'uniformisante $\pi$.
\end{theo}

\begin{preuve}
L'indépendance du choix de l'uniformisante est une conséquence directe 
de la propriété \ref{prop:tstunif}.

\medskip

Supposons $k$ algébriquement clos. On sait, par le lemme 
\ref{lem:imsimple-clos}, que l'image par le foncteur $\Tst$ d'un objet
simple de $\Mr$ est une représentation irréductible. Le foncteur $\Tst$ 
vérifie les conditions de la propriété précédente, ce qui conclut.

Pour le cas général, notons $K^\nr$ le complété $p$-adique de l'extension 
maximale non ramifiée de $K$ et $\Mrtilde_\nr$ la catégorie $\Mrtilde$ 
construite à partir de $K^\nr$.
Soit $\calM$ un objet simple de $\Mr$. Il est tué par $p$ et donc peut 
être vu comme un objet de $\Mrtilde$. Il suffit de prouver que $\Tst \pa 
\calM$ est une représentation irréductible. Notons $\calM_\nr = \bar k 
\otimes_k \calM$. L'application :
$$\begin{array}{rcl}
\Tst\pa \calM & \to & \Tst(\calM_\nr) \\
f & \mapsto & \cro{\lambda \otimes x \mapsto \cro\lambda f \pa x}
\end{array}$$
(où $\cro \lambda \in W(\bar k) \subset \O_{\bar K}$ est le représentant 
de Teichmüller de $\lambda \in \bar k$) est un isomorphisme commutant à 
l'action de $G_{K^\nr} = \gal\pa{\bar K / K^\nr}$.

Supposons par l'absurde qu'il existe $V$ un sous-$\Z_p$-module de 
$\Tst\pa \calM$, strict, non nul et $G_K$-équivariant. C'est 
aussi un sous-$\Z_p$-module de $\Tst(\calM_\nr)$ $G_{K^\nr}$-équivariant 
et donc d'après le cas précédent, on peut écrire $V = \Tst(\calC_\nr)$ 
(égalité de représentations de $G_{K^\nr}$) où $\calC_\nr$ est un 
quotient de $\calM_\nr$ dans la catégorie $\Mr_\nr$. Soit $\calM'_\nr$ 
le noyau de la projection $\calM_\nr \to \calC_\nr$, c'est un 
sous-objet strict et non nul de $\calM_\nr$ dans la catégorie 
$\Mrtilde_\nr$.

Soient $\sigma \in \gal\pa{K^{\nr}/K}$ et $\hat \sigma \in G_K$ un 
prolongement de $\sigma$. Soient $\psi \in V \subset \Tst\pa \calM$ et 
$\psi_\nr$ son image dans $\Tst(\calM_\nr)$. L'élément $\hat \sigma$ 
agit sur $\psi_\nr$ de la façon suivante :
$$\begin{array}{rcl}
\hat \sigma \psi_\nr : \quad \calM_\nr & \to & \hat A \\
s \otimes x & \mapsto & s \, \hat \sigma \psi \pa x
\end{array}$$
De plus, $\sigma$ définit un morphisme $\sigma : \calM_\nr \to 
\calM_\nr$ dans la catégorie $\Mrtilde_{K^\nr}$. On vérifie que le 
diagramme suivant commute :
$$\xymatrix @C=50pt {
\calM_\nr \ar[r]^-{\sigma} \ar[d]_-{\hat \sigma^{-1} \psi_\nr} & 
\calM_\nr \ar[d]^-{\psi_\nr} \\
\xyhat A \ar[r]_-{\hat \sigma} & \xyhat A }$$
Comme $\psi \in V$, on a $\psi_{|\calM'_\nr} = 0$ et par le diagramme 
précédent, $\psi_{|\sigma \calM'_\nr} = 0$.

On obtient un diagramme de la forme :
$$\xymatrix @C=50pt {
0 \ar[r] & V \ar[r] \ar@{=}[d] & \Tst(\calM_\nr) \ar[r] \ar@{=}[d] & 
\Tst(\calM'_\nr) \ar[r] & 0 \\
0 \ar[r] & V \ar[r] & \Tst(\calM_\nr) \ar[r] & \Tst(\sigma \calM'_\nr) 
\ar[r] & 0 \\
}$$
qui fournit un isomorphisme $\Tst(\calM'_\nr) \to \Tst(\sigma
\calM'_\nr)$, se relevant par pleine fidélité en un isomorphisme $\sigma 
\calM'_\nr \to \calM'_\nr$ faisant commuter le diagramme suivant :
$$\xymatrix @C=50pt {
0 \ar[r] & \sigma \calM'_\nr \ar[d] \ar[r] & \calM \ar@{=}[d] \\
0 \ar[r] & \calM'_\nr \ar[r] & \calM } $$
On en déduit $\sigma \calM'_\nr = \calM'_\nr$ pour tout $\sigma \in 
\gal\pa{K^{\nr}/K}$.

\medskip

On pose $\calM' = \calM'_\nr \cap \calM = \bar 
\calM'_\nr{}^{\gal\pa{\bar K^\nr/K}}$. On va montrer que $\calM'$ est un 
sous-objet strict et non nul de $\calM$ dans la catégorie $\Mrtilde$, ce 
qui est une contradiction. Soit $\pa{e_1, \ldots, e_d}$ une $k \cro u / 
u^{ep}$-base de $\calM$. Soit $y \in \bar \calM'_\nr$, $y \neq 0$. 
On peut écrire :
$$y = P_1\pa u e_1 + \ldots + P_d \pa u e_d$$
où les $P_i$ sont des polynômes à coefficients dans $\ell \cro u / 
u^{ep}$ pour $\ell$ une extension finie de $k$. D'autre part, si $P \in 
\ell \cro u / u^{ep}$, on peut définir $\tr_{\ell/k} \pa P$ en calculant 
la trace de chacun des coefficients. En outre, comme $\ell/k$ est 
séparable, on peut supposer $\tr_{\ell/k} \pa {P_1} \neq 0$, quitte à 
multiplier $\bar x$ par un élément non nul de $\ell$. Posons :
$$x = \tr_{\ell/k} \pa {P_1\pa u} e_1 + \ldots + \tr_{\ell/k} \pa {P_d 
\pa u} e_d.$$
C'est un élément de $\calM$ et, puisque $\sigma \calM'_\nr = 
\calM'_\nr$ pour tout $\sigma \in \gal\pa{K^{\nr}/K}$, c'est aussi un 
élément de $\calM'_\nr$. Comme on a supposé $\tr_{\ell/k} \pa {P_1\pa u} 
\neq 0$, on a $x \neq 0$, puis $\calM' \neq 0$ comme on voulait.

On pose $\Fil^r \calM' = \calM' \cap \Fil^r \calM'_\nr$. L'opérateur 
$\Phi_r : \Fil^r \calM'_\nr \to \calM'_\nr$ (resp. $N : \calM'_\nr \to 
\calM'_\nr$) induit une application $\phi_r : \Fil^r \calM' \to \calM'$ 
(resp. $N : \calM' \to \calM'$). Ces applications vérifient les bonnes 
conditions pour définir un objet de $\Mrtilde$. Le seul point délicat 
est le fait que $\Phi_r \pa{\Fil^r \calM'}$ engendre $\calM'$ en tant 
que $k\cro u / u^{ep}$-module. Soit $x \in \calM'$. On sait qu'il existe 
$\lambda_i \in \bar k \cro u / u^{ep}$ et $y_i \in \Fil^r 
\calM'_\nr$ tels que :
$$x = \lambda_1 \phi_r (y_1) + \ldots + \lambda_n \phi_r (y_n).$$
De plus, quitte à rentrer les constantes à l'intérieur des $\phi_r$, on 
peut supposer que $\lambda_i = u^{s_i}$ pour certains entiers $s_i$. 
Soit $\pa{e_1, \ldots, e_d}$ une $k\cro u / u^{ep}$-base de $\calM$. 
Écrivons :
$$y_j = P_{1,j} \pa u e_1 + \ldots + P_{d,j} \pa u e_d$$
où $P_{i,j} \in \bar k \cro u / u^{ep}$. Soit $\ell$ une extension 
finie de $k$ contenant tous les coefficients des polynômes $P_{i,j}$ 
définis ci-dessus. Comme précédemment, on peut définir $\tr_{\ell/k} \pa 
P$ pour $P \in \ell\cro u / u^{ep}$. Soit $\alpha \in \ell$ un élément 
tel que $\tr_{\ell/k}\pa \alpha = 1$. On pose :
$$x_j = \tr_{\ell/k} \pa{\alpha P_{1,j} \pa u} e_1 + \ldots + 
\tr_{\ell/k} \pa {\alpha P_{d,j} \pa u} e_d.$$
On a alors $x_i \in \Fil^r \calM'$ et :
$$x = u^{s_1} \phi_r \pa{x_1} +  \ldots + u^{s_n} \phi_r \pa{x_n}$$
ce qui prouve bien que $\Phi_r\pa{\Fil^r \calM'}$ engendre $\calM'$ en tant que 
$S$-module.
\end{preuve}

\bigskip

\noindent
{\it Remarque.} Comme conséquence du théorème précédent et de la pleine 
fidélité de $\Tst$, un objet $\calM \in \Mr$ est semi-simple si et 
seulement si $\Tst\pa\calM$ est une représentation semi-simple.

\bigskip

Il résulte de la démonstration précédente et de l'exactitude du foncteur 
$\Tst$ le théorème suivant :

\begin{theo}
\label{th:longueur}
Si $\calM$ est un objet de $\Mr$, on a :
$$\lg\pa{\calM} = \lg \pa{\Tst\pa \calM}$$
\end{theo}

\begin{prop}
Soit $\calM$ un objet de $\Mr$ isomorphe en tant que $S$-module à 
$S/p^{n_1} 
S \oplus \ldots \oplus S/p^{n_d}S$ pour certains entiers $n_i$. Alors 
en 
tant que $\Z_p$-module, $\Tst\pa \calM$ est isomorphe à $\Z_p/p^{n_1} \Z_p 
\oplus \ldots \oplus \Z_p/p^{n_d} \Z_p$.
\end{prop}

\begin{preuve}
Le lemme 2.3.1.2 de \cite{breuil-invent} dit que si $\calM$ est un objet 
de $\Mrtilde$, alors $\Tst\pa \calM$ est un $\F_p$-espace 
vectoriel de dimension $\rg X$. On en déduit par exactitude du foncteur 
$\Tst$ que :
$$\lg_S\pa{\calM} = \lg_{\Z_p} \pa{\Tst\pa \calM}$$
où les longueurs sont calculées respectivement dans la catégorie des 
$S$-modules et dans celle des $\Z_p$-modules.

Soit $\calM$ un objet de $\Mr$ isomorphe en tant que $S$-module à
$S/p^{n_1} S \oplus \ldots \oplus S/p^{n_d}S$. La représentation 
galoisienne $\Tst\pa \calM$ est un $\Z_p$-module de longueur finie et 
donc est isomorphe en tant que $\Z_p$-modules à $\Z_p/p^{n'_1} \Z_p 
\oplus \ldots \oplus \Z_p/p^{n'_d} \Z_p$ pour certains entiers $n'_i$.
Soit $n$ un entier. Le noyau de la multiplication par $p^n$ sur 
$\calM$ s'envoie par le foncteur exact $\Tst$ sur le conoyau de la 
multiplication par $p^n$ sur $\Tst\pa\calM$. On en déduit en regardant 
les longueurs que :
$$\sum_{i=1}^d \min \pa{n_i, n} = \sum_{i=1}^{d'} \min\pa{n'_i, n}.$$
Cela permet de conclure.
\end{preuve}

\section{Conséquences}

\subsection{Modules filtrés et modules fortement divisibles}
\label{sec:fortdiv}
\subsubsection*{Définitions}

On reprend dans ce paragraphe les définitions et propriétés du 
paragraphe 4.1.1 de \cite{breuil-ens}.

\bigskip

On rappelle que $K_0$ désigne le corps des fractions de $W$, anneau des 
vecteurs de Witt à coefficients dans $k$. On définit $S_{K_0} = S 
\otimes_W {K_0}$. C'est l'ensemble suivant :
$$S_{K_0} = \acco{\sum_{i=0}^\infty w_i \frac {\pa{E\pa u}^i}{i!}, \,
w_i \in K_0\cro u, \, \lim_{i\to \infty} w_i = 0}.$$
On munit $S_{K_0}$ d'une filtration en posant $\Fil^n S_{K_0} = \Fil^n S 
\otimes_W {K_0}$, ou encore :
$$\Fil^n S_{K_0} = \acco{\sum_{i=n}^\infty w_i \frac {\pa{E\pa 
u}^i}{i!}, \, w_i \in K_0\cro u, \, \lim_{i\to \infty} w_i = 0}.$$
On prolonge de manière évidente le Frobenius et l'opérateur de 
monodromie définis sur $S$ à tout $S_{K_0}$.

\medskip

On définit un \emph{module fortement divisible} (resp. un \emph{module
filtré sur $S_{K_0}$}) comme la donnée suivante :
\begin{enumerate}
\item un $S$-module (resp. un $S_{K_0}$-module) $\calM$ libre de rang fini ;
\item un sous-$S$-module (resp. un sous-$S_{K_0}$-module) de $\calM$, noté 
$\Fil^r \calM$ contenant $\Fil^r S \cdot \calM$ (resp. contenant $\Fil^r 
S_{K_0} \cdot \calM$) et
tel que $\calM/\Fil^r \calM$ soit sans $p$-torsion (cette dernière condition est 
automatique pour les modules filtrés sur $S_{K_0}$) ;
\item d'une flèche $\phi$-semi-linéaire $\phi_r : \Fil^r \calM \to \calM$ 
vérifiant la condition :
$$\phi_r\pa {sx} = \frac 1 {c^r} \phi_r \pa{s} \phi_r\pa{\pa{E\pa u}^r 
x}$$
et ce pour tout élément $s \in \Fil^r S$ (resp. tout élément $s \in 
\Fil^r S_{K_0}$) et tout élément $x \in \calM$ telle que $\im \phi_r$ 
engendre $\calM$ en tant que $S$-module (resp. en tant que 
$S_{K_0}$-module) ;
\item une application $W$-linéaire (resp. une application 
$K_0$-linéaire) $N : \calM \to \calM$ vérifiant les trois conditions :
\begin{itemize}
\item pour tout $s \in S$ (resp. pour tout $s \in S_{K_0}$) et tout $x 
\in \M$, $N\pa{sx} = N\pa s x + s N\pa x$
\item $E\pa u N\pa{\Fil^r \M} \subset \Fil^r \M$
\item le diagramme suivant commute :
$$\xymatrix @C=50pt {
\Fil^r \calM \ar[r]^{\phi_r} \ar[d]_{E\pa u N} & \calM \ar[d]^{cN} \\
\Fil^r \calM \ar[r]^{\phi_r} & \calM}$$
\end{itemize}
\end{enumerate}

\medskip

Suivant toujours \cite{breuil-ens}, on définit de manière évidente 
la catégorie des modules filtrés sur $S_{K_0}$ et celle des modules 
fortement divisibles. Elles sont équipées d'un foncteur vers les 
représentations galoisiennes. Précisément, si $\calM$ est un module filtré 
sur $S_{K_0}$, on pose $\Tst\pa \calM = \hom(\calM, \Bst)$ où par définition 
$\Bst = \Ast \otimes_W {K_0}$ muni des structures induites et où $\hom$ 
est compatible à toutes les structures ; on obtient une 
$\Q_p$-représentation de $G_K$. De même si $\calM$ est un module fortement 
divisible, on définit $\Tst\pa \calM = \hom(\calM, \Ast)$, le $\hom$ étant 
encore compatible à toutes les structures. On obtient une 
$\Z_p$-représentation libre de $G_K$. Les rangs des représentations 
obtenues coïncident avec les rangs des objets $\calM$.

\medskip

Si $\calM$ est un module fortement divisible, on vérifie immédiatement que 
$\calM \otimes_W K_0$ est un module filtré sur $S_{K_0}$ et que pour tout 
entier $n \geq 1$, $\calM/p^n \calM$ est un objet de la catégorie $\Mr$. De 
plus, on montre que $\calM$ s'identifie à la limite projective de $\calM/p^n \calM$, 
puis que $\Tst\pa \calM$ s'identifie à la limite 
projective de $\Tst\pa{\calM/p^n \calM}$. On déduit de la pleine 
fidélité prouvée précédemment le corollaire suivant :

\begin{theo}
\label{th:fortdivpleinfid}
Le foncteur $\Tst$ de la catégorie des modules fortement divisibles dans 
la catégories des $\Z_p$-représentations (libres) de $G_K$ est 
pleinement fidèle.
\end{theo}

\subsection{Modules fortement divisibles et foncteur $\Tst$}
Nous démontrons dans ce paragraphe le théorème \ref{theo:fortdivtst}. En 
fait, nous démontrons la formulation équivalente mais légèrement 
différente suivante :

\begin{theo}
\label{theo:fortdivtstlem}
On suppose $er < p-1$.
Soit $\calM$ un module fortement divisible sur $S$, et soit $V$ la 
représentation galoisienne associée \emph{via} le foncteur $\Tst$
à $\calM_{K_0} = \calM \otimes_W K_0$ qui est un module filtré sur $S_{K_0}$. Le 
foncteur $\Tst$ réalise une anti-équivalence de catégories entre la 
catégorie des sous-modules fortement divisibles de $\calM_{K_0}$ et celle 
des sous-$\Z_p$-réseaux de $V$ stables par $G_K$.
\end{theo}

\begin{preuve}
Nous suivons pas à pas la preuve de la proposition 3 de 
\cite{breuil-bullsoc}, qui n'utilise essentiellement que la pleine
fidélité du foncteur $\Tst$ et un équivalent du théorème 
\ref{theo:tststable}.

Dans un premier temps, la pleine fidélité du foncteur $\Tst$ considérée 
dans l'énoncé du théorème se déduit directement du théorème 
\ref{th:fortdivpleinfid}.

\medskip

Reste l'essentielle surjectivité. Soit $T$ un $\Z_p$-réseau de $V$ 
stable par $G_K$. Il existe un entier $n_0$ tel que :
$$p^{n_0} T \subset \Tst\pa \calM \subset \pa{1/p^{n_0}} T$$
On en déduit que pour $n \geq n_0$, $p^{n_0} T / p^n T$ est un 
sous-objet de $\Tst\pa \calM / p^n T$, ce dernier étant un quotient de
$\Tst\pa{\calM / p^{n+n_0} \calM}$. Le théorème \ref{theo:tststable} assure 
alors que $p^{n_0} T / p^n T$ s'écrit $\Tst\pa{\calM_n}$ pour $\calM_n$ un 
certain objet de $\Mr$.

La pleine fidélité de $\Tst$ assure l'existence d'une unique flèche 
$\calM_n \to \calM_{n+1}$ relevant la projection $p^{n_0} T / p^{n+1} T \to 
p^{n_0} T / p^n T$, et la limite inductive de ce système s'identifie 
à $\calM_\infty \otimes_{\Z_p} \Q_p/\Z_p$ pour un certain module fortement 
divisible $\calM_\infty$ qui répond à la question.
\end{preuve}

\bigskip

\noindent
{\it Remarque.} Noter que si $\calM$ est un $S_{K_0}$-module filtré \og 
faiblement admissible \fg, alors il contient toujours un module 
fortement divisible par \cite{breuil-invent}.

\subsection{Variante d'une conjecture de Serre}
\label{sec:serre}
Dans ce paragraphe, on se propose d'expliquer comment le théorème donné
dans l'introduction et que nous rappelons ci-dessous est conséquence de
la théorie présentée précédemment. Noter qu'ici, on ne suppose \emph{a 
priori} plus rien ni sur $e$, ni sur $r$.

\begin{theo}
Soit $X$ un schéma propre et lisse sur $K$ et à réduction semi-stable
sur l'anneau des entiers $\O_K$. On fixe $r$ un entier. Les exposants
qui décrivent l'action de l'inertie modérée sur la semi-simplifiée
modulo $p$ de $H^r_\et\pa{X_{\bar K}, \Q_p}^\star$ (où $X_{\bar K}$ est
l'extension des scalaires de $X$ à $\bar K$ et où \og $\star$ \fg\
signifie que l'on prend le dual) sont tous compris entre $0$ et $er$.
\end{theo}

\begin{preuve}
Dans un premier temps, il est clair que l'on peut supposer $er < p-1$, 
le théorème étant trivialement vérifié dans le cas contraire. On peut 
donc utiliser les résultats précédents.

\medskip

D'après les résultats de \cite{tsuji} et du paragraphe 2.2 de 
\cite{breuil-azumino}, la $\Q_p$-représentation $V = H^r_\et\pa{X_{\bar 
K}, \Q_p}^\star$ (le dual étant cette fois-ci le $\Q_p$-dual) provient 
\emph{via} le foncteur $\Tst$ d'un module filtré $\calM_{K_0}$ sur 
$S_{K_0}$, et d'après les résultats de \cite{breuil-invent}, ce module 
admet un sous-module fortement divisible $\calM$.

\medskip

D'autre part, la $\Z_p$-représentation $T$ est un réseau de $V$ stable 
par Galois, et donc d'après le théorème \ref{theo:fortdivtstlem}, il 
existe un module fortement divisible inclus dans $\calM_{K_0}$ dont l'image 
par $\Tst$ est isomorphe à $T$. Appelons $\calM$ un tel module.

\medskip

La représentation quotient $T/p$ correspond \emph{via} le foncteur 
$\Tst$ à $\calM/p$ qui est un objet de $\Mrtilde$. La semi-simplifiée de 
$T/p$ est la somme directe de ses quotients de Jordan-Hölder, et chacun 
de ces quotients correspond à un objet simple de $\Mrtilde$. Le théorème
\ref{th:tstsimple} permet de conclure.
\end{preuve}

\bigskip

\noindent
{\it Remarque.} Si l'on préfère, on peut ne pas utiliser le théorème
\ref{theo:fortdivtstlem}, mais dire à la place que si $T$ et $T'$ sont
deux $\Z_p$-réseaux de $V$ stables par Galois, alors les 
semi-simplifiées des réductions modulo $p$ de ces deux représentations
sont isomorphes. On aurait donc pu garder le premier module fortement
divisible $\calM$.

\nocite{fontaine-ast1} \nocite{fontaine-ast2}
\nocite{faltings} \nocite{faltings2}
\bibliography{fl3}
\bibliographystyle{amsalpha}

\end{document}